\documentclass[12pt]{article}
\usepackage{amssymb}
\usepackage{amsfonts}
\usepackage{amsmath}
\usepackage{eurosym}
\usepackage{makeidx}
\usepackage{amssymb}
\usepackage{times}
\usepackage{anysize}
\usepackage{cite}
\usepackage[polutonikogreek,english]{babel}
\usepackage[utf8x]{inputenx}

\marginsize{2cm}{2cm}{2.5cm}{2.5cm}
\newtheorem{theorem}{Theorem}[section]

\newtheorem{condition}[theorem]{Condition}

\newtheorem{corollary}[theorem]{Corollary}

\newtheorem{definition}[theorem]{Definition}

\newtheorem{lemma}[theorem]{Lemma}

\newtheorem{proposition}[theorem]{Proposition}

\newenvironment{proof}[1][Proof]{\noindent\textbf{#1.} }{\ \rule{0.5em}{0.5em}}
\input{tcilatex}
\begin{document}

\title{Weak$^{\ast }$ Hypertopologies with Application to Genericity of
Convex Sets }
\author{J.-B. Bru \and W. de Siqueira Pedra}
\date{\today }
\maketitle

\begin{abstract}
We propose a new class of hypertopologies, called here weak$^{\ast }$
hypertopologies, on the dual space $\mathcal{X}^{\ast }$ of a real or
complex topological vector space $\mathcal{X}$. The most well-studied and
well-known hypertopology is the one associated with the Hausdorff metric for
closed sets in a complete metric space. Therefore, we study in detail its
corresponding weak$^{\ast }$ hypertopology, constructed from the Hausdorff
distance on the field (i.e. $\mathbb{R}$ or $\mathbb{C}$) of the vector
space $\mathcal{X}$ and named here the weak$^{\ast }$-Hausdorff
hypertopology. It has not been considered so far and we show that it can
have very interesting mathematical connections with other mathematical
fields, in particular with mathematical logics. We explicitly demonstrate
that weak$^{\ast }$ hypertopologies are very useful and natural structures\
by using again the weak$^{\ast }$-Hausdorff hypertopology in order to study
generic convex weak$^{\ast }$-compact sets in great generality. We show that
convex weak$^{\ast }$-compact sets have generically weak$^{\ast }$-dense set
of extreme points in infinite dimensions. An extension of the well-known
Straszewicz theorem to Gateaux-differentiability (non necessarily Banach)
spaces is also proven in the scope of this application.\bigskip\

 \bigskip

\end{abstract}

\tableofcontents%

\section{Introduction}

Topologies for sets of closed subsets of topological spaces have been
studied since the beginning of the last century. When such topologies,
restricted to singletons, coincide with the original topology of the
underlying space, we talk about hypertopologies and hyperspaces of closed
sets. In the literature, there exist various hypertopologies and related
set-convergence notions, like the Fell, Vietoris, Wijsman, proximal,
Hausdorff metric, or locally finite hypertopologies, to name a few
well-known examples. For a review on this field, see, e.g., \cite{Beer}.
Hypertopologies are important topological structures in applied mathematics,
for instance to study stability of a minimization problem, as explained in
\cite[Chapter 8]{Lucchetti}. In this paper, we propose a new class of
hypertopologies whose underlying topological space is the dual space $%
\mathcal{X}^{\ast }$ of a real or complex topological vector space $\mathcal{%
X}$, endowed with the weak$^{\ast }$ topology.

The weak$^{\ast }$ topology of $\mathcal{X}^{\ast }$ is generally not
metrizable, but locally convex and Hausdorff. Denote by $\mathbf{F}(\mathcal{%
X}^{\ast })$ the set of all weak$^{\ast }$-closed subsets of $\mathcal{X}%
^{\ast }$. The hypertopologies on $\mathbf{F}(\mathcal{X}^{\ast })$ we
propose here are all constructed from the following scheme: Any element $%
A\in \mathcal{X}$ defines a weak$^{\ast }$-continuous linear mapping $\hat{A}
$ from $\mathcal{X}^{\ast }$ to the field $\mathbb{K}=\mathbb{R},\mathbb{C}$
of the vector space $\mathcal{X}$, in turn implying a canonical mapping $%
\mathbf{A}$ from $\mathbf{F}(\mathcal{X}^{\ast })$ to the set $\mathbf{F}(%
\mathbb{K})$ of all closed subsets of $\mathbb{K}$. Then, providing $\mathbf{%
F}(\mathbb{K})$ with some hypertopology, a \emph{weak}$^{\ast }$
hypertopology on $\mathbf{F}(\mathcal{X}^{\ast })$ is the so-called initial
topology of the family $\{\mathbf{A}\}_{A\in \mathcal{X}}$. This
construction is analogous to the one of the weak$^{\ast }$ topology, which
is the initial topology of the family $\{\hat{A}\}_{A\in \mathcal{X}}$. This
new class of hypertopologies, called here weak$^{\ast }$ hypertopologies,
does not seem to have been systematically considered in the past. The scalar
topology described in \cite[Section 4.3]{Beer}, adapted for a locally convex
space like $\mathcal{X}^{\ast }$, is retrospectively a first example of a
weak$^{\ast }$ hypertopology.

The most well-studied and well-known hypertopology on a complete metric
space is the Hausdorff metric topology. Therefore, we study in detail the
weak$^{\ast }$ hypertopology associated with the Hausdorff metric in the
field $\mathbb{K}=\mathbb{R},\mathbb{C}$ of the vector space $\mathcal{X}$.
It is named here the \emph{weak}$^{\ast }$-\emph{Hausdorff} hypertopology.
To our knowledge, it has not been considered so far and we show that it can
have very interesting mathematical connections with other mathematical
fields, like mathematical logics:

\begin{itemize}
\item We use the weak$^{\ast }$-closed convex hull operator $\overline{%
\mathrm{co}}$ in order to study the weak$^{\ast }$-Hausdorff topology of
hyperspaces. Such a closure operator has certainly been considered in the
past and combines (i) an algebraic (or finitary) closure operator and (ii) a
topological (or Kuratowski) closure operator. It gives a first connection
with mathematical logics where fascinating applications of closure operators
have been developed, already by Tarski himself during the 1930's.

\item We introduce the notion of immeasurable hyperspaces, by using the sets
of topologically bounded weak$^{\ast }$-closed subsets of $\mathcal{X}^{\ast
}$. It is associated with an infinite collection of weak$^{\ast }$%
-Hausdorff-clopen sets, which can be used to form a Boolean algebra, as is
usual in mathematical logics\footnote{%
See Stone's representation theorem for Boolean algebras.}. Such a study has
for instance been performed in \cite{Bankston} for the hyperspace associated
with a Stonean compact topological space and the Vietoris hypertopology.
\end{itemize}

We demonstrate that the weak$^{\ast }$-Hausdorff hypertopology cannot
distinguish a set from its weak$^{\ast }$-closed convex hull, as it appears
for other well-established hypertopologies like the scalar topology (see
\cite[Section 4.3]{Beer}). Similar to the scalar topology, only the
hyperspace $\overline{\mathrm{co}}(\mathbf{F}(\mathcal{X}^{\ast }))$ of all $%
\overline{\mathrm{co}}$-closed, or equivalently convex weak$^{\ast }$%
-closed, subsets of $\mathcal{X}^{\ast }$ is in general a Hausdorff space.
We also relate the weak$^{\ast }$-Hausdorff convergence with the celebrated
lower and upper limits of sets \`{a} la Painlev\'{e} \cite[\S\ 29]%
{topology-painleve} and strongly strengthen our results when the topological
vector space $\mathcal{X}$ is separable. This is performed by showing that
the Banach-Alaoglu theorem and the metrizability of the weak$^{\ast }$
topology for absolute polars in the dual space $\mathcal{X}^{\ast }$ of a
separable space $\mathcal{X}$ are inherited by $\mathbf{F}(\mathcal{X}^{\ast
})$. This leads to a large class of metrizable weak$^{\ast }$%
-Hausdorff-compact hyperspaces. Other aspects of the weak$^{\ast }$%
-Hausdorff hypertopology are also studied in detail.

We explicitly demonstrate that weak$^{\ast }$ hypertopologies are very
useful and natural structures\ by using again the weak$^{\ast }$-Hausdorff
hypertopology in order to study generic convex weak$^{\ast }$-compact sets
in great generality. In fact, the present paper has been originally inspired
by the fairly complicated geometrical structure of the convex weak$^{\ast }$%
-compact state space of antiliminal\emph{\ }and simple $C^{\ast }$-algebras.
In this case, by \cite[Lemma 11.2.4]{Dixmier}, the state space is known to
have a weak$^{\ast }$-dense set of extreme points. For some UHF (uniformly
hyperfinite) $C^{\ast }$-algebras, this space can even be represented by the
weak$^{\ast }$-closure of a strictly increasing, countable, family of
Poulsen simplices\footnote{%
It is a (unique, up to an affine homeomorphism, universal, homogeneous)
metrizable simplex with dense extreme boundary.} \cite{Poulsen}, as
explained in \cite{BruPedra4} for CAR $C^{\ast }$-algebras. It is
interesting to know whether these disconcerting properties of \emph{%
physically relevant}\footnote{%
E.g., quantum spin systems, fermions (like electrons) on lattices.} systems
are purely accidental, or not.

This question is reminiscent of the celebrated Wonderland theorem, by Simon
\cite{Simon} in 1995, and papers \cite{cesar1,cesar2,cesar3,cesar4}, by
Carvalho and de Oliveira in 2016-2018, showing that apparently peculiar
features of Schr\"{o}dinger operators in quantum mechanics, like the
existence of a purely singular spectrum, are in fact generic. See, e.g.,
\cite{cesar1}. I.e., the exceptional-looking property turn out to be the
rule from the topological viewpoint. We study a similar genericity issue on
convex weak$^{\ast }$-compact sets, by using the weak$^{\ast }$-Hausdorff
hypertopology introduced by us in the present paper.

In 1959, Klee shows \cite{Klee} that, for convex norm-compact sets within a
Banach space, the property of having a dense set of extreme points is \emph{%
generic} in infinite dimensions. More precisely, by \cite[Proposition 2.1,
Theorem 2.2]{Klee}, the set of all such convex compact subsets of an
infinite-dimensional separable\footnote{\cite[Proposition 2.1, Theorem 2.2]%
{Klee} seem to lead to the asserted property for all (possibly
non-separable) Banach spaces, as claimed in \cite%
{Klee,FonfLindenstrauss,infinite dim convexity}. However, \cite[Theorem 1.5]%
{Klee}, which assumes the separability of the Banach space, is clearly
invoked to prove the corresponding density stated in \cite[Theorem 2.2]{Klee}%
. We do not know how to remove the a priori separability condition.} Banach
space $\mathcal{Y}$ is generic\footnote{%
That is, the complement of a meagre set, i.e., a nowhere dense set.} in the
complete metric space of compact convex subsets of $\mathcal{Y}$, endowed
with the well-known Hausdorff metric (hyper)topology. Klee's result has been
refined in 1998 by Fonf and Lindenstraus \cite[Section 4]{FonfLindenstrauss}
for bounded norm-closed (but not necessarily norm-compact) convex subsets of
$\mathcal{Y}$ having so-called empty quasi-interior (as a necessary
condition). In this case, \cite[Theorem 4.3]{FonfLindenstrauss} shows that
such sets can be approximated in the Hausdorff metric topology by closed
convex sets with a norm-dense set of strongly exposed points\footnote{$x\in
K $ is a strongly exposed point of a convex set $K\subseteq \mathcal{Y}$
when there is $f\in \mathcal{Y}^{\ast }$ satisfying $f(x)=1$ and such that
the diameter of $\{y\in K:f(y)\geq 1-\varepsilon \}$ tends to $0$ as $%
\varepsilon \rightarrow 0^{+}$. (Strongly) exposed points are extreme
elements of $K$.}. See, e.g., \cite[Section 7]{infinite dim convexity} for a
recent review on this subject.

In the present paper, we demonstrate the same genericity in the dual space $%
\mathcal{X}^{\ast }$, endowed with its weak$^{\ast }$-topology, of an
infinite-dimensional, separable topological vector space $\mathcal{X}$. The
hypertopology used in \cite{Klee,FonfLindenstrauss} is the Hausdorff metric
topology, induced by the Hausdorff distance associated with the norm of a
Banach space. Here, $\mathcal{X}^{\ast }$ is only endowed with the weak$%
^{\ast }$ topology and the Hausdorff metric (hyper)topology used in \cite%
{Klee,FonfLindenstrauss} is naturally replaced by the weak$^{\ast }$%
-Hausdorff hypertopology.

If $\mathcal{X}$ is a separable Banach space, note that one can use the
standard norm topology on $\mathcal{X}^{\ast }$ for continuous linear
functionals and directly apply previous results \cite{Klee,FonfLindenstrauss}
to the separable Banach space $\mathcal{X}^{\ast }$. This is not anymore
possible if one considers the weak$^{\ast }$-topology. In particular, \cite[%
Theorem 4.3]{FonfLindenstrauss} cannot be invoked in this specific case
because, in general, weak$^{\ast }$-compact sets do not have an empty
interior, in the sense of the norm topology. Anyway, our study is done for
locally convex spaces $\mathcal{X}^{\ast }$, as the dual spaces -- endowed
with the weak$^{\ast }$ topology -- of separable, infinite-dimensional
topological vector spaces $\mathcal{X}$, and considerably extends previous
ones \cite{Klee,FonfLindenstrauss}. We also prove similar results within weak%
$^{\ast }$-closed sets of \emph{positive} functionals of $\mathcal{X}^{\ast
} $, provided some positive (convex) cone in $\mathcal{X}$ is given. We
summarize our main results on infinite-dimensional convexity as follow:

\begin{itemize}
\item The set $\mathcal{D}_{0}$ of convex weak$^{\ast }$-compact subsets
with weak$^{\ast }$-dense set of \emph{exposed} points is weak$^{\ast }$%
-Hausdorff dense in the set of all convex weak$^{\ast }$-compact subsets of
an infinite-dimensional absolute polar in $\mathcal{X}^{\ast }$. $\mathcal{X}
$ is here a real or complex topological vector space and if $\mathcal{X}$ is
additionally a separable Gateaux-differentiability space then $\mathcal{D}%
_{0}$ is also a $G_{\delta }$ set. Recall that a Gateaux-differentiability
space $\mathcal{X}$ is defined to be a topological vector space on which
every continuous convex real-valued function with a nonempty open convex
subset as domain is Gateaux-differentiable on a dense subset of that domain.
See \cite{LarmanPhelps,Phelps-Asplund}.

\item Meanwhile, the set of \emph{exposed} points is proven to be weak$%
^{\ast }$-dense in the set of \emph{extreme} points of a convex weak$^{\ast
} $-compact subset of a Gateaux-differentiability space $\mathcal{X}$. This
corresponds to an extension of \cite[Theorem 6.2]{Phelps-Asplund}, which
only refers to Banach spaces $\mathcal{X}$. Our proof is quite direct and
thus, pedagogical while being very general. This result extends the
Straszewicz theorem proven in 1935 \cite{Straszewicz} by Straszewicz in $%
\mathbb{R}^{n}$, generalized by Klee in 1958 \cite[Theorems (2.1), (2.3)]%
{Klee-exposed} and in 1976 by Bair \cite[Theorem 1]{Bair} in arbitrary real
linear space for algebraically closed convex sets with so-called finite
\textquotedblleft copointure\textquotedblright\ \cite[Section II.5.1]{BairI}%
. Note that this last condition cannot be satisfied by weak$^{\ast }$%
-compact sets.

\item We refine our analysis in order to get similar results for the special
case of positive functionals of $\mathcal{X}^{\ast }$, by using the
decomposition of equicontinuous linear functionals into positive
equicontinuous components \cite{decomposition}, as originally proven by
Grosberg and Krein \cite{decomposition1} in 1939 for normed spaces and by
Bonsall \cite{decomposition2} in 1957 for locally convex real vector spaces.
\end{itemize}

As an application to the Banach case, we show, for an infinite-dimensional,
separable and unital $C^{\ast }$-algebra $\mathcal{X}$, that convex weak$%
^{\ast }$-compact sets of the corresponding state space, i.e., the space of
positive and normalized functionals of $\mathcal{X}^{\ast }$, have
generically a weak$^{\ast }$-dense set of extreme points. This demonstrates
that the fairly complicated geometrical structure of the state space of
physically relevant systems is \emph{not} accidental, but \emph{generic} for
infinite-dimensional, separable unital $C^{\ast }$-algebras. This is done by
using a very fine weak$^{\ast }$-type topology, namely the weak$^{\ast }$%
-Hausdorff hypertopology.

To conclude, it is not so easy to highlight a particular result of this
paper since various independent statements are proven along the current
studies on the weak$^{\ast }$-Hausdorff hypertopology and convex sets.
Propositions \ref{dddddddddddddddddd copy(1)}, \ref{convexity corrolary
copy(1)}, Corollaries \ref{connected sub-hyperspace copy(2)}, \ref{convexity
corrolary}, \ref{convexity corrolary copy(2)}, \ref{Solution selfbaby
copy(4)} and \ref{Solution selfbaby copy(13)} as well as Theorems \ref%
{Solution selfbaby copy(4)+1}, \ref{stra theorem}, \ref{genereic theorem}, %
\ref{Solution selfbaby copy(12)} are probably the most important of those
results. New mathematical concepts are also introduced in this paper: weak$%
^{\ast }$ hypertopologies (Definition \ref{weak hypertopology}), weak$^{\ast
}$-Hausdorff hypertopology (Definition \ref{hypertopology0}), and
immeasurable hypersets (Definition \ref{hypertopology0 copy(1)}). The paper
is organized as follows: Section \ref{Hypersets} gathers all the necessary
definitions to introduce in Section \ref{Weak* Hypertopologies} the class of
weak$^{\ast }$ hypertopologies as well as the important example of the weak$%
^{\ast }$-Hausdorff hypertopology. We next study, in Section \ref{Hausdorff
Hypertopology}, this particular weak$^{\ast }$ hypertopology in detail. In
Section \ref{Generic section}, using the results of Section \ref{Hausdorff
Hypertopology}, we study generic convex weak$^{\ast }$-compact sets in great
generality.

\section{Hyperspaces from Dual Spaces\label{Hypersets}}

\subsection{Hypersets and Hypermappings\label{Hypersets copy(1)}}

\noindent \underline{($\mathcal{X}^{\ast }$):} We denote by $\mathcal{X}$ a
topological $\mathbb{K}$-vector space with $\mathbb{K}=\mathbb{R},\mathbb{C}$%
. We define $\mathcal{O}\subseteq 2^{\mathcal{X}}$ to be the set of all\ $0$%
-neighborhoods $\mathcal{U}\subseteq \mathcal{X}$. Note that the definitions
of topological vector spaces found in the literature differ slightly from
each other. Those differences mostly concern the Hausdorff property. Here,
we use Rudin's definition \cite[Section 1.6]{Rudin}. In this case, the space
$\mathcal{X}$ is\emph{\ }Hausdorff, by \cite[Theorem 1.12]{Rudin}.

Let $\mathcal{X}^{\ast }$ be the (topological) dual space of $\mathcal{X}$.
To avoid peculiar cases, we assume throughout the paper that $\mathcal{X}%
^{\ast }$ \emph{separates points on} $\mathcal{X}$. This is for instance the
case for any (Hausdorff) locally convex space $\mathcal{X}$. See also \cite[%
Exercice 5 of Chap. 3]{Rudin} for other examples. As explained in \cite[%
Sections 3.8, 3.10, 3.14]{Rudin}, recall that the $\sigma (\mathcal{X}^{\ast
},\mathcal{X})$-topology -- called the weak$^{\ast }$ topology of $\mathcal{X%
}^{\ast }$ -- is the initial topology of the family $\{\hat{A}\}_{A\in
\mathcal{X}}$ of linear mappings from $\mathcal{X}^{\ast }$ to $\mathbb{K}$,
defined by
\begin{equation}
\hat{A}\left( \sigma \right) \doteq \sigma \left( A\right) \text{ },\text{%
\qquad }\sigma \in \mathcal{X}^{\ast },\ A\in \mathcal{X}\ .
\label{sdfsdfkljsdlfkj}
\end{equation}%
It is, by definition, the coarsest topology on $\mathcal{X}^{\ast }$ that
makes the mapping $\hat{A}$ continuous for every $A\in \mathcal{X}$. See
\cite[Section 3.8]{Rudin}. Here, the topology of $\mathcal{X}^{\ast }$ is,
by default, the weak$^{\ast }$ topology. In this case, $\mathcal{X}^{\ast }$
is a (Hausdorff)\textit{\ }locally convex space and its (topological) dual
space is $\mathcal{X}$: Any element of $\mathcal{X}^{\ast \ast }\equiv
\mathcal{X}$ is of the form (\ref{sdfsdfkljsdlfkj}). See \cite[Theorem 3.10]%
{Rudin}. By \cite[Theorem 1.37 and Remark 1.38 (b)]{Rudin}, the collection
of all finite intersections of the sets
\begin{equation}
\mathcal{V}_{A}\doteq \left\{ \sigma \in \mathcal{X}^{\ast }:\left\vert
\sigma \left( A\right) \right\vert <1\right\} \ ,\qquad A\in \mathcal{X},
\label{cylindre1}
\end{equation}%
is a convex balanced local base of neighborhoods at $0\in \mathcal{X}^{\ast }
$ for the weak$^{\ast }$ topology. In the sequel we use the notation%
\begin{equation}
\mathcal{V}_{\mathcal{A}_{n}}\doteq \bigcap\limits_{j=1}^{n}\mathcal{V}%
_{A_{j}}\ ,\qquad n\in \mathbb{N},\ \mathcal{A}_{n}\doteq \left(
A_{j}\right) _{j=1}^{n}\subseteq \mathcal{X}\ ,  \label{cylindre2}
\end{equation}%
and
\begin{equation}
\mathcal{C}\doteq \left\{ \mathcal{V}_{\mathcal{A}_{n}}:n\in \mathbb{N},\
\mathcal{A}_{n}\doteq \left( A_{j}\right) _{j=1}^{n}\subseteq \mathcal{X}%
\right\}   \label{cylinder}
\end{equation}%
for the set of cylinder sets constructed from finite sequences of $\mathcal{X%
}$.

\noindent \underline{(F):} As is usual in the theory of hyperspaces \cite%
{Beer}, we study the set of all nonempty closed subsets of $\mathcal{X}%
^{\ast }$ denoted by%
\begin{equation}
\mathbf{F}\left( \mathcal{X}^{\ast }\right) \doteq \left\{ F\subseteq
\mathcal{X}^{\ast }:F\neq \emptyset \text{ is weak}^{\ast }\text{-closed}%
\right\} \ .  \label{hyperspace F}
\end{equation}%
Similarly, $\mathbf{F}(\mathbb{K})$ is the set of nonempty closed (in the
sense of the absolute value) subsets of the field $\mathbb{K}=\mathbb{R},%
\mathbb{C}$. The complement of any subset $\mathbf{F}_{0}\subseteq \mathbf{F}%
(\mathcal{X}^{\ast })$ is denoted by%
\begin{equation}
\mathbf{F}^{c}\doteq \mathbf{F}\left( \mathcal{X}^{\ast }\right) \backslash
\mathbf{F}_{0}\ .  \label{complement}
\end{equation}%
The subset of all nonempty convex weak$^{\ast }$-closed subsets of $\mathcal{%
X}^{\ast }$ is defined by%
\begin{equation}
\mathbf{CF}\left( \mathcal{X}^{\ast }\right) \doteq \left\{ F\in \mathbf{F}%
\left( \mathcal{X}^{\ast }\right) :F\text{ is convex}\right\} \subseteq
\mathbf{F}\left( \mathcal{X}^{\ast }\right) \ .  \label{hyperspace C}
\end{equation}%
There is a natural mapping from $\mathbf{F}(\mathcal{X}^{\ast })$ to $%
\mathbf{CF}\left( \mathcal{X}^{\ast }\right) $: the \emph{weak}$^{\ast }$-%
\emph{closed convex hull operator}, $\overline{\mathrm{co}}:\mathbf{F}(%
\mathcal{X}^{\ast })\rightarrow \mathbf{F}(\mathcal{X}^{\ast })$, defined by
\begin{equation}
\overline{\mathrm{co}}\left( F\right) \doteq \overline{\mathrm{co}F}\
,\qquad F\in \mathbf{F}(\mathcal{X}^{\ast })\ ,  \label{closure operator}
\end{equation}%
where $\overline{\mathrm{co}F}$ is the weak$^{\ast }$-closure of the convex
hull of $F$ or, equivalently, the intersection of all weak$^{\ast }$-closed
convex subsets of $\mathcal{X}^{\ast }$ containing $F$. It is a \emph{%
closure\ }(or hull)\emph{\ }operator \cite[Definition 5.1]{Universal Algebra}
since it satisfies the following properties:

\begin{itemize}
\item For any $F\in \mathbf{F}(\mathcal{X}^{\ast })$, $F\subseteq \overline{%
\mathrm{co}}(F)$ (extensive);

\item For any $F\in \mathbf{F}(\mathcal{X}^{\ast })$, $\overline{\mathrm{co}}%
\left( \overline{\mathrm{co}}(F)\right) =\overline{\mathrm{co}}(F)$
(idempotent);

\item For any $F_{1},F_{2}\in \mathbf{F}(\mathcal{X}^{\ast })$ such that $%
F_{1}\subseteq F_{2}$, $\overline{\mathrm{co}}\left( F_{1}\right) \subseteq
\overline{\mathrm{co}}(F_{2})$ (isotone).
\end{itemize}

\noindent Such a closure operator has certainly been used in the past. It is
a composition of (i) an \emph{algebraic} (or finitary) closure operator \cite%
[Definition 5.4]{Universal Algebra} defined by $F\mapsto \mathrm{co}F$ with
(ii) a \emph{topological} (or Kuratowski) closure operator \cite[Chapter 1,
p. 43]{topology} defined by $F\mapsto \overline{F}$ on $\mathbf{F}(\mathcal{X%
}^{\ast })$.

As is usual, weak$^{\ast }$-closed subsets $F\in \mathbf{F}(\mathcal{X}%
^{\ast })$ satisfying $F=\overline{\mathrm{co}}(F)$ are, by definition, $%
\overline{\mathrm{co}}$\emph{-closed }sets and, obviously,
\begin{equation}
\mathbf{CF}\left( \mathcal{X}^{\ast }\right) =\overline{\mathrm{co}}\left(
\mathbf{F}\left( \mathcal{X}^{\ast }\right) \right)  \label{surjective}
\end{equation}%
is the set of all $\overline{\mathrm{co}}$-closed\emph{\ }sets. Below, the
weak$^{\ast }$-closed convex hull operator naturally appears in the study of
the Hausdorff property associated with the weak$^{\ast }$-Hausdorff
hypertopology (Definition \ref{hypertopology0}). \smallskip

\noindent \underline{(B):} For any cylinder $\mathcal{V\in C}$, define the
set of all nonempty, \textquotedblleft $\mathcal{V}$-bounded%
\textquotedblright , weak$^{\ast }$-closed subsets of $\mathcal{X}^{\ast }$
by%
\begin{equation}
\mathbf{B}_{\mathcal{V}}\left( \mathcal{X}^{\ast }\right) \doteq \left\{
B\in \mathbf{F}\left( \mathcal{X}^{\ast }\right) :B\subseteq \lambda
\mathcal{V}\text{ for some }\lambda \in \mathbb{R}^{+}\right\} \subseteq
\mathbf{F}\left( \mathcal{X}^{\ast }\right) \ .  \label{Hyperspace BV}
\end{equation}%
Clearly, $\mathbf{B}_{\mathcal{X}^{\ast }}(\mathcal{X}^{\ast })=\mathbf{F}(%
\mathcal{X}^{\ast })$ and
\begin{equation}
\mathbf{B}_{\mathcal{V}_{1}\cap \mathcal{V}_{2}}\left( \mathcal{X}^{\ast
}\right) =\mathbf{B}_{\mathcal{V}_{1}}\left( \mathcal{X}^{\ast }\right) \cap
\mathbf{B}_{\mathcal{V}_{2}}\left( \mathcal{X}^{\ast }\right) \ ,\qquad
\mathcal{V}_{1},\mathcal{V}_{2}\in \mathcal{C}\ .  \label{dddddd}
\end{equation}%
Then, the (nonempty) set of all nonempty (topologically) bounded weak$^{\ast
}$-closed subsets of $\mathcal{X}^{\ast }$ is defined by%
\begin{equation}
\mathbf{B}\left( \mathcal{X}^{\ast }\right) \doteq \bigcap_{\mathcal{V}\in
\mathcal{C}}\mathbf{B}_{\mathcal{V}}\left( \mathcal{X}^{\ast }\right)
\subseteq \mathbf{F}\left( \mathcal{X}^{\ast }\right) \ .
\label{hyperspace B}
\end{equation}%
See, e.g., \cite[Section 1.6 and Theorem 1.37 (b)]{Rudin}. Let
\begin{equation}
\mathbf{CB}\left( \mathcal{X}^{\ast }\right) \doteq \mathbf{CF}\left(
\mathcal{X}^{\ast }\right) \cap \mathbf{B}\left( \mathcal{X}^{\ast }\right) =%
\overline{\mathrm{co}}\left( \mathbf{B}\left( \mathcal{X}^{\ast }\right)
\right) \ ,  \label{CB}
\end{equation}%
where the last equality follows from the fact that the weak$^{\ast }$%
-closure of the convex hull of a weak$^{\ast }$-closed bounded set is
bounded, by the triangle inequality.

If $\mathcal{X}$ is a Banach space, then $\mathbf{B}(\mathcal{X}^{\ast })$
is nothing else than the set of all nonempty norm-bounded weak$^{\ast }$%
-closed subsets of $\mathcal{X}^{\ast }$, by the uniform boundedness
principle \cite[Theorems 2.4 and 2.5]{Rudin}. \medskip

\noindent \underline{(K):} The set of all nonempty weak$^{\ast }$-compact
subsets of $\mathcal{X}^{\ast }$ is denoted by%
\begin{equation}
\mathbf{K}\left( \mathcal{X}^{\ast }\right) \doteq \left\{ K\in \mathbf{F}%
\left( \mathcal{X}^{\ast }\right) :K\text{ is weak}^{\ast }\text{-compact}%
\right\} \subseteq \mathbf{F}\left( \mathcal{X}^{\ast }\right) \ .
\label{hyperspace K}
\end{equation}%
(Recall that in a Hausdorff space every compact set is closed.) By \cite[%
Theorem 1.15 (b)]{Rudin}, weak$^{\ast }$-compact subsets of $\mathcal{X}%
^{\ast }$ are bounded:
\begin{equation*}
\mathbf{K}\left( \mathcal{X}^{\ast }\right) \subseteq \mathbf{B}\left(
\mathcal{X}^{\ast }\right) \ .
\end{equation*}%
Since any convex weak$^{\ast }$-compact set $K$ satisfies $\overline{\mathrm{%
co}}\left( K\right) =K$, note that
\begin{equation}
\mathbf{CK}\left( \mathcal{X}^{\ast }\right) \doteq \mathbf{CF}\left(
\mathcal{X}^{\ast }\right) \cap \mathbf{K}\left( \mathcal{X}^{\ast }\right)
\subseteq \overline{\mathrm{co}}\left( \mathbf{K}\left( \mathcal{X}^{\ast
}\right) \right) \subseteq \mathbf{CB}\left( \mathcal{X}^{\ast }\right) \ .
\label{hyperspace CK}
\end{equation}%
In general, as discussed for instance in \cite[Section 3.18]{Rudin}, the
equality
\begin{equation*}
\mathbf{CK}\left( \mathcal{X}^{\ast }\right) =\overline{\mathrm{co}}\left(
\mathbf{K}\left( \mathcal{X}^{\ast }\right) \right)
\end{equation*}%
does \emph{not} hold true: For instance, take as topological $\mathbb{R}$%
-vector space the space $\mathcal{X}\subseteq \mathbb{R}^{\mathbb{N}}$ of
sequences $A\equiv (A_{n})_{n\in \mathbb{N}}\subseteq \mathbb{R}$ that
eventually vanishes, along with the supremum norm. For any $m\in \mathbb{N}$%
, define the continuous linear functionals $\sigma _{m}\in \mathcal{X}^{\ast
}$ by%
\begin{equation*}
\sigma _{m}\left( A\right) =2^{m}A_{m}\ ,\qquad A\equiv (A_{n})_{n\in
\mathbb{N}}\in \mathbb{R}^{\mathbb{N}}\ .
\end{equation*}%
As $m\rightarrow \infty $, $\sigma _{m}$ converges in the weak$^{\ast }$
topology to zero. In particular, the set
\begin{equation*}
K=\left\{ \sigma _{m}:m\in \mathbb{N}\right\} \cup \left\{ 0\right\}
\subseteq \mathcal{X}^{\ast }
\end{equation*}%
is weak$^{\ast }$-compact, i.e., $K\in \mathbf{K}(\mathcal{X}^{\ast })$,
but, by \cite[Theorem 2.9]{Rudin} with $\Gamma =\{\hat{\delta}_{m}\}_{m\in
\mathbb{N}}\subseteq \mathcal{X}^{\ast \ast }$, where $\hat{\delta}%
_{m}\left( \sigma \right) =\sigma \left( \delta _{m}\right) $ and $\delta
_{m}\equiv (\delta _{m,n})_{n\in \mathbb{N}}\in \mathcal{X}$ is defined by $%
\delta _{m,m}=1$ and $0$ otherwise, one deduces that $\overline{\mathrm{co}}%
K\notin \mathbf{K}(\mathcal{X}^{\ast })$. This example is taken from \cite[%
Chapter II, Section 10]{Khaleelulla}.

Observe that the weak$^{\ast }$-closed convex hull operator $\overline{%
\mathrm{co}}$ yields a notion of compactness, defined as follows: A set $%
K\in \mathbf{F}(\mathcal{X}^{\ast })$ is $\overline{\mathrm{co}}$\emph{%
-compact} if it is $\overline{\mathrm{co}}$-closed and each family of $%
\overline{\mathrm{co}}$-closed subsets of $K$ which has the finite
intersection property has a non-empty intersection. Compare this definition
with \cite[Chapter 5, Theorem 1]{topology}. The set $\mathbf{CK}(\mathcal{X}%
^{\ast })$ of all nonempty convex weak$^{\ast }$-compact sets belongs to the
set%
\begin{equation*}
\mathbf{K}_{\overline{\mathrm{co}}}\left( \mathcal{X}^{\ast }\right) \doteq
\left\{ K\in \mathbf{F}\left( \mathcal{X}^{\ast }\right) :K\text{ is }%
\overline{\mathrm{co}}\text{-compact}\right\}
\end{equation*}%
of $\overline{\mathrm{co}}$-compact sets:

\begin{proposition}[Space of $\overline{\mathrm{co}}$-compact sets]
\label{convexity lemma copy(3)}\mbox{ }\newline
Let $\mathcal{X}$ be a topological $\mathbb{K}$-vector space with $\mathbb{K}%
=\mathbb{R},\mathbb{C}$. Then,
\begin{equation*}
\mathbf{CK}\left( \mathcal{X}^{\ast }\right) \subseteq \mathbf{K}_{\overline{%
\mathrm{co}}}\left( \mathcal{X}^{\ast }\right) \subseteq \mathbf{CB}\left(
\mathcal{X}^{\ast }\right) \ .
\end{equation*}
\end{proposition}

\begin{proof}
By \cite[Chapter 5, Theorem 1]{topology}, $\mathbf{CK}(\mathcal{X}^{\ast
})\subseteq \mathbf{K}_{\overline{\mathrm{co}}}(\mathcal{X}^{\ast })$. Now,
take any $\overline{\mathrm{co}}$-compact element $K\in \mathbf{K}_{%
\overline{\mathrm{co}}}(\mathcal{X}^{\ast })$. If $K$ is not bounded, then
there is $A\in \mathcal{X}$ such that $\hat{A}(K)\subseteq \mathbb{K}$ is
not bounded, where we recall that $\hat{A}:\mathcal{X}^{\ast }\rightarrow
\mathbb{K}$ is the weak$^{\ast }$-continuous linear functional defined by (%
\ref{sdfsdfkljsdlfkj}). Without loss of generality, assume that $\mathrm{Re}%
\{\hat{A}(K)\}$ is not bounded from above. Define, for every $n\in \mathbb{N}
$, the set
\begin{equation*}
K_{n}\doteq \left\{ \sigma \in K:\mathrm{Re}\{\hat{A}(\sigma )\}\geq
n\right\} \ .
\end{equation*}%
Clearly, by convexity of $K$, $K_{n}$ is a convex weak$^{\ast }$-closed
subset of $K$ and the family $(K_{n})_{n\in \mathbb{N}}$ has the finite
intersection property, but, by construction,
\begin{equation*}
\bigcap_{n\in \mathbb{N}}K_{n}=\emptyset \ .
\end{equation*}%
(The intersection of preimages is the preimage of the intersection.) This
contradicts the fact that $K$ is $\overline{\mathrm{co}}$-compact.
Therefore, $K\in \mathbf{K}_{\overline{\mathrm{co}}}(\mathcal{X}^{\ast })$
is bounded and, being $\overline{\mathrm{co}}$-compact, it is also convex.
Consequently, $\mathbf{K}_{\overline{\mathrm{co}}}(\mathcal{X}^{\ast
})\subseteq \mathbf{CB}(\mathcal{X}^{\ast })$.
\end{proof}

\noindent Unless $\mathcal{X}$ is a Banach space (see (\ref{Banach})), we do
not expect the equality $\mathbf{CK}(\mathcal{X}^{\ast })=\mathbf{K}_{%
\overline{\mathrm{co}}}(\mathcal{X}^{\ast })$ to hold true. \medskip

\noindent \underline{(U):} Recall that the closed convex hull of a compact
set of an infinite-dimensional topological space need not be compact and so,
the set $\mathbf{K}(\mathcal{X}^{\ast })$ of all nonempty weak$^{\ast }$%
-compact subsets of $\mathcal{X}^{\ast }$ is generally not invariant under
the weak$^{\ast }$-closed convex hull operator $\overline{\mathrm{co}}$, as
explicitly demonstrated after Equation (\ref{hyperspace CK}). This motivates
the introduction of a specific, albeit still large, class of weak$^{\ast }$%
-compact sets, in relation with \emph{absolute polars} of $0$-neighborhoods
in $\mathcal{X}$: the absolute polar $\mathcal{U}^{\circ }$ of any $\mathcal{%
U}\in \mathcal{O}$ is defined by%
\begin{equation}
\mathcal{U}^{\circ }\doteq \left\{ \sigma \in \mathcal{X}^{\mathcal{\ast }%
}:\left\vert \sigma \left( A\right) \right\vert \leq 1\text{ for every }A\in
\mathcal{U}\right\} \in \mathbf{CF}\left( \mathcal{X}^{\ast }\right) \ ,
\label{polar}
\end{equation}%
where we recall that $\mathcal{O}\subseteq 2^{\mathcal{X}}$ is the set of
all (open) $0$-neighborhoods. As is well-known, it is a nonempty (balanced)
convex weak$^{\ast }$-closed set containing $0\in \mathcal{X}^{\ast }$, but,
what's more, a weak$^{\ast }$-compact set, by the Banach-Alaoglu theorem
\cite[Theorem 3.15]{Rudin}. For any $\mathcal{U}\in \mathcal{O}$, let
\begin{equation}
\mathbf{U}_{\mathcal{U}}\left( \mathcal{X}^{\ast }\right) \doteq \left\{
U\in \mathbf{F}\left( \mathcal{X}^{\ast }\right) :U\subseteq \mathcal{U}%
^{\circ }\right\} \subseteq \mathbf{K}\left( \mathcal{X}^{\ast }\right) \ .
\label{U0}
\end{equation}%
Clearly, one has%
\begin{equation*}
\mathbf{U}_{\mathcal{U}_{1}}\left( \mathcal{X}^{\ast }\right) \cup \mathbf{U}%
_{\mathcal{U}_{2}}\left( \mathcal{X}^{\ast }\right) \subseteq \mathbf{U}_{%
\mathcal{U}_{1}\cap \mathcal{U}_{2}}\left( \mathcal{X}^{\ast }\right) \
,\qquad \mathcal{U}_{1},\mathcal{U}_{2}\in \mathcal{O}\ .
\end{equation*}%
We next denote by
\begin{equation}
\mathbf{U}\left( \mathcal{X}^{\ast }\right) \doteq \bigcup\limits_{\mathcal{U%
}\in \mathcal{O}}\mathbf{U}_{\mathcal{U}}\left( \mathcal{X}^{\ast }\right)
\subseteq \mathbf{K}\left( \mathcal{X}^{\ast }\right) \subseteq \mathbf{B}%
\left( \mathcal{X}^{\ast }\right)  \label{U}
\end{equation}%
the set of all nonempty, uniformly bounded in a $0$-neighborhood, weak$%
^{\ast }$-closed subsets of $\mathcal{X}^{\ast }$.

By the triangle inequality and the closure operator property of $\overline{%
\mathrm{co}}$, together with the convexity and weak$^{\ast }$-closedness of
absolute polars (\ref{polar}), observe that the weak$^{\ast }$-closed convex
hull operator $\overline{\mathrm{co}}$ preserves $\mathbf{U}_{\mathcal{U}}(%
\mathcal{X}^{\ast })$ for every $\mathcal{U}\in \mathcal{O}$:
\begin{equation}
\mathbf{CU}_{\mathcal{U}}\left( \mathcal{X}^{\ast }\right) \doteq \mathbf{CF}%
\left( \mathcal{X}^{\ast }\right) \cap \mathbf{U}_{\mathcal{U}}\left(
\mathcal{X}^{\ast }\right) =\overline{\mathrm{co}}\left( \mathbf{U}_{%
\mathcal{U}}\left( \mathcal{X}^{\ast }\right) \right) \subseteq \mathbf{CB}%
\left( \mathcal{X}^{\ast }\right) \ .  \label{CU0}
\end{equation}%
Hence,%
\begin{equation}
\mathbf{CU}\left( \mathcal{X}^{\ast }\right) \doteq \mathbf{CF}\left(
\mathcal{X}^{\ast }\right) \cap \mathbf{U}\left( \mathcal{X}^{\ast }\right) =%
\overline{\mathrm{co}}\left( \mathbf{U}\left( \mathcal{X}^{\ast }\right)
\right) \subseteq \mathbf{CK}\left( \mathcal{X}^{\ast }\right) \subseteq
\mathbf{CB}\left( \mathcal{X}^{\ast }\right) \ .  \label{CU}
\end{equation}%
Compare with (\ref{hyperspace CK}).

If $\mathcal{X}$ is a Banach space then $\mathbf{B}(\mathcal{X}^{\ast })$
equals the set of all nonempty norm-bounded weak$^{\ast }$-closed subsets of
$\mathcal{X}^{\ast }$, because, in this case, a subset of $\mathcal{X}^{\ast
}$ is norm-bounded and weak$^{\ast }$-closed iff it is weak$^{\ast }$%
-compact, as is well-known. See, e.g., \cite[Proposition 1.2.9]{Beer}. This
fact is a consequence of the uniform boundedness principle \cite[Theorems
2.4 and 2.5]{Rudin} and the Banach-Alaoglu theorem \cite[Theorem 3.15]{Rudin}%
, since the absolute polar of a \emph{norm-}closed ball of radius $D$ in $%
\mathcal{X}$ is a a \emph{norm-}closed ball of radius $D^{-1}$ in $\mathcal{X%
}^{\ast }$. In particular, in this situation, absolute polars can be
replaced with norm-closed balls in $\mathcal{X}^{\ast }$. In fact, if $%
\mathcal{X}$ is a Banach space then%
\begin{equation}
\mathbf{U}\left( \mathcal{X}^{\ast }\right) =\mathbf{K}\left( \mathcal{X}%
^{\ast }\right) =\mathbf{B}\left( \mathcal{X}^{\ast }\right)  \label{Banach0}
\end{equation}%
is nothing else than the set of all nonempty norm-bounded weak$^{\ast }$%
-closed subsets of $\mathcal{X}^{\ast }$ and, using the weak$^{\ast }$%
-closed convex hull operator $\overline{\mathrm{co}}$, we deduce that%
\begin{equation}
\mathbf{CU}\left( \mathcal{X}^{\ast }\right) =\mathbf{CK}\left( \mathcal{X}%
^{\ast }\right) =\overline{\mathrm{co}}\left( \mathbf{K}\left( \mathcal{X}%
^{\ast }\right) \right) =\mathbf{K}_{\overline{\mathrm{co}}}\left( \mathcal{X%
}^{\ast }\right) =\mathbf{CB}\left( \mathcal{X}^{\ast }\right) \ ,
\label{Banach}
\end{equation}%
by (\ref{CB}), (\ref{hyperspace CK}), (\ref{CU}) and Proposition \ref%
{convexity lemma copy(3)}. Thus, if $\mathcal{X}$ is a Banach space,
Proposition \ref{convexity lemma copy(3)} gives an elegant abstract
characterization of $\mathbf{CK}(\mathcal{X}^{\ast })$, only expressed in
terms of a closure operator, namely the weak$^{\ast }$-closed convex hull
operator. It demonstrates a first connection with other mathematical fields,
in particular with mathematical logics where fascinating applications of
closure operators have been developed, already by Tarski himself during the
1930's.

\subsection{Weak$^{\ast }$ Hypertopologies\label{Weak* Hypertopologies}}

All hypersets (i.e., sets of closed sets) in Section \ref{Hypersets copy(1)}
can be endowed with hypertopologies. An \emph{hyper}topology is a topology
such that any net $(\sigma _{j})_{j\in J}$ in the primordial space (here the
dual space $\mathcal{X}^{\ast }$ or the field $\mathbb{K}=\mathbb{R},\mathbb{%
C}$) converges to an element $\sigma $ iff the net $(\{\sigma _{j}\})_{j\in
J}$ converges to $\{\sigma \}$ in the corresponding hyperspace (here $%
\mathbf{F}(\mathcal{X}^{\ast })$ or $\mathbf{F}(\mathbb{K})$). Recall that
there are various standard hypertopologies on general sets of nonempty
closed subsets of a complete metric space $(\mathcal{Y},d)$: the Fell,
Vietoris, Wijsman, proximal or locally finite hypertopologies, to name a few
well-known examples. See, e.g., \cite{Beer}.

None of these well-known hypertopologies is used here for $\mathbf{F}(%
\mathcal{X}^{\ast })$, for the weak$^{\ast }$ topology of $\mathcal{X}^{\ast
}$ is generally not metrizable. However, all these hypertopologies
associated with a complete metric space can be used to define, in a
systematic and very natural way that is similar to the weak$^{\ast }$
topology of $\mathcal{X}^{\ast }$, a \emph{new} class of hypertopologies on
the set $\mathbf{F}(\mathcal{X}^{\ast })$ of all nonempty weak$^{\ast }$%
-closed subsets of $\mathcal{X}^{\ast }$: At any fixed $A\in \mathcal{X}$,
we define the mapping $\mathbf{A}:\mathbf{F}(\mathcal{X}^{\ast })\rightarrow
\mathbf{F}(\mathbb{K})$ by%
\begin{equation}
\mathbf{A}\left( F\right) \doteq \overline{\hat{A}\left( F\right) }\doteq
\overline{\left\{ \sigma \left( A\right) :\sigma \in F\right\} }\text{ },%
\text{\qquad }F\in \mathbf{F}\left( \mathcal{X}^{\ast }\right) \ .
\label{map A}
\end{equation}%
See Equation (\ref{sdfsdfkljsdlfkj}). For a fixed hypertopology on $\mathbf{F%
}(\mathbb{K})$, the associated weak$^{\ast }$ hypertopology of $\mathbf{F}(%
\mathcal{X}^{\ast })$ is the coarsest topology on $\mathbf{F}(\mathcal{X}%
^{\ast })$ that makes the mapping $\mathbf{A}$ continuous for every $A\in
\mathcal{X}$:

\begin{definition}[Weak$^{\ast }$ hypertopologies]
\label{weak hypertopology}\mbox{ }\newline
Pick some hypertopology on $\mathbf{F}(\mathbb{K})$. The associated weak$%
^{\ast }$ hypertopology $\tau $ on $\mathbf{F}(\mathcal{X}^{\ast })$ is the
initial topology of the family $\{\mathbf{A}\}_{A\in \mathcal{X}}$ of
mappings from $\mathbf{F}(\mathcal{X}^{\ast })$ to $\mathbf{F}(\mathbb{K})$
defined by (\ref{map A}). That is, $\tau $ is the collection of all unions
of finite intersections of sets $\mathbf{A}^{-1}(V)$ with $A\in \mathcal{X}$
and $V$ open in $\mathbf{F}(\mathbb{K})$.
\end{definition}

\noindent By construction, $\tau $ is obviously an hypertopology on $\mathbf{%
F}(\mathcal{X}^{\ast })$, keeping in mind that the topology on $\mathcal{X}%
^{\ast }$ is, by default, the weak$^{\ast }$ topology, i.e., the initial
topology of the family $\{\hat{A}\}_{A\in \mathcal{X}}$ of linear mappings
from $\mathcal{X}^{\ast }$ to $\mathbb{K}$ defined by (\ref{sdfsdfkljsdlfkj}%
).

The Fell, Vietoris, Wijsman, proximal or locally finite hypertopologies, the
Hausdorff metric topology, etc., lead to various weak$^{\ast }$
hypertopologies. They belong to a new class of hypertopologies, that is, the
class of weak$^{\ast }$ hypertopologies, which does not seem to have been
systematically considered in the past. Note however that the scalar topology
described in \cite[Section 4.3]{Beer}, when defined, mutatis mutandis, for a
locally convex space like $\mathcal{X}^{\ast }$, is retrospectively an
example\footnote{%
And the unique one we are aware of.} of a weak$^{\ast }$ hypertopology.

In the case that the hypertopology on $\mathbf{F}(\mathbb{K})$ is
metrizable, the corresponding weak$^{\ast }$ hypertopology on $\mathbf{F}(%
\mathcal{X}^{\ast })$ is the topology generated by the family of
pseudometrics:%
\begin{equation}
d^{(A)}(F_{1},F_{2})\doteq \mathfrak{d}\left( \mathbf{A}\left( F_{1}\right) ,%
\mathbf{A}\left( F_{2}\right) \right) \ ,\qquad F_{1},F_{2}\in \mathbf{F}(%
\mathcal{X}^{\ast }),\ A\in \mathcal{X}\ ,  \label{ddddd}
\end{equation}%
where%
\begin{equation}
\mathfrak{d}:\mathbf{F}(\mathbb{K})\times \mathbf{F}(\mathbb{K})\rightarrow
\mathbb{R}_{0}^{+}\cup \left\{ \infty \right\}  \label{sdfsdfsf}
\end{equation}%
is the metric in $\mathbf{F}(\mathbb{K})$ generating its hypertopology.
(Recall also (\ref{map A}).) In this case, the weak$^{\ast }$ hypertopology
is a uniform topology, see, e.g., \cite[Chapter 6]{topology}. It is the
coarsest topology on $\mathbf{F}(\mathcal{X}^{\ast })$ that makes every
mapping $\mathbf{A}$, defined by (\ref{map A}) for $A\in \mathcal{X}$,
continuous. Equivalently, in this topology, an arbitrary net $(F_{j})_{j\in
J}\subseteq \mathbf{F}(\mathcal{X}^{\ast })$ converges to $F$ iff, for all $%
A\in \mathcal{X}$,
\begin{equation*}
\lim_{J}d^{(A)}(F_{j},F)=0\ .
\end{equation*}%
This condition defines a unique topology in $\mathbf{F}(\mathcal{X}^{\ast })$%
, by \cite[Chapter 2, Theorem 9]{topology}. In the context of
hypertopologies, note that it is natural to consider metrics $\mathfrak{d}$
taking values in the extended positive reals $\mathbb{R}_{0}^{+}\cup \left\{
\infty \right\} $, like (\ref{Hausdorf}) below.

The most well-studied and well-known hypertopology associated with a metric
space $(\mathcal{Y},d)$ is the Hausdorff metric topology \cite[Definition
3.2.1]{Beer}. It is generated by the Hausdorff distance between two sets $%
Y_{1},Y_{2}\subseteq \mathcal{Y}$:
\begin{equation}
\mathfrak{d}_{H}\left( Y_{1},Y_{2}\right) \doteq \max \left\{ \sup_{x_{1}\in
Y_{1}}\inf_{x_{2}\in Y_{2}}d\left( x_{1},x_{2}\right) ,\sup_{x_{2}\in
Y_{2}}\inf_{x_{1}\in Y_{1}}d\left( x_{1},x_{2}\right) \right\} \in \mathbb{R}%
_{0}^{+}\cup \left\{ \infty \right\} \ .  \label{Hausdorf}
\end{equation}%
In this case, the corresponding hyperspace of nonempty closed subsets of $%
\mathcal{Y}$ is complete iff the metric space $(\mathcal{Y},d)$ is complete.
See, e.g., \cite[Theorem 3.2.4]{Beer}. The Hausdorff metric topology is the
hypertopology used in \cite{Klee,FonfLindenstrauss}, the metric $d$ being
the one associated with the norm of a separable Banach space $\mathcal{Y}$,
in order to prove the density of the set of convex compact subsets of $%
\mathcal{Y}$ with dense extreme boundary. As a first and instructive example
of a weak$^{\ast }$ hypertopology, it is thus natural to study the weak$%
^{\ast }$ version of the Hausdorff metric topology. It corresponds to the
weak$^{\ast }$ hypertopology of Definition \ref{weak hypertopology} with $%
\mathbf{F}(\mathbb{K})$ endowed with the Hausdorff metric topology:

\begin{definition}[Weak$^{\ast }$-Hausdorff hypertopology]
\label{hypertopology0}\mbox{ }\newline
The weak$^{\ast }$-Hausdorff hypertopology on $\mathbf{F}(\mathcal{X}^{\ast
})$ is the topology induced by the family of Hausdorff pseudometrics $%
d_{H}^{(A)}$ defined, for all $A\in \mathcal{X}$, by%
\begin{equation}
d_{H}^{(A)}(F,\tilde{F})=\max \left\{ \sup_{\sigma \in F}\inf_{\tilde{\sigma}%
\in \tilde{F}}\left\vert \left( \sigma -\tilde{\sigma}\right) \left(
A\right) \right\vert ,\sup_{\tilde{\sigma}\in \tilde{F}}\inf_{\sigma \in
F}\left\vert \left( \sigma -\tilde{\sigma}\right) \left( A\right)
\right\vert \right\} \in \mathbb{R}_{0}^{+}\cup \left\{ \infty \right\}
,\quad F,\tilde{F}\in \mathbf{F}\left( \mathcal{X}^{\ast }\right) \ .
\label{def}
\end{equation}
\end{definition}

\noindent Note that, for any $F,\tilde{F}\in \mathbf{F}\left( \mathcal{X}%
^{\ast }\right) $,%
\begin{equation*}
\sup_{\sigma \in F}\inf_{\tilde{\sigma}\in \tilde{F}}\left\vert \left(
\sigma -\tilde{\sigma}\right) \left( A\right) \right\vert =\sup_{x\in \hat{A}%
\left( F\right) }\inf_{\tilde{x}\in \hat{A}(\tilde{F})}\left\vert x-\tilde{x}%
\right\vert =\sup_{x\in \hat{A}\left( F\right) }\inf_{\tilde{x}\in \mathbf{A}%
(\tilde{F})}\left\vert x-\tilde{x}\right\vert =\sup_{x\in \mathbf{A}\left(
F\right) }\inf_{\tilde{x}\in \mathbf{A}(\tilde{F})}\left\vert x-\tilde{x}%
\right\vert ,
\end{equation*}%
by the triangle inequality for the absolute value. (See, e.g., the arguments
justifying (\ref{dd}).)

Definition \ref{hypertopology0} is equivalent to Definition \ref{weak
hypertopology} with $\mathbf{F}(\mathbb{K})$ endowed with the Hausdorff
metric topology, as explained above for the more general case where $\mathbf{%
F}(\mathbb{K})$ is metrizable. To our knowledge, this hypertopology has not
been considered so far.

Here, $\mathbf{F}(\mathcal{X}^{\ast })$ and the subspaces $\mathbf{B}(%
\mathcal{X}^{\ast })$, $\mathbf{CB}(\mathcal{X}^{\ast })$, $\mathbf{K}(%
\mathcal{X}^{\ast })$, etc., are, by default, all endowed with the weak$%
^{\ast }$-Hausdorff hypertopology.

\section{The Weak$^{\ast }$-Hausdorff Hypertopology\label{Hausdorff
Hypertopology}}

\subsection{Boolean Algebras Associated with Immeasurable Hyperspaces}

Observe that one can only ensure that (\ref{def}) is finite only if $F,%
\tilde{F}\in \mathbf{B}\left( \mathcal{X}^{\ast }\right) \subseteq \mathbf{F}%
\left( \mathcal{X}^{\ast }\right) $. We show below that the weak$^{\ast }$%
-Hausdorff family of pseudometric $(d_{H}^{(A)})_{A\in \mathcal{X}}$ \emph{%
immeasurably} separates unbounded sets from bounded ones. This motivates the
following definition:

\begin{definition}[Immeasurable hypersets]
\label{hypertopology0 copy(1)}\mbox{ }\newline
Assume the existence of a metric $\mathfrak{d}$ satisfying (\ref{sdfsdfsf})
and generating the hypertopology in $\mathbf{F}(\mathbb{K})$. Let $\tau $ be
the corresponding weak$^{\ast }$ hypertopology on $\mathbf{F}(\mathcal{X}%
^{\ast })$. Two subsets\ $\mathbf{F}_{1},\mathbf{F}_{2}\subseteq \mathbf{F}(%
\mathcal{X}^{\ast })$ are said to be $\mathfrak{d}$-immeasurable if, for any
$F_{1}\in \mathbf{F}_{1}$ and $F_{2}\in \mathbf{F}_{2}$, there is $A\in
\mathcal{X}$ such that
\begin{equation*}
d^{(A)}(F_{1},F_{2})=\infty
\end{equation*}%
with the pseudometrics $d^{(A)}$ defined by (\ref{ddddd}). $\mathfrak{d}_{H}$%
-immeasurable sets are named here weak$^{\ast }$-Hausdorff-immeasurable
sets. See (\ref{Hausdorf}).
\end{definition}

A generally infinite collection of weak$^{\ast }$-Hausdorff-immeasurable
subspaces of $\mathbf{F}(\mathcal{X}^{\ast })$ is given by the subspaces $%
\mathbf{B}_{\mathcal{V}}(\mathcal{X}^{\ast })$ of all nonempty, $\mathcal{V}$%
-bounded, weak$^{\ast }$-closed subsets of $\mathcal{X}^{\ast }$, defined by
(\ref{Hyperspace BV}) for each cylinder $\mathcal{V\in C}$ (see (\ref%
{cylinder})).

\begin{proposition}[Pairs of immeasurable subhyperspaces]
\label{dddddddddddddddddd copy(1)}\mbox{ }\newline
Let $\mathcal{X}$ be a topological $\mathbb{K}$-vector space\footnote{%
Recall that all topological $\mathbb{K}$-vector spaces $\mathcal{X}$ in this
paper are Hausdorff, by \cite[Theorem 1.12]{Rudin}.} with $\mathbb{K}=%
\mathbb{R},\mathbb{C}$. For all cylinders $\mathcal{V}_{1},\mathcal{V}%
_{2}\in \mathcal{C}$ such that $\mathbf{B}_{\mathcal{V}_{1}}(\mathcal{X}%
^{\ast })\varsubsetneq \mathbf{B}_{\mathcal{V}_{2}}(\mathcal{X}^{\ast })$, $%
\mathbf{B}_{\mathcal{V}_{1}}(\mathcal{X}^{\ast })$ and $\mathbf{B}_{\mathcal{%
V}_{2}}(\mathcal{X}^{\ast })\backslash \mathbf{B}_{\mathcal{V}_{1}}(\mathcal{%
X}^{\ast })$ are weak$^{\ast }$-Hausdorff-immeasurable. In particular, $%
\mathbf{B}(\mathcal{X}^{\ast })$ and its complement $\mathbf{B}^{c}\left(
\mathcal{X}^{\ast }\right) $, defined by (\ref{complement}), are weak$^{\ast
}$-Hausdorff-immeasurable.
\end{proposition}

\begin{proof}
Take any $F\in \mathbf{B}_{\mathcal{V}_{2}}(\mathcal{X}^{\ast })\backslash
\mathbf{B}_{\mathcal{V}_{1}}(\mathcal{X}^{\ast })$, where $\mathcal{V}_{1},%
\mathcal{V}_{2}$ are two cylinders of the set $\mathcal{C}$ defined by (\ref%
{cylinder}) such that $\mathbf{B}_{\mathcal{V}_{1}}(\mathcal{X}^{\ast
})\varsubsetneq \mathbf{B}_{\mathcal{V}_{2}}(\mathcal{X}^{\ast })$. Then,
there is a net $(\sigma _{j})_{j\in J}\subseteq F$ and $A\in \mathcal{X}$
such that
\begin{equation}
\lim_{J}\left\vert \sigma _{j}\left( A\right) \right\vert =\infty \qquad
\text{and}\qquad \sup_{\tilde{\sigma}\in B}\left\vert \tilde{\sigma}\left(
A\right) \right\vert <\infty  \label{infinite}
\end{equation}%
for any $B\in \mathbf{B}_{\mathcal{V}_{1}}(\mathcal{X}^{\ast })$. By\
Definition \ref{hypertopology0} and the triangle inequality, for any $B\in
\mathbf{B}_{\mathcal{V}_{1}}(\mathcal{X}^{\ast })$ and $j\in J$,%
\begin{equation}
d_{H}^{(A)}\left( F,B\right) \geq \inf_{\tilde{\sigma}\in B}\left\vert
\left( \sigma _{j}-\tilde{\sigma}\right) \left( A\right) \right\vert \geq
\left\vert \sigma _{j}\left( A\right) \right\vert -\sup_{\tilde{\sigma}\in
B}\left\vert \tilde{\sigma}\left( A\right) \right\vert \ .  \label{infinite2}
\end{equation}%
By (\ref{infinite}), it follows that
\begin{equation*}
d_{H}^{(A)}(F,B)=\infty \ ,\qquad B\in \mathbf{B}_{\mathcal{V}_{1}}\left(
\mathcal{X}^{\ast }\right) \ .
\end{equation*}%
In other words, $\mathbf{B}_{\mathcal{V}_{1}}(\mathcal{X}^{\ast })$ and $%
\mathbf{B}_{\mathcal{V}_{2}}(\mathcal{X}^{\ast })\backslash \mathbf{B}_{%
\mathcal{V}_{1}}(\mathcal{X}^{\ast })$ are weak$^{\ast }$%
-Hausdorff-immeasurable.

If $F\in \mathbf{B}^{c}(\mathcal{X}^{\ast })$ then there is a cylinder\ $%
\mathcal{V}_{1}\neq \mathcal{X}^{\ast }$ such that $F\notin \mathbf{B}_{%
\mathcal{V}_{1}}(\mathcal{X}^{\ast })$ while any bounded set $B\in \mathbf{B}%
\left( \mathcal{X}^{\ast }\right) $ belongs to $\mathbf{B}_{\mathcal{V}_{1}}(%
\mathcal{X}^{\ast })$, by Equation (\ref{hyperspace B}). Using the previous
arguments for the special case $\mathcal{V}_{2}=\mathcal{X}^{\ast }\supseteq
\mathcal{V}_{1}$, we thus deduce that $\mathbf{B}(\mathcal{X}^{\ast })$ and
its complement $\mathbf{B}^{c}(\mathcal{X}^{\ast })$ are weak$^{\ast }$%
-Hausdorff-immeasurable.
\end{proof}

Observe that $\mathbf{B}_{\mathcal{X}^{\ast }}(\mathcal{X}^{\ast })=\mathbf{F%
}(\mathcal{X}^{\ast })$ and Proposition \ref{dddddddddddddddddd copy(1)}
applied to $\mathcal{V}_{2}=\mathcal{X}^{\ast }$ yield that, for each
cylinder $\mathcal{V\in C}\backslash \{\mathcal{X}^{\ast }\}$ such that $%
\mathbf{B}_{\mathcal{V}}(\mathcal{X}^{\ast })\varsubsetneq \mathbf{F}(%
\mathcal{X}^{\ast })$, the (nonempty) subhyperspaces $\mathbf{B}_{\mathcal{V}%
}(\mathcal{X}^{\ast })$ and its complement%
\begin{equation}
\mathbf{B}_{\mathcal{V}}^{c}\left( \mathcal{X}^{\ast }\right) \doteq \mathbf{%
F}\left( \mathcal{X}^{\ast }\right) \backslash \mathbf{B}_{\mathcal{V}%
}\left( \mathcal{X}^{\ast }\right)  \label{set of unbounded}
\end{equation}%
are weak$^{\ast }$-Hausdorff-immeasurable, like $\mathbf{B}(\mathcal{X}%
^{\ast })$ and $\mathbf{B}^{c}(\mathcal{X}^{\ast })$.

Additionally, the subspaces $\mathbf{B}_{\mathcal{V}}(\mathcal{X}^{\ast })$,
$\mathcal{V\in C}$, form a family of weak$^{\ast }$-Hausdorff clopen sets:

\begin{proposition}[Weak$^{\ast }$-Hausdorff-clopen subhyperspaces]
\label{convexity corrolary copy(1)}\mbox{ }\newline
Let $\mathcal{X}$ be a topological $\mathbb{K}$-vector space with $\mathbb{K}%
=\mathbb{R},\mathbb{C}$. Then,
\begin{equation}
\mathcal{C}\ell \doteq \left\{ \mathbf{B}_{\mathcal{V}}\left( \mathcal{X}%
^{\ast }\right) ,\mathbf{B}_{\mathcal{V}}^{c}\left( \mathcal{X}^{\ast
}\right) \right\} _{\mathcal{V\in C}}  \label{clopen}
\end{equation}%
is a family of weak$^{\ast }$-Hausdorff-closed subsets of $\mathbf{F}(%
\mathcal{X}^{\ast })$. In other words, $\mathcal{C}\ell $ is a family of
(nonempty) weak$^{\ast }$-Hausdorff-clopen\footnote{%
I.e., they are both open and closed in the weak$^{\ast }$-Hausdorff
hypertopology.} subsets of $\mathbf{F}(\mathcal{X}^{\ast })$.
\end{proposition}

\begin{proof}
Let $\mathcal{V=V}_{\mathcal{A}_{n}}$ with $\mathcal{A}_{n}\doteq \left(
A_{k}\right) _{k=1}^{n}\subseteq \mathcal{X}$, as defined by (\ref{cylindre1}%
)-(\ref{cylindre2}). Clearly, for any $B\in \mathbf{B}_{\mathcal{V}}\left(
\mathcal{X}^{\ast }\right) $, the set%
\begin{equation*}
\mathcal{F}\doteq \left\{ F\in \mathbf{F}\left( \mathcal{X}^{\ast }\right)
:\forall k\in \left\{ 1,\ldots ,n\right\} ,\ d_{H}^{(A_{k})}\left(
F,B\right) <1\right\}
\end{equation*}%
is a weak$^{\ast }$-Hausdorff neighborhood of $B$ in $\mathbf{F}(\mathcal{X}%
^{\ast })$. Additionally, by definition (\ref{def}) of the Hausdorff
pseudometric, $\mathcal{F}\subseteq \mathbf{B}_{\mathcal{V}}\left( \mathcal{X%
}^{\ast }\right) $. Therefore, $\mathbf{B}_{\mathcal{V}}\left( \mathcal{X}%
^{\ast }\right) $ is a weak$^{\ast }$-Hausdorff open set. Take now a net $%
\left( B_{j}\right) _{j\in J}\subseteq \mathbf{B}_{\mathcal{V}}\left(
\mathcal{X}^{\ast }\right) $ converging to $F\in \mathbf{F}(\mathcal{X}%
^{\ast })$ in the weak$^{\ast }$-Hausdorff topology. In particular, for some
$j\in J$,
\begin{equation*}
d_{H}^{(A_{k})}\left( F,B_{j}\right) <\infty \ ,\qquad k\in \{1,\ldots ,n\}\
.
\end{equation*}%
Again by (\ref{def}), it follows that $F\in \mathbf{B}_{\mathcal{V}}\left(
\mathcal{X}^{\ast }\right) $. Hence, $\mathbf{B}_{\mathcal{V}}\left(
\mathcal{X}^{\ast }\right) $ is a weak$^{\ast }$-Hausdorff closed set. Note
that a subset of a topological space is closed iff it contains all of its
limit points, by \cite[Chapter 1, Theorem 5]{topology}, and limit points
(also named cluster or accumulation points) of a set are precisely the
limits of (convergent) nets whose elements are in this set, by \cite[Chapter
2, Theorem 2]{topology}.
\end{proof}

In particular, any subspace $\mathbf{B}_{\mathcal{V}}(\mathcal{X}^{\ast })$,
$\mathcal{V\in C}$, has empty boundary\footnote{%
I.e., there is no element which is interior to neither $\mathbf{B}_{\mathcal{%
V}}(\mathcal{X}^{\ast })$ nor $\mathbf{B}_{\mathcal{V}}^{c}(\mathcal{X}%
^{\ast })$.} and thus, for infinite-dimensional topological vector spaces $%
\mathcal{X}$, the topological hyperspace $\mathbf{F}(\mathcal{X}^{\ast })$
has an infinite number of connected components in the weak$^{\ast }$%
-Hausdorff hypertopology. This leads to a whole collection of weak$^{\ast }$%
-Hausdorff-clopen sets, which could be used to form a Boolean algebra whose
lattice operations are given by the union and intersection, as is usual in
mathematical logics\footnote{%
See Stone's representation theorem for Boolean algebras.}. Such a study has
been performed in \cite{Bankston} for the hyperspace associated with a
Boolean compact topological space and the Vietoris hypertopology.

Note that the hyperspace $\mathbf{B}(\mathcal{X}^{\ast })$ of all nonempty
bounded weak$^{\ast }$-closed subsets of $\mathcal{X}^{\ast }$ is weak$%
^{\ast }$-Hausdorff-closed, the intersection of closed set being always
closed, but it is generally not weak$^{\ast }$-Hausdorff-open, even if $%
\mathbf{B}(\mathcal{X}^{\ast })$ and its complement $\mathbf{B}^{c}(\mathcal{%
X}^{\ast })$ are weak$^{\ast }$-Hausdorff-immeasurable. Additionally, for
any fixed cylinder $\mathcal{V\in C}$, observe that $\mathbf{B}_{\mathcal{V}%
}(\mathcal{X}^{\ast })$ is generally not a connected hyperspace. By (\ref%
{dddddd}), the weak$^{\ast }$-Hausdorff-clopen set $\mathbf{B}_{\mathcal{V}}(%
\mathcal{X}^{\ast })$ contains (possibly infinitely many) proper weak$^{\ast
}$-Hausdorff-clopen subsets, leading to many connected components. However,
the infimum (with respect to inclusion) of the family $\left\{ \mathbf{B}_{%
\mathcal{V}}\left( \mathcal{X}^{\ast }\right) \right\} _{\mathcal{V\in C}}$,
that is, $\mathbf{B}(\mathcal{X}^{\ast })$, is connected:

\begin{proposition}[$\mathbf{B}(\mathcal{X}^{\ast })$ as connected
subhyperspace]
\label{connected sub-hyperspace}\mbox{ }\newline
Let $\mathcal{X}$ be a topological $\mathbb{K}$-vector space with $\mathbb{K}%
=\mathbb{R},\mathbb{C}$. Then, the weak$^{\ast }$-Hausdorff-closed set $%
\mathbf{B}(\mathcal{X}^{\ast })$ is convex and path-connected. Moreover, it
is a connected component\footnote{%
That is, a maximal connected subset.} of $\mathbf{F}(\mathcal{X}^{\ast })$.
\end{proposition}

\begin{proof}
Take any $B_{0},B_{1}\in \mathbf{B}(\mathcal{X}^{\ast })$. Define the
mapping $f$ from $[0,1]$ to $\mathbf{B}(\mathcal{X}^{\ast })$ by%
\begin{equation}
f\left( \lambda \right) \doteq \left\{ \left( 1-\lambda \right) \sigma
_{0}+\lambda \sigma _{1}:\sigma _{0}\in B_{0},\ \sigma _{1}\in B_{1}\right\}
\ ,\qquad \lambda \in \left[ 0,1\right] \ .  \label{function}
\end{equation}%
(This already demonstrates that $\mathbf{B}(\mathcal{X}^{\ast })$ is
convex.) By Definition \ref{hypertopology0}, for any $\lambda _{1},\lambda
_{2}\in \lbrack 0,1]$,%
\begin{equation*}
d_{H}^{(A)}\left( f\left( \lambda _{1}\right) ,f\left( \lambda _{2}\right)
\right) \leq \left\vert \lambda _{2}-\lambda _{1}\right\vert \sup_{\sigma
\in (B_{0}-B_{1})}\left\vert \sigma \left( A\right) \right\vert \ ,\qquad
A\in \mathcal{X}\ .
\end{equation*}%
Note that, for all $B_{0},B_{1}\in \mathbf{B}(\mathcal{X}^{\ast })$ and $%
A\in \mathcal{X}$,
\begin{equation*}
\sup_{\sigma \in (B_{0}-B_{1})}\left\vert \sigma \left( A\right) \right\vert
<\infty \ .
\end{equation*}
So, the mapping $f$ is a continuous function from $[0,1]$ to $\mathbf{B}(%
\mathcal{X}^{\ast })$ with $f\left( 0\right) =B_{0}$ and $f\left( 1\right)
=B_{1}$. Therefore, $\mathbf{B}(\mathcal{X}^{\ast })$ is path-connected. The
image under a continuous mapping of a connected set is connected and, by
\cite[Chapter 1, Theorem 21]{topology}, $\mathbf{B}(\mathcal{X}^{\ast })$,
being path-connected, is connected. In particular, $\mathbf{B}(\mathcal{X}%
^{\ast })$ belongs to the connected component of any element $B\in \mathbf{B}%
(\mathcal{X}^{\ast })$, denoted by $\mathcal{B}$.

For any $B\in \mathbf{B}(\mathcal{X}^{\ast })$, define
\begin{equation*}
\mathcal{\tilde{B}}\doteq \bigcap \left\{ \mathcal{F}\subseteq \mathbf{F}(%
\mathcal{X}^{\ast }):\mathcal{F}\text{ is weak}^{\ast }\text{%
-Hausdorff-clopen and }B\in \mathcal{F}\right\} \ ,
\end{equation*}%
the so-called pseudocomponent of $B$. Since a clopen is never a proper
subset of a connected component, $\mathcal{B}\subseteq \mathcal{\tilde{B}}$.
By Proposition \ref{convexity corrolary copy(1)}, note that%
\begin{equation*}
\mathcal{\tilde{B}}\subseteq \bigcap_{\mathcal{V}\in \mathcal{C}}\mathbf{B}_{%
\mathcal{V}}\left( \mathcal{X}^{\ast }\right) =\mathbf{B}\left( \mathcal{X}%
^{\ast }\right) \subseteq \mathcal{B}\ .
\end{equation*}
This means that $\mathcal{B}=\mathbf{B}\left( \mathcal{X}^{\ast }\right) $.
\end{proof}

\begin{corollary}[$\mathbf{F}(\mathcal{X}^{\ast })$ as non-locally connected
hyperspace]
\label{connected sub-hyperspace copy(2)}\mbox{ }\newline
Let $\mathcal{X}$ be an infinite-dimensional topological $\mathbb{K}$-vector
space with $\mathbb{K}=\mathbb{R},\mathbb{C}$. Then $\mathbf{F}(\mathcal{X}%
^{\ast })$ is not locally connected.
\end{corollary}

\begin{proof}
If $\mathbf{B}(\mathcal{X}^{\ast })$\ is not weak$^{\ast }$-Hausdorff-open
then $\mathbf{F}(\mathcal{X}^{\ast })$ is not locally connected: Assume by
contradiction that $\mathbf{F}(\mathcal{X}^{\ast })$ is locally connected.
Then, because of \cite[Chap. 1, Problem (S), (a), p. 61]{topology}, any
connected component of $\mathbf{F}(\mathcal{X}^{\ast })$ is a weak$^{\ast }$%
-Hausdorff-clopen subset. This is not possible if $\mathbf{B}(\mathcal{X}%
^{\ast })$\ is not weak$^{\ast }$-Hausdorff-open, because $\mathbf{B}(%
\mathcal{X}^{\ast })$\ is a connected component, by Proposition \ref%
{connected sub-hyperspace}.

So, it remains to prove that $\mathbf{B}(\mathcal{X}^{\ast })$\ is not weak$%
^{\ast }$-Hausdorff-open when $\mathcal{X}$ is infinite-dimensional. To this
end, assume by contradiction that $\mathbf{B}(\mathcal{X}^{\ast })$\ is weak$%
^{\ast }$-Hausdorff-open. Then, for any $B\in \mathbf{B}(\mathcal{X}^{\ast
}) $, there exists $\left( A_{k}\right) _{k=1}^{n}\subseteq \mathcal{X}$ and
$\varepsilon \in \mathbb{R}^{+}$ such that%
\begin{equation}
\left\{ F\in \mathbf{F}\left( \mathcal{X}^{\ast }\right) :\forall k\in
\left\{ 1,\ldots ,n\right\} ,\ d_{H}^{(A_{k})}\left( F,B\right) <\varepsilon
\right\} \subseteq \mathbf{B}\left( \mathcal{X}^{\ast }\right) \ .
\label{qq1}
\end{equation}%
Take any
\begin{equation}
\sigma \in \bigcap\limits_{k=1}^{n}\mathrm{ker}(\hat{A}_{k})\backslash
\left\{ 0\right\}  \label{qq2}
\end{equation}%
with $\hat{A}$ being defined by (\ref{sdfsdfkljsdlfkj}) for any $A\in
\mathcal{X}$. Such an element always exists because $\mathcal{X}^{\ast }$ is
infinite-dimensional, if $\mathcal{X}$ is infinite-dimensional. By (\ref{qq2}%
), for any $k\in \left\{ 1,\ldots ,n\right\} $,
\begin{equation*}
d_{H}^{(A_{k})}\left( \overline{B+\mathbb{R}\sigma },B\right) =0
\end{equation*}%
and (\ref{qq1}) yields $\overline{B+\mathbb{R}\sigma }\in \mathbf{B}(%
\mathcal{X}^{\ast })$, which is clearly not possible since $\sigma \neq 0$.
Consequently, $\mathbf{B}(\mathcal{X}^{\ast })$\ is not weak$^{\ast }$%
-Hausdorff-open when $\mathcal{X}$ is infinite-dimensional.
\end{proof}

In contrast with $\mathbf{B}(\mathcal{X}^{\ast })$, within the set $\mathbf{B%
}^{c}(\mathcal{X}^{\ast })$ there are possibly many connected components of $%
\mathbf{F}(\mathcal{X}^{\ast })$, as one can see from Proposition \ref%
{convexity corrolary copy(1)}.

Finally, for any $0$-neighborhood $\mathcal{U}\in \mathcal{O}$, note that
the hyperspace $\mathbf{U}_{\mathcal{U}}(\mathcal{X}^{\ast })$ of all
nonempty, uniformly bounded in $\mathcal{U}$, weak$^{\ast }$-closed subsets
of $\mathcal{X}^{\ast }$ defined by (\ref{U0}) has topological properties
that are similar to the set $\mathbf{B}(\mathcal{X}^{\ast })$:

\begin{lemma}[Hyperconvergence of uniformly bounded near $0$ sets]
\label{lemma close unif boundedness}\mbox{ }\newline
Let $\mathcal{X}$ be a topological $\mathbb{K}$-vector space with $\mathbb{K}%
=\mathbb{R},\mathbb{C}$. Take any weak$^{\ast }$-Hausdorff convergent net $%
(U_{j})_{j\in J}\subseteq \mathbf{U}(\mathcal{X}^{\ast })$ with $U_{j}\in
\mathbf{U}_{\mathcal{U}_{j}}\left( \mathcal{X}^{\ast }\right) $ and $%
\mathcal{U}_{j}\in \mathcal{O}$ for $j\in J$. If
\begin{equation}
\mathcal{U}_{\infty }\doteq \bigcap\limits_{j\in J}\mathcal{U}_{j}\in
\mathcal{O}  \label{defined Uinfinity}
\end{equation}%
then $(U_{j})_{j\in J}$ converges to $U_{\infty }\in \mathbf{U}_{\mathcal{U}%
_{\infty }}(\mathcal{X}^{\ast })\subseteq \mathbf{U}(\mathcal{X}^{\ast })$.
\end{lemma}

\begin{proof}
Take any weak$^{\ast }$-Hausdorff convergent net $(U_{j})_{j\in J}\subseteq
\mathbf{U}(\mathcal{X}^{\ast })$, as stated in the lemma. Assume that the
limit $U_{\infty }\notin \mathbf{U}_{\mathcal{U}_{\infty }}(\mathcal{X}%
^{\ast })$ with $\mathcal{U}_{\infty }\in \mathcal{O}$ defined by (\ref%
{defined Uinfinity}). Then, there is $\sigma _{\infty }\in U_{\infty }$ and $%
A\in \mathcal{U}_{\infty }$ such that $\left\vert \sigma _{\infty }\left(
A\right) \right\vert >1$ and so,
\begin{equation*}
d_{H}^{(A)}\left( U_{\infty },U_{j}\right) \geq \inf_{\tilde{\sigma}\in
U_{j}}\left\vert \left( \sigma _{\infty }-\tilde{\sigma}\right) \left(
A\right) \right\vert \geq \left\vert \sigma _{\infty }\left( A\right)
\right\vert -\sup_{\tilde{\sigma}\in U_{j}}\left\vert \tilde{\sigma}\left(
A\right) \right\vert >0\ ,
\end{equation*}%
by\ Definition \ref{hypertopology0} and the triangle inequality. But this
contradicts the fact that $(U_{j})_{j\in J}$ converges to $U_{\infty }$ in
the weak$^{\ast }$-Hausdorff hypertopology. Therefore, $U_{\infty }\in
\mathbf{U}_{\mathcal{U}_{\infty }}(\mathcal{X}^{\ast })$.
\end{proof}

\begin{corollary}[Weak$^{\ast }$-Hausdorff-closed subhyperspaces $\mathbf{U}%
_{\mathcal{U}}(\mathcal{X}^{\ast })$]
\label{lemma close unif boundedness-II}\mbox{ }\newline
Let $\mathcal{X}$ be a topological $\mathbb{K}$-vector space with $\mathbb{K}%
=\mathbb{R},\mathbb{C}$. Then, for any $\mathcal{U}\in \mathcal{O}$, $%
\mathbf{U}_{\mathcal{U}}(\mathcal{X}^{\ast })$ is a convex, path-connected,
weak$^{\ast }$-Hausdorff-closed subset of $\mathbf{U}\left( \mathcal{X}%
^{\ast }\right) \subseteq \mathbf{K}(\mathcal{X}^{\ast })$.
\end{corollary}

\begin{proof}
Convexity is obvious and path connectedness is proven by using the function (%
\ref{function}), observing that weak$^{\ast }$-compact subsets of $\mathcal{X%
}^{\ast }$ are bounded, by \cite[Theorem 1.15 (b)]{Rudin}. By Lemma \ref%
{lemma close unif boundedness}, $\mathbf{U}_{\mathcal{U}}(\mathcal{X}^{\ast
})$ is weak$^{\ast }$-Hausdorff-closed for any fixed $\mathcal{U}\in
\mathcal{O}$.
\end{proof}

\subsection{Hausdorff Property and Convexity}

One fundamental question one shall ask regarding the hyperspace $\mathbf{F}(%
\mathcal{X}^{\ast })$ is whether it is a Hausdorff space, with respect to
the weak$^{\ast }$-Hausdorff hypertopology, or not. The answer is \emph{%
negative} for real Banach spaces of dimension greater than 1, as
demonstrated in the next lemma by using elements of the set $\mathbf{K}(%
\mathcal{X}^{\ast })$ of all nonempty weak$^{\ast }$-compact sets defined by
(\ref{hyperspace K}):

\begin{lemma}[Non-weak$^{\ast }$-Hausdorff-separable points]
\label{convexity lemma copy(2)}\mbox{ }\newline
Let $\mathcal{X}$ be a topological $\mathbb{R}$-vector space. Take any
convex weak$^{\ast }$-compact set $K\in \mathbf{CK}(\mathcal{X}^{\ast })$
with weak$^{\ast }$-path-connected weak$^{\ast }$-closed set $\mathcal{E}%
(K)\subseteq K$ of extreme points\footnote{%
Cf. the Krein-Milman theorem \cite[Theorem 3.23]{Rudin}.}. Then, $\mathcal{E}%
(K)\in \mathbf{K}(\mathcal{X}^{\ast })$ and $d_{H}^{(A)}(K,\mathcal{E}(K))=0$
for any $A\in \mathcal{X}$.
\end{lemma}

\begin{proof}
Let $\mathcal{X}$ be a topological $\mathbb{R}$-vector space. Recall that
any $A\in \mathcal{X}$ defines a weak$^{\ast }$-continuous linear functional
$\hat{A}:\mathcal{X}^{\ast }\rightarrow \mathbb{R}$, by Equation (\ref%
{sdfsdfkljsdlfkj}). Observe next that
\begin{equation}
d_{H}^{(A)}(K,\mathcal{E}(K))=\max \left\{ \max_{x_{1}\in \hat{A}\left(
K\right) }\min_{x_{2}\in \hat{A}\left( \mathcal{E}(K)\right) }\left\vert
x_{1}-x_{2}\right\vert ,\max_{x_{2}\in \hat{A}\left( \mathcal{E}(K)\right)
}\min_{x_{1}\in \hat{A}\left( K\right) }\left\vert x_{1}-x_{2}\right\vert
\right\} \ .  \label{rewritte0}
\end{equation}%
We obviously have the inclusions
\begin{equation}
\hat{A}\left( \mathcal{E}\left( K\right) \right) \subseteq \hat{A}\left(
K\right) \subseteq \left[ \min \hat{A}\left( K\right) ,\max \hat{A}\left(
K\right) \right] \ .  \label{rewritte}
\end{equation}%
By the Bauer maximum principle \cite[Lemma 10.31]{BruPedra2} together with
the affinity and weak$^{\ast }$-continuity of $\hat{A}$,
\begin{equation*}
\min \hat{A}\left( K\right) =\min \hat{A}\left( \mathcal{E}\left( K\right)
\right) \qquad \text{and}\qquad \max \hat{A}\left( K\right) =\max \hat{A}%
\left( \mathcal{E}\left( K\right) \right) \ .
\end{equation*}%
In particular, we can rewrite (\ref{rewritte}) as
\begin{equation}
\hat{A}\left( \mathcal{E}(K)\right) \subseteq \hat{A}\left( K\right)
\subseteq \left[ \min \hat{A}\left( \mathcal{E}\left( K\right) \right) ,\max
\hat{A}\left( \mathcal{E}\left( K\right) \right) \right] \ .
\label{rewritte2}
\end{equation}%
Since $\mathcal{E}(K)$ is, by assumption, path-connected in the weak$^{\ast
} $ topology, there is a weak$^{\ast }$-continuous path $\gamma
:[0,1]\rightarrow \mathcal{E}(K)$ from a minimizer to a maximizer of $\hat{A}
$ in $\mathcal{E}(K)$. By weak$^{\ast }$-continuity of $\hat{A}$, it follows
that
\begin{equation*}
\left[ \min \hat{A}\left( \mathcal{E}(K)\right) ,\max \hat{A}\left( \mathcal{%
E}(K)\right) \right] =\hat{A}\circ \gamma \left( \left[ 0,1\right] \right)
\subseteq \hat{A}\left( \mathcal{E}(K)\right)
\end{equation*}%
and we infer from (\ref{rewritte2}) that
\begin{equation*}
\hat{A}\left( \mathcal{E}(K)\right) =\hat{A}\left( K\right) =\left[ \min
\hat{A}\left( K\right) ,\max \hat{A}\left( K\right) \right] =\left[ \min
\hat{A}\left( \mathcal{E}(K)\right) ,\max \hat{A}\left( \mathcal{E}%
(K)\right) \right] \ .
\end{equation*}%
Together with (\ref{rewritte0}), this last equality obviously leads to the
assertion. Note that $\mathcal{E}(K)\in \mathbf{K}(\mathcal{X}^{\ast })$
since it is, by assumption, a weak$^{\ast }$-closed subset of the weak$%
^{\ast }$-compact set $K$.
\end{proof}

\begin{corollary}[Non-Hausdorff hyperspaces]
\label{Non-Hausdorff hyperspaces}\mbox{ }\newline
Let $\mathcal{X}$ be a topological $\mathbb{R}$-vector space of dimension
greater than 1. Then, $\mathbf{F}(\mathcal{X}^{\ast })$, $\mathbf{B}(%
\mathcal{X}^{\ast })$, $\mathbf{K}(\mathcal{X}^{\ast })$, $\mathbf{U}(%
\mathcal{X}^{\ast })$ are all non-Hausdorff spaces.
\end{corollary}

\begin{proof}
This corollary is a direct consequence of Lemma \ref{convexity lemma copy(2)}
by observing that the dual space of a topological $\mathbb{R}$-vector space
of dimension greater than 1 contains a two-dimensional closed disc.
\end{proof}

In fact, the weak$^{\ast }$-Hausdorff hypertopology cannot distinguish a set
from its weak$^{\ast }$-closed convex hull, as it also appears for other
well-established hypertopologies, like the so-called scalar topology for
closed sets (see \cite[Section 4.3]{Beer}). Similar to the scalar topology,
only $\mathbf{CF}(\mathcal{X}^{\ast })$ is a Hausdorff hyperspace. To get an
intuition of this, consider the following result:

\begin{proposition}[Separation of the weak$^{\ast }$-closed convex hull]
\label{convexity lemma copy(1)}\mbox{ }\newline
Let $\mathcal{X}$ be a topological $\mathbb{K}$-vector space with $\mathbb{K}%
=\mathbb{R},\mathbb{C}$. Take $F_{1},F_{2}\in \mathbf{F}(\mathcal{X}^{\ast
}) $. If $d_{H}^{(A)}(F_{1},F_{2})=0$ for all $A\in \mathcal{X}$, then $%
\overline{\mathrm{co}}F_{1}=\overline{\mathrm{co}}F_{2}$, where $\overline{%
\mathrm{co}}$ is the weak$^{\ast }$-closed convex hull operator defined by (%
\ref{closure operator}).
\end{proposition}

\begin{proof}
Pick any weak$^{\ast }$-closed sets $F_{1},F_{2}$ satisfying $%
d_{H}^{(A)}(F_{1},F_{2})=0$ for all $A\in \mathcal{X}$. Let $\sigma _{1}\in
F_{1}$. By\ Definition \ref{hypertopology0}, it follows that%
\begin{equation}
\inf_{\sigma _{2}\in F_{2}}\left\vert \left( \sigma _{1}-\sigma _{2}\right)
\left( A\right) \right\vert =0\ ,\qquad A\in \mathcal{X}\ .
\label{contradition00}
\end{equation}%
Recall that the dual space $\mathcal{X}^{\ast }$ of $\mathcal{X}$ is a
locally convex (Hausdorff) space in the weak$^{\ast }$ topology and its dual
space is $\mathcal{X}$. By (\ref{closure operator}), $\overline{\mathrm{co}}%
F_{2}$ is convex and weak$^{\ast }$-closed and $\{\sigma _{1}\}$ is a convex
weak$^{\ast }$-compact set. If $\sigma _{1}\notin \overline{\mathrm{co}}%
F_{2} $ then we infer from the Hahn-Banach separation theorem \cite[Theorem
3.4 (b)]{Rudin} the existence of $A_{0}\in \mathcal{X}$ and $x_{1},x_{2}\in
\mathbb{R}$ such that
\begin{equation}
\sup_{\sigma _{2}\in \overline{\mathrm{co}}F_{2}}\mathrm{Re}\left\{ \sigma
_{2}\left( A_{0}\right) \right\} <x_{1}<x_{2}<\mathrm{Re}\left\{ \sigma
_{1}\left( A_{0}\right) \right\} \ ,  \label{hahn banach}
\end{equation}%
which contradicts (\ref{contradition00}) for $A=A_{0}$. As a consequence, $%
\sigma _{1}\in \overline{\mathrm{co}}F_{2}$ and hence, $F_{1}\subseteq
\overline{\mathrm{co}}F_{2}$. This in turn yields $\overline{\mathrm{co}}%
F_{1}\subseteq \overline{\mathrm{co}}F_{2}$. By switching the role of the
weak$^{\ast }$-closed sets, we thus deduce the assertion.
\end{proof}

\begin{corollary}[$\mathbf{CF}(\mathcal{X}^{\ast })$ as an Hausdorff
hyperspace]
\label{convexity corrolary}\mbox{ }\newline
Let $\mathcal{X}$ be a topological $\mathbb{K}$-vector space with $\mathbb{K}%
=\mathbb{R},\mathbb{C}$. Then, $\mathbf{CF}(\mathcal{X}^{\ast })$ is a
Hausdorff hyperspace.
\end{corollary}

\begin{proof}
This is a direct consequence of Proposition \ref{convexity lemma copy(1)}.
\end{proof}

Note that the weak$^{\ast }$-closed convex hull operator $\overline{\mathrm{%
co}}$ is a weak$^{\ast }$-Hausdorff continuous mapping:

\begin{proposition}[Weak$^{\ast }$-Hausdorff continuity of the weak$^{\ast }$%
-closed convex hull operator]
\label{convexity lemma}\mbox{ }\newline
Let $\mathcal{X}$ be a topological $\mathbb{K}$-vector space with $\mathbb{K}%
=\mathbb{R},\mathbb{C}$. Then, $\overline{\mathrm{co}}$ is a weak$^{\ast }$%
-Hausdorff continuous mapping from $\mathbf{F}(\mathcal{X}^{\ast })$ onto $%
\mathbf{CF}(\mathcal{X}^{\ast })$.
\end{proposition}

\begin{proof}
Let $\mathcal{X}$ be a topological $\mathbb{K}$-vector space with $\mathbb{K}%
=\mathbb{R},\mathbb{C}$. The surjectivity of $\overline{\mathrm{co}}$ seen
as a mapping from $\mathbf{F}(\mathcal{X}^{\ast })$ to $\mathbf{CF}(\mathcal{%
X}^{\ast })$ is obvious, by (\ref{surjective}). Now, take any weak$^{\ast }$%
-Hausdorff convergent net $(F_{j})_{j\in J}\subseteq \mathbf{F}(\mathcal{X}%
^{\ast })$ with limit $F_{\infty }\in \mathbf{F}(\mathcal{X}^{\ast })$. Note
that
\begin{equation}
\sup_{\sigma \in \overline{\mathrm{co}}\left( F_{\infty }\right) }\inf_{%
\tilde{\sigma}\in \overline{\mathrm{co}}\left( F_{j}\right) }\left\vert
\left( \sigma -\tilde{\sigma}\right) \left( A\right) \right\vert
=\sup_{\sigma \in \mathrm{co}\left( F_{\infty }\right) }\inf_{\tilde{\sigma}%
\in \overline{\mathrm{co}}\left( F_{j}\right) }\left\vert \left( \sigma -%
\tilde{\sigma}\right) \left( A\right) \right\vert \ ,\qquad A\in \mathcal{X}%
\ ,  \label{dd}
\end{equation}%
because, for any $A\in \mathcal{X}$, $j\in J$, $\sigma _{1},\sigma _{2}\in
\overline{\mathrm{co}}\left( F_{\infty }\right) $ and $\tilde{\sigma}\in
\overline{\mathrm{co}}\left( F_{j}\right) $,%
\begin{equation*}
\left\vert \left\vert \left( \sigma _{1}-\tilde{\sigma}\right) \left(
A\right) \right\vert -\left\vert \left( \sigma _{2}-\tilde{\sigma}\right)
\left( A\right) \right\vert \right\vert \leq \left\vert \left( \sigma
_{1}-\sigma _{2}\right) \left( A\right) \right\vert \ ,
\end{equation*}%
which yields%
\begin{equation*}
\left\vert \inf_{\tilde{\sigma}\in \overline{\mathrm{co}}\left( F_{j}\right)
}\left\vert \left( \sigma _{1}-\tilde{\sigma}\right) \left( A\right)
\right\vert -\inf_{\tilde{\sigma}\in \overline{\mathrm{co}}\left(
F_{j}\right) }\left\vert \left( \sigma _{2}-\tilde{\sigma}\right) \left(
A\right) \right\vert \right\vert \leq \left\vert \left( \sigma _{1}-\sigma
_{2}\right) \left( A\right) \right\vert
\end{equation*}%
for any $A\in \mathcal{X}$, $j\in J$ and $\sigma _{1},\sigma _{2}\in
\overline{\mathrm{co}}\left( F_{\infty }\right) $. Fix $n\in \mathbb{N}$, $%
\sigma _{1},\ldots ,\sigma _{n}\in F_{\infty }$ and parameters $\lambda
_{1},\ldots ,\lambda _{n}\in \left[ 0,1\right] $ such that
\begin{equation*}
\sum\limits_{k=1}^{n}\lambda _{k}=1\ .
\end{equation*}%
Pick any parameter $\varepsilon \in \mathbb{R}^{+}$. For any $A\in \mathcal{X%
}$ and $k\in \{1,\ldots ,n\}$, we define $\tilde{\sigma}_{k,j}\in F_{j}$
such that%
\begin{equation*}
\left\vert \left( \sigma _{k}-\tilde{\sigma}_{k,j}\right) \left( A\right)
\right\vert \leq \inf_{\tilde{\sigma}\in F_{j}}\left\vert \left( \sigma _{k}-%
\tilde{\sigma}\right) \left( A\right) \right\vert +\varepsilon \ .
\end{equation*}%
Then, for all $j\in J$ and $\varepsilon \in \mathbb{R}^{+}$,%
\begin{equation*}
\inf_{\tilde{\sigma}\in \overline{\mathrm{co}}\left( F_{j}\right)
}\left\vert \left( \sum\limits_{k=1}^{n}\lambda _{k}\sigma _{k}-\tilde{\sigma%
}\right) \left( A\right) \right\vert \leq \sum\limits_{k=1}^{n}\lambda
_{k}\left\vert \left( \sigma _{k}-\tilde{\sigma}_{k,j}\right) \left(
A\right) \right\vert \leq \varepsilon +\sup_{\sigma \in F_{\infty }}\inf_{%
\tilde{\sigma}\in F_{j}}\left\vert \left( \sigma -\tilde{\sigma}\right)
\left( A\right) \right\vert \ .
\end{equation*}%
Using (\ref{dd}), we thus deduce that, for all $j\in J$,
\begin{equation}
\sup_{\sigma \in \overline{\mathrm{co}}\left( F_{\infty }\right) }\inf_{%
\tilde{\sigma}\in \overline{\mathrm{co}}\left( F_{j}\right) }\left\vert
\left( \sigma -\tilde{\sigma}\right) \left( A\right) \right\vert \leq
\sup_{\sigma \in F_{\infty }}\inf_{\tilde{\sigma}\in F_{j}}\left\vert \left(
\sigma -\tilde{\sigma}\right) \left( A\right) \right\vert \ ,\qquad A\in
\mathcal{X}\ .  \label{convex1}
\end{equation}%
By switching the role of $F_{\infty }$ and $F_{j}$ for every $j\in J$, we
also arrive at the inequality
\begin{equation}
\sup_{\tilde{\sigma}\in \overline{\mathrm{co}}\left( F_{j}\right)
}\inf_{\sigma \in \overline{\mathrm{co}}\left( F_{\infty }\right)
}\left\vert \left( \sigma -\tilde{\sigma}\right) \left( A\right) \right\vert
\leq \sup_{\tilde{\sigma}\in F_{j}}\inf_{\sigma \in F_{\infty }}\left\vert
\left( \sigma -\tilde{\sigma}\right) \left( A\right) \right\vert \ ,\qquad
A\in \mathcal{X}\ .  \label{convex2}
\end{equation}%
Since $(F_{j})_{j\in J}$ converges in the weak$^{\ast }$-Hausdorff
hypertopology to $F_{\infty }$, Inequalities (\ref{convex1})-(\ref{convex2})
combined with Definition \ref{hypertopology0} yield the weak$^{\ast }$%
-Hausdorff convergence of $(\overline{\mathrm{co}}\left( F_{j}\right)
)_{j\in J}$ to $\overline{\mathrm{co}}\left( F_{\infty }\right) $. By \cite[%
Chapter 3, Theorem 1]{topology}, $\overline{\mathrm{co}}$ is a weak$^{\ast }$%
-Hausdorff continuous mapping onto $\mathbf{CF}(\mathcal{X}^{\ast })$.
\end{proof}

Proposition \ref{convexity lemma} has a direct consequence on the
topological properties of hyperspaces of convex weak$^{\ast }$-closed sets:

\begin{corollary}[Weak$^{\ast }$-Hausdorff-closed hyperspaces of convex sets]

\label{convexity corrolary copy(2)}\mbox{ }\newline
Let $\mathcal{X}$ be a topological $\mathbb{K}$-vector space with $\mathbb{K}%
=\mathbb{R},\mathbb{C}$. \newline
\emph{(i)} $\mathbf{CF}(\mathcal{X}^{\ast })=\overline{\mathrm{co}}(\mathbf{F%
}(\mathcal{X}^{\ast }))$ is a convex, weak$^{\ast }$-Hausdorff-closed subset
of $\mathbf{F}(\mathcal{X}^{\ast })$. \newline
\emph{(ii)} $\mathbf{CB}(\mathcal{X}^{\ast })\doteq \overline{\mathrm{co}}(%
\mathbf{B}(\mathcal{X}^{\ast }))$ is a convex, path-connected, weak$^{\ast }$%
-Hausdorff-closed subset of $\mathbf{CF}(\mathcal{X}^{\ast })$.\newline
\emph{(iii)} For any $\mathcal{U}\in \mathcal{O}$, $\mathbf{CU}_{\mathcal{U}%
}(\mathcal{X}^{\ast })$ is a convex, path-connected, weak$^{\ast }$%
-Hausdorff-closed subset of $\mathbf{CU}\left( \mathcal{X}^{\ast }\right)
\subseteq \mathbf{CK}(\mathcal{X}^{\ast })\subseteq \mathbf{CF}(\mathcal{X}%
^{\ast })$.
\end{corollary}

\begin{proof}
By Corollary \ref{convexity corrolary}, $\mathbf{CF}(\mathcal{X}^{\ast })$
endowed with the weak$^{\ast }$-Hausdorff hypertopology is a Hausdorff
space. Hence, by \cite[Chapter 2, Theorem 3]{topology}, each convergent net
in this space converges in the weak$^{\ast }$-Hausdorff hypertopology to at
most one point, which, by Proposition \ref{convexity lemma}, must be a
convex weak$^{\ast }$-closed set. Assertion (i) is thus proven. Convexity of
$\mathbf{CF}(\mathcal{X}^{\ast })$ is obvious.

To prove (ii), recall that $\mathbf{B}(\mathcal{X}^{\ast })$ is a weak$%
^{\ast }$-Hausdorff-closed set, by Proposition \ref{convexity corrolary
copy(1)}. Using this together with Proposition \ref{convexity lemma} and (%
\ref{CB}), we deduce that $\mathbf{CB}(\mathcal{X}^{\ast })$ is also weak$%
^{\ast }$-Hausdorff-closed. By Propositions \ref{connected sub-hyperspace}
and \ref{convexity lemma} and the fact that the image under a continuous
mapping of a path-connected space is path-connected, $\mathbf{CB}(\mathcal{X}%
^{\ast })$ is also path-connected. Convexity of $\mathbf{CB}(\mathcal{X}%
^{\ast })$ is obvious.

Assertion (iii) follows from Corollary \ref{lemma close unif boundedness-II}
and Proposition \ref{convexity lemma} together with (\ref{CU0}). In
particular, the convexity of $\mathbf{CU}_{\mathcal{U}}(\mathcal{X}^{\ast })$
is obvious. Recall also the Banach-Alaoglu theorem \cite[Theorem 3.15]{Rudin}%
, leading to $\mathbf{CU}\left( \mathcal{X}^{\ast }\right) \subseteq \mathbf{%
CK}(\mathcal{X}^{\ast })$.
\end{proof}

An extension of Corollary \ref{convexity corrolary copy(2)} (iii) to the set
$\mathbf{CK}(\mathcal{X}^{\ast })$ of all nonempty convex weak$^{\ast }$%
-compact sets, defined by (\ref{hyperspace CK}), is not a priori clear
because the weak$^{\ast }$-closed convex hull operator $\overline{\mathrm{co}%
}$ does not necessarily maps weak$^{\ast }$-compact sets to weak$^{\ast }$%
-compact sets for general (real or complex) topological vector space $%
\mathcal{X}$. We do not know a priori whether $\mathbf{K}(\mathcal{X}^{\ast
})$ and, hence, $\mathbf{CK}(\mathcal{X}^{\ast })$, are weak$^{\ast }$%
-Hausdorff-closed subset of $\mathbf{F}(\mathcal{X}^{\ast })$. This property
is at least true when $\mathcal{X}$ is a Banach space, since in this case, $%
\mathbf{B}(\mathcal{X}^{\ast })=\mathbf{K}(\mathcal{X}^{\ast })$, by the
Banach-Alaoglu theorem \cite[Theorem 3.15]{Rudin} and the uniform
boundedness principle \cite[Theorems 2.4 and 2.5]{Rudin}. See, e.g., \cite[%
Proposition 1.2.9]{Beer}.

\subsection{Weak$^{\ast }$-Hausdorff Hyperconvergence\label%
{Hyperconvergences}}

It is instructive to relate weak$^{\ast }$-Hausdorff limits of nets to lower
and upper limits of sets \`{a} la Painlev\'{e} \cite[\S\ 29]%
{topology-painleve}: The \emph{lower limit} of any net $(F_{j})_{j\in J}$ of
subsets of $\mathcal{X}^{\ast }$ is defined by
\begin{equation}
\mathrm{Li}_{j\in J}F_{j}\doteq \left\{ \sigma \in \mathcal{X}^{\ast
}:\sigma \text{ is a weak}^{\ast }\text{ limit of a net }(\sigma _{j})_{j\in
J}\text{ with }\sigma _{j}\in F_{j}\text{, }j\succ j_{0}\text{, for some }%
j_{0}\in J\right\} \ ,  \label{Li}
\end{equation}%
while its \emph{upper limit} equals
\begin{equation}
\mathrm{Ls}_{j\in J}F_{j}\doteq \left\{ \sigma \in \mathcal{X}^{\ast }:\text{%
there is a subnet }(F_{s(i)})_{i\in I}\text{ and a net }\sigma
_{i}\rightarrow \sigma \text{ with }\sigma _{i}\in F_{s(i)}\text{ for all }%
i\in I\right\} \ .  \label{Ls}
\end{equation}%
Clearly, $\mathrm{Li}_{j\in J}F_{j}\subseteq \mathrm{Ls}_{j\in J}F_{j}$. If $%
\mathrm{Li}_{j\in J}F_{j}=\mathrm{Ls}_{j\in J}F_{j}$ then $\left(
F_{j}\right) _{j\in J}$ is said to be convergent to this set. See \cite[\S\ %
29, I, III, VI]{topology-painleve}, which however defines $\mathrm{Li}$ and $%
\mathrm{Ls}$ within metric spaces. This refers in the literature to the
\emph{Kuratowski} or \emph{Kuratowski-Painlev\'{e}}\footnote{%
The idea of upper and lower limits is due to Painlev\'{e}, as acknowledged
by Kuratowski himself in \cite[\S\ 29, Footnote 1, p. 335]{topology-painleve}%
. We thus use the name Kuratowski-Painlev\'{e} convergence.} convergence,
see e.g. \cite[Appendix B]{Lucchetti} and \cite[Section 5.2]{Beer}. By \cite[%
Theorem 1.22]{Rudin}, if $\mathcal{X}$ is an infinite-dimensional space,
then its dual space $\mathcal{X}^{\ast }$ is not locally compact. In this
case, the Kuratowski-Painlev\'{e} convergence is not topological \cite[%
Theorem B.3.2]{Lucchetti}. See also \cite[Chapter 5]{Beer}, in particular
\cite[Theorem 5.2.6 and following discussions]{Beer} which relates the
Kuratowski-Painlev\'{e} convergence to the so-called Fell\emph{\ }topology.

We start by proving the weak$^{\ast }$-Hausdorff convergence of
monotonically increasing nets which are bounded from above within the
subspace $\mathbf{K}(\mathcal{X}^{\ast })$ of all nonempty weak$^{\ast }$%
-compact subsets of $\mathcal{X}^{\ast }$ defined by (\ref{hyperspace K}).

\begin{proposition}[Weak$^{\ast }$-Hausdorff hyperconvergence of increasing
nets]
\label{Solution selfbaby copy(5)+1}\mbox{ }\newline
Let $\mathcal{X}$ be a topological $\mathbb{K}$-vector space with $\mathbb{K}%
=\mathbb{R},\mathbb{C}$. Any increasing net $(K_{j})_{j\in J}\subseteq
\mathbf{K}(\mathcal{X}^{\ast })$ such that
\begin{equation}
K\doteq \overline{\bigcup\limits_{j\in J}K_{j}}\in \mathbf{K}\left( \mathcal{%
X}^{\ast }\right) \varsubsetneq \mathbf{F}\left( \mathcal{X}^{\ast }\right)
\label{sdsdsdsd}
\end{equation}%
(with respect to the weak$^{\ast }$ closure) converges in the weak$^{\ast }$%
-Hausdorff hypertopology to
\begin{equation}
K=\overline{\mathrm{Li}_{j\in J}K_{j}}=\overline{\mathrm{Ls}_{j\in J}K_{j}}\
.  \label{sdsdsdsdsdsdsdsd}
\end{equation}%
Additionally, $K$ is the Kuratowski-Painlev\'{e} limit of $(K_{j})_{j\in J}$
whenever $\mathcal{X}$ is separable.
\end{proposition}

\begin{proof}
Let $(K_{j})_{j\in J}\subseteq \mathbf{K}(\mathcal{X}^{\ast })$ be any
increasing net, i.e., $K_{j_{1}}\subseteq K_{j_{2}}$ whenever $j_{1}\prec
j_{2}$, satisfying (\ref{sdsdsdsd}). Assume without loss of generality that $%
K_{j}\neq \emptyset $ for all $j\in J$. Because $K\in \mathbf{K}(\mathcal{X}%
^{\ast })$, it is bounded, see \cite[Theorem 1.15 (b)]{Rudin}. By the
convergence of increasing bounded nets of real numbers, it follows that, for
any $A\in \mathcal{X}$,
\begin{equation*}
\lim_{J}\max_{\tilde{\sigma}\in K_{j}}\min_{\sigma \in K}\left\vert \left(
\tilde{\sigma}-\sigma \right) \left( A\right) \right\vert =\sup_{j\in
J}\max_{\tilde{\sigma}\in K_{j}}\min_{\sigma \in K}\left\vert \left( \tilde{%
\sigma}-\sigma \right) \left( A\right) \right\vert \leq \max_{\tilde{\sigma}%
\in K}\min_{\sigma \in K}\left\vert \left( \tilde{\sigma}-\sigma \right)
\left( A\right) \right\vert =0\ .
\end{equation*}%
Therefore, by Definition \ref{hypertopology0}, if
\begin{equation}
\limsup_{J}\max_{\sigma \in K}\min_{\tilde{\sigma}\in K_{j}}\left\vert
\left( \tilde{\sigma}-\sigma \right) \left( A\right) \right\vert =0\ ,\qquad
A\in \mathcal{X}\ ,  \label{limit prove}
\end{equation}%
then the increasing net $(K_{j})_{j\in J}$ converges in $\mathbf{K}(\mathcal{%
X}^{\ast })$ to $K$. To prove (\ref{limit prove}), assume by contradiction
the existence of $\varepsilon \in \mathbb{R}^{+}$ such that%
\begin{equation}
\limsup_{J}\max_{\sigma \in K}\min_{\tilde{\sigma}\in K_{j}}\left\vert
\left( \tilde{\sigma}-\sigma \right) \left( A\right) \right\vert \geq
\varepsilon \in \mathbb{R}^{+}  \label{contradiction}
\end{equation}%
for some fixed $A\in \mathcal{X}$. For any $j\in J$, take $\sigma _{j}\in K$
such that%
\begin{equation}
\max_{\sigma \in K}\min_{\tilde{\sigma}\in K_{j}}\left\vert \left( \tilde{%
\sigma}-\sigma \right) \left( A\right) \right\vert =\min_{\tilde{\sigma}\in
K_{j}}\left\vert \left( \tilde{\sigma}-\sigma _{j}\right) \left( A\right)
\right\vert \ .  \label{to prove2}
\end{equation}%
By weak$^{\ast }$-compactness of $K$, there is a subnet $(\sigma
_{j_{l}})_{l\in L}$ converging in the weak$^{\ast }$ topology to $\sigma
_{\infty }\in K$. Via Equation (\ref{to prove2}) and the triangle
inequality, we then get that, for any $l\in L$,
\begin{equation*}
\max_{\sigma \in K}\min_{\tilde{\sigma}\in K_{j_{l}}}\left\vert \left(
\tilde{\sigma}-\sigma \right) \left( A\right) \right\vert \leq \left\vert
\left( \sigma _{j_{l}}-\sigma _{\infty }\right) \left( A\right) \right\vert
+\min_{\tilde{\sigma}\in K_{j_{l}}}\left\vert \left( \tilde{\sigma}-\sigma
_{\infty }\right) \left( A\right) \right\vert \ .
\end{equation*}%
By (\ref{sdsdsdsd}) and the fact that $(K_{j})_{j\in J}\subseteq \mathbf{K}(%
\mathcal{X}^{\ast })$ is an increasing net, it follows that%
\begin{equation}
\lim_{L}\max_{\sigma \in K}\min_{\tilde{\sigma}\in K_{j_{l}}}\left\vert
\left( \tilde{\sigma}-\sigma \right) \left( A\right) \right\vert =0\ .
\label{contradiction0}
\end{equation}%
By the convergence of decreasing bounded nets of real numbers, note that
\begin{equation*}
\limsup_{J}\max_{\sigma \in K}\min_{\tilde{\sigma}\in K_{j}}\left\vert
\left( \tilde{\sigma}-\sigma \right) \left( A\right) \right\vert
=\liminf_{J}\max_{\sigma \in K}\min_{\tilde{\sigma}\in K_{j}}\left\vert
\left( \tilde{\sigma}-\sigma \right) \left( A\right) \right\vert
\end{equation*}%
and hence, (\ref{contradiction0}) contradicts (\ref{contradiction}). As a
consequence, Equation (\ref{limit prove}) holds true.

In order to prove (\ref{sdsdsdsdsdsdsdsd}), observe that%
\begin{equation}
\bigcup\limits_{j\in J}K_{j}\subseteq \mathrm{Li}_{j\in J}K_{j}\subseteq
\mathrm{Ls}_{j\in J}K_{j}\subseteq \overline{\bigcup\limits_{j\in J}K_{j}}\ ,
\label{combined}
\end{equation}%
because the net $(K_{j})_{j\in J}$ is increasing. If $\mathcal{X}$ is
separable then it is well-known that the weak$^{\ast }$ topology of any
compact set is metrizable \cite[Theorem 3.16]{Rudin}, see Equation (\ref%
{metrics0}) below. Let $d$ be any metric generating the weak$^{\ast }$
topology on $K$. For any fixed $\sigma \in K$, consider the net $(\sigma
_{j})_{j\in J}$ where $\sigma _{j}$ is some minimizer in $K_{j}$ of $%
d(\sigma _{j},\sigma )$. Note that a minimizer always exists because $K_{j}$
is weak$^{\ast }$ compact. Clearly, this net converges in the weak$^{\ast }$
topology to $\sigma \in K$, for $(K_{j})_{j\in J}$ is an increasing net.
Hence, $K\subseteq \mathrm{Li}_{j\in J}K_{j}$ and therefore, using (\ref%
{combined}),%
\begin{equation*}
K=\mathrm{Li}_{j\in J}K_{j}=\mathrm{Ls}_{j\in J}K_{j}
\end{equation*}%
when $\mathcal{X}$ is separable.
\end{proof}

Non-monotonic, weak$^{\ast }$-Hausdorff convergent nets in $\mathbf{F}(%
\mathcal{X}^{\ast })$ are not trivial to study, in general. In the next
proposition, we give preliminary results on limits of convergent nets.

\begin{proposition}[Weak$^{\ast }$-Hausdorff hypertopology vs. upper and
lower limits]
\label{Solution selfbaby copy(5)+000}\mbox{ }\newline
Let $\mathcal{X}$ be a topological $\mathbb{K}$-vector space with $\mathbb{K}%
=\mathbb{R},\mathbb{C}$. For any weak$^{\ast }$-Hausdorff convergent net $%
(F_{j})_{j\in J}\subseteq \mathbf{F}(\mathcal{X}^{\ast })$ with limit $%
F_{\infty }\in \mathbf{F}(\mathcal{X}^{\ast })$,
\begin{equation*}
\mathrm{Li}_{j\in J}F_{j}\subseteq \overline{\mathrm{co}}\left( F_{\infty
}\right)
\end{equation*}%
and if $(F_{j}\equiv U_{j})_{j\in J}\subseteq \mathbf{U}_{\mathcal{U}}(%
\mathcal{X}^{\ast })$ for some $0$-neighborhood $\mathcal{U}\in \mathcal{O}$%
,
\begin{equation*}
F_{\infty }\equiv U_{\infty }\subseteq \overline{\mathrm{co}}\left( \mathrm{%
Ls}_{j\in J}U_{j}\right) \ ,
\end{equation*}%
where we recall that $\overline{\mathrm{co}}$ is the weak$^{\ast }$-closed
convex hull operator defined by (\ref{closure operator}).
\end{proposition}

\begin{proof}
Fix all parameters of the proposition. Assume without loss of generality
that $\mathrm{Li}_{j\in J}F_{j}$ is nonempty. Let $\sigma _{\infty }\in
\mathrm{Li}_{j\in J}F_{j}$, which is, by definition, the weak$^{\ast }$
limit of a net $(\sigma _{j})_{j\in J}$ such that $\sigma _{j}\in F_{j}$ for
all $j\succ j_{0}$, for some $j_{0}\in J$. Then, for any $A\in \mathcal{X}$
and $j\succ j_{0}$,
\begin{equation*}
\inf_{\sigma \in F_{\infty }}\left\vert \left( \sigma -\sigma _{\infty
}\right) \left( A\right) \right\vert \leq \left\vert \left( \sigma
_{j}-\sigma _{\infty }\right) \left( A\right) \right\vert +\inf_{\sigma \in
F_{\infty }}\left\{ \left\vert \left( \sigma -\sigma _{j}\right) \left(
A\right) \right\vert \right\} \ .
\end{equation*}%
Taking this last inequality in the limit with respect to $J$ and using
Definition \ref{hypertopology0}, we deduce that
\begin{equation}
\inf_{\sigma \in F_{\infty }}\left\vert \left( \sigma -\sigma _{\infty
}\right) \left( A\right) \right\vert =0\ ,\qquad A\in \mathcal{X}\ .
\label{contradition0}
\end{equation}%
If $\sigma _{\infty }\notin \overline{\mathrm{co}}\left( F_{\infty }\right) $
then, as it is done to prove (\ref{hahn banach}), we infer from the
Hahn-Banach separation theorem \cite[Theorem 3.4 (b)]{Rudin} the existence
of $A_{0}\in \mathcal{X}$ and $x_{1},x_{2}\in \mathbb{R}$ such that
\begin{equation*}
\sup_{\sigma \in \overline{\mathrm{co}}\left( F_{\infty }\right) }\mathrm{Re}%
\left\{ \sigma \left( A_{0}\right) \right\} <x_{1}<x_{2}<\mathrm{Re}\left\{
\sigma _{\infty }\left( A_{0}\right) \right\} \ ,
\end{equation*}%
which contradicts (\ref{contradition0}) for $A=A_{0}$. As a consequence, $%
\sigma _{\infty }\in \overline{\mathrm{co}}\left( F_{\infty }\right) $ and,
hence, $\mathrm{Li}_{j\in J}F_{j}\subseteq \overline{\mathrm{co}}\left(
F_{\infty }\right) $.

Assume now that $(F_{j}\equiv U_{j})_{j\in J}\subseteq \mathbf{U}_{\mathcal{U%
}}(\mathcal{X}^{\ast })$ for some $0$-neighborhood $\mathcal{U}\in \mathcal{O%
}$ with limit $F_{\infty }\equiv U_{\infty }$. By Definition \ref%
{hypertopology0} and the Banach-Alaoglu theorem \cite[Theorem 3.15]{Rudin},
we deduce that, for any $\sigma _{\infty }\in U_{\infty }$,
\begin{equation*}
\lim_{J}\inf_{\sigma \in U_{j}}\left\vert \left( \sigma -\sigma _{\infty
}\right) \left( A\right) \right\vert =\lim_{J}\min_{\sigma \in
U_{j}}\left\vert \left( \sigma -\sigma _{\infty }\right) \left( A\right)
\right\vert =0\ ,\qquad A\in \mathcal{X}\ .
\end{equation*}%
From this equality and the Banach-Alaoglu theorem \cite[Theorem 3.15]{Rudin}%
, for any $A\in \mathcal{X}$ and $\sigma _{\infty }\in U_{\infty }$, there
is $\tilde{\sigma}\in \mathrm{Ls}_{j\in J}U_{j}$ such that
\begin{equation*}
\tilde{\sigma}\left( A\right) =\sigma _{\infty }\left( A\right) \ .
\end{equation*}%
Consequently, one infers from the Hahn-Banach separation theorem \cite[%
Theorem 3.4 (b)]{Rudin} that $\sigma _{\infty }\in U_{\infty }$ belongs to
the weak$^{\ast }$-closed convex hull of the upper limit $\mathrm{Ls}_{j\in
J}U_{j}$.
\end{proof}

\begin{corollary}[Weak$^{\ast }$-Hausdorff hypertopology and convexity vs.
upper and lower limits]
\label{Solution selfbaby copy(2)}\mbox{ }\newline
Let $\mathcal{X}$ be a topological $\mathbb{K}$-vector space with $\mathbb{K}%
=\mathbb{R},\mathbb{C}$, $\mathcal{U}\in \mathcal{O}$ and $U_{\infty }\in
\mathbf{CU}_{\mathcal{U}}(\mathcal{X}^{\ast })$ be any weak$^{\ast }$%
-Hausdorff limit of a convergent net $(U_{j})_{j\in J}\subseteq \mathbf{CU}_{%
\mathcal{U}}(\mathcal{X}^{\ast })$. Then,
\begin{equation*}
\overline{\mathrm{Li}_{j\in J}U_{j}}=\overline{\mathrm{co}}\left( \mathrm{Li}%
_{j\in J}U_{j}\right) \subseteq U_{\infty }\subseteq \overline{\mathrm{co}}%
\left( \mathrm{Ls}_{j\in J}U_{j}\right) \ .
\end{equation*}
\end{corollary}

\begin{proof}
The assertion is an obvious application of Proposition \ref{Solution
selfbaby copy(5)+000} to the subset $\mathbf{CU}_{\mathcal{U}}(\mathcal{X}%
^{\ast })\subseteq \mathbf{U}_{\mathcal{U}}(\mathcal{X}^{\ast })$ together
with the idempotency of the weak$^{\ast }$-closed convex hull operator $%
\overline{\mathrm{co}}$. Note that $\mathrm{Li}_{j\in J}U_{j}$\ is in this
case a convex set.
\end{proof}

\subsection{Metrizable Hyperspaces}

We are interested in investigating \emph{metrizable} subspaces of $\mathbf{F}%
(\mathcal{X}^{\ast })$. Metrizable topological spaces are Hausdorff, so, in
the light of Corollaries \ref{Non-Hausdorff hyperspaces} and \ref{convexity
corrolary}, we restrict our analysis on subspaces of the Hausdorff
hyperspace $\mathbf{CF}(\mathcal{X}^{\ast })$ of all nonempty convex weak$%
^{\ast }$-closed subsets of $\mathcal{X}^{\ast }$ defined by Equation (\ref%
{hyperspace C}).

For a separable topological $\mathbb{K}$-vector space $\mathcal{X}$, recall
that the weak$^{\ast }$ topology on any compact set $K\in \mathbf{K}(%
\mathcal{X}^{\ast })$ is metrizable, see \cite[Theorem 3.16]{Rudin}. Here,
we use the following metric on $K$: Fix a countable dense set $(A_{n})_{n\in
\mathbb{N}}\subseteq \mathcal{X}$ and define%
\begin{equation}
d\left( \sigma _{1},\sigma _{2}\right) \doteq \sum\limits_{n\in \mathbb{N}}%
\frac{2^{-n}}{1+\max_{\sigma \in K}\left\vert \sigma \left( A_{n}\right)
\right\vert }\left\vert \left( \sigma _{1}-\sigma _{2}\right) \left(
A_{n}\right) \right\vert \ ,\qquad \sigma _{1},\sigma _{2}\in K\ .
\label{metrics0}
\end{equation}%
This metric is well-defined and induces the weak$^{\ast }$ topology on the
weak$^{\ast }$-compact set $K$. Absolute polars of $\mathcal{X}^{\ast }$
(cf. (\ref{polar})) are special example of compact sets, see \cite[Theorems\
3.15]{Rudin}. We show how (\ref{metrics0}) leads to the metrizability of the
weak$^{\ast }$-Hausdorff hypertopology on the hyperspace $\mathbf{CU}_{%
\mathcal{U}}(\mathcal{X}^{\ast })$ of all nonempty convex, uniformly bounded
in a $0$-neighborhood $\mathcal{U}\in \mathcal{O}$, weak$^{\ast }$-closed
subsets of $\mathcal{X}^{\ast }$.

Using the above metric $d$ in (\ref{Hausdorf}), denote by $\mathfrak{d}_{H}$
the Hausdorff distance between two elements $U_{1},U_{2}\in \mathbf{CU}_{%
\mathcal{U}}(\mathcal{X}^{\ast })$, that is\footnote{%
Minima in (\ref{metric1}) come from the compactness of sets and the
continuity of $d$. The following maxima in (\ref{metric1}) result from the
compactness of sets and the fact that the minimum over a continuous map
defines an upper semicontinuous function.},
\begin{equation}
\mathfrak{d}_{H}\left( U_{1},U_{2}\right) \doteq \max \left\{ \max_{\sigma
_{1}\in U_{1}}\min_{\sigma _{2}\in U_{2}}d\left( \sigma _{1},\sigma
_{2}\right) ,\max_{\sigma _{2}\in U_{2}}\min_{\sigma _{1}\in U_{1}}d\left(
\sigma _{1},\sigma _{2}\right) \right\} \ .  \label{metric1}
\end{equation}%
This Hausdorff distance induces the weak$^{\ast }$-Hausdorff hypertopology
on $\mathbf{CU}_{\mathcal{U}}(\mathcal{X}^{\ast })$:

\begin{theorem}[Complete metrizability of the weak$^{\ast }$-Hausdorff
hypertopology]
\label{Solution selfbaby copy(4)+1}\mbox{ }\newline
Let $\mathcal{X}$ be a separable topological $\mathbb{K}$-vector space with $%
\mathbb{K}=\mathbb{R},\mathbb{C}$ and $\mathcal{U}\in \mathcal{O}$. The
family
\begin{equation*}
\left\{ \left\{ U_{2}\in \mathbf{CU}_{\mathcal{U}}\left( \mathcal{X}^{\ast
}\right) :\mathfrak{d}_{H}\left( U_{1},U_{2}\right) <r\right\} :r\in \mathbb{%
R}^{\mathbb{+}},\ U_{1}\in \mathbf{CU}_{\mathcal{U}}\left( \mathcal{X}^{\ast
}\right) \right\}
\end{equation*}%
is a basis of the weak$^{\ast }$-Hausdorff hypertopology of $\mathbf{CU}_{%
\mathcal{U}}(\mathcal{X}^{\ast })$. Additionally, $\mathbf{CU}_{\mathcal{U}}(%
\mathcal{X}^{\ast })$ is weak$^{\ast }$-Hausdorff-compact and completely
metrizable.
\end{theorem}

\begin{proof}
Recall that a topology is finer than a second one iff any convergent net of
the first topology converges also in the second topology to the same limit.
See, e.g., \cite[Chapter 2, Theorems 4, 9]{topology}. We first show that the
topology induced by the Hausdorff metric $\mathfrak{d}_{H}$ is finer than
the weak$^{\ast }$-Hausdorff hypertopology of $\mathbf{CU}_{\mathcal{U}}(%
\mathcal{X}^{\ast })$ at fixed $\mathcal{U}\in \mathcal{O}$: Take any net $%
(U_{j})_{j\in J}$ converging in $\mathbf{CU}_{\mathcal{U}}(\mathcal{X}^{\ast
})$ to $U$ in the topology induced by the Hausdorff metric (\ref{metric1}).
Let $A\in \mathcal{X}$. By density of $(A_{n})_{n\in \mathbb{N}}$ in $%
\mathcal{X}$, for any $\varepsilon \in \mathbb{R}^{\mathbb{+}}$, there is $%
n\in \mathbb{N}$ such that $(A_{n}-A)\in 2^{-1}\varepsilon \mathcal{U}\in
\mathcal{O}$. In particular, by the definition of $\mathcal{U}^{\circ }$
(see (\ref{polar})) and (\ref{metrics0}), for all $j\in J$,
\begin{equation*}
d_{H}^{(A)}(U,U_{j})\leq \varepsilon +d_{H}^{(A_{n})}(U,U_{j})\leq
\varepsilon +2^{n}\left( 1+\max_{\sigma \in \mathcal{U}^{\circ }}\left\vert
\sigma \left( A_{n}\right) \right\vert \right) \mathfrak{d}_{H}(U,U_{j})\ .
\end{equation*}%
Thus, the net $(U_{j})_{j\in J}$ converges to $U$ also in the weak$^{\ast }$%
-Hausdorff hypertopology.

Endowed with the Hausdorff metric topology, the space of closed subsets of a
compact metric space is compact, by \cite[Theorem 3.2.4]{Beer}. In
particular, by weak$^{\ast }$ compactness of absolute polars (the
Banach-Alaoglu theorem \cite[Theorem 3.15]{Rudin}), $\mathbf{U}_{\mathcal{U}%
}(\mathcal{X}^{\ast })$ endowed with the Hausdorff metric $\mathfrak{d}_{H}$
is a compact hyperspace. By Corollary \ref{convexity corrolary copy(2)}, $%
\mathbf{CU}_{\mathcal{U}}(\mathcal{X}^{\ast })$ is closed with respect to
the weak$^{\ast }$-Hausdorff hypertopology, and thus closed with respect to
the topology induced by $\mathfrak{d}_{H}$, because this topology is finer
than the weak$^{\ast }$-Hausdorff hypertopology, as proven above. Hence, $%
\mathbf{CU}_{\mathcal{U}}(\mathcal{X}^{\ast })$ is also compact with respect
to the topology induced by $\mathfrak{d}_{H}$. Since the weak$^{\ast }$%
-Hausdorff hypertopology is a Hausdorff topology (Corollary \ref{convexity
corrolary}), as is well-known \cite[Section 3.8 (a)]{Rudin}, both topologies
must coincide.
\end{proof}

\noindent Note that Theorem \ref{Solution selfbaby copy(4)+1} is similar to
the assertion \cite[End of p. 91]{Beer}. It leads to a strong improvement of
Proposition \ref{Solution selfbaby copy(5)+1} and Corollary \ref{Solution
selfbaby copy(2)}:

\begin{corollary}[Weak$^{\ast }$-Hausdorff hypertopology and
Kuratowski-Painlev\'{e} convergence]
\label{Solution selfbaby copy(4)}\mbox{ }\newline
Let $\mathcal{X}$ be a separable topological $\mathbb{K}$-vector space with $%
\mathbb{K}=\mathbb{R},\mathbb{C}$ and $\mathcal{U}\in \mathcal{O}$. Then any
weak$^{\ast }$-Hausdorff convergent net $(U_{j})_{j\in J}\subseteq \mathbf{CU%
}_{\mathcal{U}}(\mathcal{X}^{\ast })$ converges to the Kuratowski-Painlev%
\'{e} limit%
\begin{equation*}
U_{\infty }=\mathrm{Li}_{j\in J}U_{j}=\mathrm{Ls}_{j\in J}U_{j}\in \mathbf{CU%
}_{\mathcal{U}}\left( \mathcal{X}^{\ast }\right) \ .
\end{equation*}
\end{corollary}

\begin{proof}
It is a direct consequence of Theorem \ref{Solution selfbaby copy(4)+1} and
\cite[\S\ 29, Section IX, Theorem 2]{topology-painleve}.
\end{proof}

\section{Generic Hypersets in Infinite Dimensions\label{Generic section}}

The Krein-Milman theorem \cite[Theorem 3.23]{Rudin} tells us that any convex
weak$^{\ast }$-compact set $K\in \mathbf{CK}(\mathcal{X}^{\ast })$ is the
weak$^{\ast }$-closure of the convex hull of the (nonempty) set $\mathcal{E}%
(K)$ of its extreme points:
\begin{equation*}
K=\overline{\mathrm{co}}\mathcal{E}\left( K\right) \ .
\end{equation*}%
The set $\mathcal{E}(K)$ is also called the extreme boundary of $K$. We are
interested in the question whether the subset of all $K\in \mathbf{CK}(%
\mathcal{X}^{\ast })$ with weak$^{\ast }$-dense set $\mathcal{E}\left(
K\right) $ of extreme points is generic, or not, when the topological space $%
\mathcal{X}$ has infinite dimension.

As is well-known, such convex compact sets exist in infinite-dimensional
topological spaces. For instance, the unit ball of any infinite-dimensional
Hilbert space has a dense extreme boundary in the weak topology. Another
example is given by the celebrated Poulsen simplex constructed in 1961 \cite%
{Poulsen}, within the Hilbert space $\ell ^{2}(\mathbb{N})$. In fact, a
convex compact set with dense extreme boundary \emph{is not an accident} in
this case: The set of all such convex compact subsets of an
infinite-dimensional separable\footnote{\cite[Proposition 2.1, Theorem 2.2]%
{Klee} seem to lead to the asserted property for all (possibly
non-separable) Banach spaces, as claimed in \cite%
{Klee,FonfLindenstrauss,infinite dim convexity}. However, \cite[Theorem 1.5]%
{Klee}, which assumes the separability of the Banach space, is clearly
invoked to prove the corresponding density stated in \cite[Theorem 2.2]{Klee}%
. We do not know how to remove the separability condition.} Banach space $%
\mathcal{Y}$ is generic\footnote{%
That is, the complement of a meagre set, i.e., a nowhere dense set.} in the
complete metric space of compact convex subsets of $\mathcal{Y}$, endowed
with the well-known Hausdorff metric topology \cite[Definition 3.2.1]{Beer}.
See \cite[Proposition 2.1, Theorem 2.2]{Klee}, which has been refined in
\cite[Section 4]{FonfLindenstrauss}. See, e.g., \cite[Section 7]{infinite
dim convexity} for a more recent review on this subject.

In this section we demonstrate the genericity in the dual space $\mathcal{X}%
^{\ast }$ of an infinite-dimensional, separable topological $\mathbb{K}$%
-vector space $\mathcal{X}$ ($\mathbb{K}=\mathbb{R},\mathbb{C}$), endowed
with its weak$^{\ast }$-topology. In this situation, results similar to \cite%
{Klee,FonfLindenstrauss} can be proven in a natural way by using the weak$%
^{\ast }$-Hausdorff hypertopology.

\subsection{Infinite Dimensionality of Absolute Polars}

A large class of convex weak$^{\ast }$-compact sets is given by absolute
polars (\ref{polar}) of any $0$-neighborhoods in $\mathcal{X}$, by the
Banach-Alaoglu theorem \cite[Theorem 3.15]{Rudin}. In fact, weak$^{\ast }$%
-closed subsets of absolute polars are the main sources of weak$^{\ast }$%
-compact sets in the dual space $\mathcal{X}^{\ast }$ of a real or complex
topological vector space $\mathcal{X}$. Therefore, it is natural to study
generic convex weak$^{\ast }$-compact sets within some absolute polar. If
the absolute polar can be embedded in a finite-dimensional subspace then
there is no convex weak$^{\ast }$-compact set with weak$^{\ast }$-dense
extreme boundary. We are thus interested in the infinite-dimensional
situation: we consider absolute polars which are \emph{infinite-dimensional}%
, that is, their (linear) spans are infinite-dimensional subspaces of $%
\mathcal{X}^{\ast }$.

Note that the infinite dimensionality of $\mathcal{X}$ does not guarantee
such a property of polars in $\mathcal{X}^{\ast }$: Take for instance $%
\mathcal{X}=\mathcal{H}$ being any infinite-dimensional Hilbert space
endowed with its weak topology and a scalar product denoted by $\left\langle
\cdot ,\cdot \right\rangle _{\mathcal{H}}$. For any $\mathcal{U}\in \mathcal{%
O}$, there are $n\in \mathbb{N}$ and $\psi _{1},\ldots ,\psi _{n}\in
\mathcal{H}$ such that%
\begin{equation*}
\left\{ \varphi \in \mathcal{H}:\left\vert \left\langle \psi _{k},\varphi
\right\rangle _{\mathcal{H}}\right\vert <1,\text{ }k=1,\ldots ,n\right\}
\subseteq \mathcal{U}\ .
\end{equation*}%
By (\ref{polar}), it follows that the absolute polar $\mathcal{U}^{\circ }$
is orthogonal to the set $\{\psi _{1},\ldots ,\psi _{n}\}^{\perp }$, leading
to%
\begin{equation*}
\mathcal{U}^{\circ }\subseteq \mathrm{span}\left\{ \psi _{1},\ldots ,\psi
_{n}\right\} \ .
\end{equation*}%
In this specific case, (weak) neighborhoods are too big, implying too small
absolute polars. If one takes instead the usual norm topology of the
separable infinite-dimensional Hilbert space to define neighborhoods, then
we obtain infinite-dimensional absolute polars, for all bounded
neighborhoods, like in any infinite-dimensional Banach space.

We give a general sufficient condition on a $0$-neighborhood $\mathcal{U}\in
\mathcal{O}$ for the infinite-dimensionality of its absolute polar $\mathcal{%
U}^{\circ }$:

\begin{condition}[Infinite dimensionality of absolute polars]
\label{condition}\mbox{ }\newline
There exists an infinite set $\{\sigma _{n}\}_{n\in \mathbb{N}}$ of linearly
independent elements $\sigma _{n}\in \mathcal{X}^{\ast }$ such that $\sup
\left\vert \sigma _{n}\left( \mathcal{U}\right) \right\vert <\infty $ for
every $n\in \mathbb{N}$.
\end{condition}

Condition \ref{condition} obviously implies the infinite-dimensionality of
the polar $\mathcal{U}^{\circ }$. In particular, $\mathcal{X}$ and $\mathcal{%
X}^{\ast }$ are infinite-dimensional like in \cite%
{Klee,FonfLindenstrauss,infinite dim convexity}. Such a condition can be
satisfied within a very \emph{large} class of topological vector
spaces:\medskip

\noindent \underline{Example 1:} If $\mathcal{X}$ is an infinite-dimensional
Banach space then any bounded neighborhood $\mathcal{U}\in \mathcal{O}$
satisfies Condition \ref{condition}, since it is contained in an open ball
of center $0$ in $\mathcal{X}$ (so that its polar $\mathcal{U}^{\circ }$
contains an open ball in $\mathcal{X}^{\ast }$ seen as a Banach
space).\medskip

\noindent \underline{Example 2:} If $\mathcal{U}\in \mathcal{O}$ is a
finite-dimensional\footnote{%
A subset of a vector space is finite-dimensional if its span is a
finite-dimensional subspace.} $0$-neighborhood in $\mathcal{X}$ then, by
\cite[Theorem 5.110]{Aliprantis},
\begin{equation*}
\mathrm{dim}\left\{ \sigma \in \mathcal{X}^{\ast }:\sigma \left( \mathcal{U}%
\right) =\left\{ 0\right\} \right\} =\mathrm{dim}\mathcal{X}^{\ast }-\mathrm{%
dim}\left( \mathrm{span}\mathcal{U}\right) \ .
\end{equation*}%
As a consequence, $\mathcal{U}$ satisfies Condition \ref{condition} whenever
$\mathcal{X}$ has infinite dimension. Existence of a finite-dimensional $%
\mathcal{U}\in \mathcal{O}$ is obviously ensured when $\mathcal{X}$ is a
locally finite-dimensional space, meaning that each point of $\mathcal{X}$
has a finite-dimensional neighborhood \cite[Definition 5]{locally finite dim}%
. A typical example of such topological vector spaces is given by the union
\begin{equation*}
\mathcal{X}=\bigcup\limits_{n\in \mathbb{N}}\mathcal{X}_{n}
\end{equation*}%
of an increasing sequence of finite-dimensional (normed) spaces $\mathcal{X}%
_{n}$ with diverging dimension \textrm{dim}$\mathcal{X}_{n}\rightarrow
\infty $, as $n\rightarrow \infty $, whose topology has as a $0$-basis the
family of all open balls of $\mathcal{X}_{n}$, $n\in \mathbb{N}$, centered
at $0$. \medskip

\noindent \underline{Example 3:} Let $\mathcal{X}$ be any vector space with
\emph{algebraic} dual space $\mathcal{X}^{\prime }$, that is, the vector
space of all linear functionals on $\mathcal{X}$. Assume the existence of an
infinite set $\{\sigma _{n}\}_{n\in \mathbb{N}}\subseteq \mathcal{X}^{\prime
}$ of linearly independent linear functionals $\sigma _{n}$, $n\in \mathbb{N}
$, which separates points in $\mathcal{X}$ and is pointwise uniformly
bounded:
\begin{equation*}
\sup_{n\in \mathbb{N}}\left\vert \sigma _{n}\left( A\right) \right\vert
<\infty \ ,\qquad A\in \mathcal{X}\ .
\end{equation*}%
Pick any topology on $\mathcal{X}$ such that $\{\sigma _{n}\}_{n\in \mathbb{N%
}}\subseteq \mathcal{X}^{\ast }$. If $\mathcal{U}\in \mathcal{O}$ is such
that
\begin{equation*}
\mathcal{U}\subseteq \left\{ A\in \mathcal{X}:\left\vert \sigma _{n}\left(
A\right) \right\vert <M\left( n\right) \right\}
\end{equation*}%
for some fixed, possibly unbounded, $M:\mathbb{N\rightarrow R}^{+}$, then $%
\mathcal{U}$ satisfies Condition \ref{condition}. As an example of such a
topology on $\mathcal{X}$, take the topology induced by the invariant metric%
\begin{equation*}
d_{\{\sigma _{n}\}_{n\in \mathbb{N}}}\left( A_{1},A_{2}\right) \doteq
\sum\limits_{n\in \mathbb{N}}\frac{1}{M\left( n\right) }\left\vert \sigma
_{n}\left( A_{1}-A_{2}\right) \right\vert <\infty \ ,\qquad A_{1},A_{2}\in
\mathcal{X}\ ,
\end{equation*}%
for some fixed $M:\mathbb{N\rightarrow R}^{+}$ with summable inverse. In
this case, any $\mathcal{U}\in \mathcal{O}$ contained in some open ball
satisfies Condition \ref{condition}.

\subsection{Weak$^{\ast }$-Hausdorff Dense Subsets of Convex Weak$^{\ast }$%
-Compact Sets\label{Dense Subsets}}

We first study generic convex weak$^{\ast }$-compact subsets of an
infinite-dimensional absolute polar. More precisely, taking a $0$%
-neighborhood $\mathcal{U}\in \mathcal{O}$ satisfying Condition \ref%
{condition}, we consider the set defined by (\ref{CU0}), that is,
\begin{equation*}
\mathbf{CU}_{\mathcal{U}}\left( \mathcal{X}^{\ast }\right) \doteq \left\{
U\in \mathbf{CF}\left( \mathcal{X}^{\ast }\right) :U\subseteq \mathcal{U}%
^{\circ }\right\} \subseteq \mathbf{CK}\left( \mathcal{X}^{\ast }\right)
\subseteq \mathbf{CF}\left( \mathcal{X}^{\ast }\right)
\end{equation*}%
with $\mathcal{U}^{\circ }\in \mathbf{CK}(\mathcal{X}^{\ast })$ being the
absolute polar (\ref{polar}) of $\mathcal{U}$.

Weak$^{\ast }$-closed subsets of absolute polars are the main source of weak$%
^{\ast }$-compact sets of $\mathcal{X}^{\ast }$ and, for any $\mathcal{U}\in
\mathcal{O}$, the hyperspace $\mathbf{CU}_{\mathcal{U}}(\mathcal{X}^{\ast })$
is a weak$^{\ast }$-Hausdorff-closed subspace of the Hausdorff hyperspace $%
\mathbf{CF}(\mathcal{X}^{\ast })$, by Corollaries \ref{convexity corrolary}
and \ref{convexity corrolary copy(2)} (iii). We study the topological
properties of the subset
\begin{equation}
\mathcal{D}_{\mathcal{U}}\doteq \left\{ U\in \mathbf{CU}_{\mathcal{U}}\left(
\mathcal{X}^{\ast }\right) :U=\overline{\mathcal{E}\left( U\right) }\right\}
\subseteq \mathbf{CU}_{\mathcal{U}}\left( \mathcal{X}^{\ast }\right)
\label{Zbis}
\end{equation}%
for any fixed $0$-neighborhood $\mathcal{U}\in \mathcal{O}$, where we recall
that $\mathcal{E}(U)$ is the extreme boundary of $U$, i.e., the (nonempty,
by the Krein-Milman theorem \cite[Theorem 3.23]{Rudin}) set of extreme
points of the convex weak$^{\ast }$-compact set $U$.

Note that the so-called \emph{exposed} points are particular examples of
extreme ones: a point $\sigma _{0}\in C$ in a convex subset $C\subseteq
\mathcal{X}^{\ast }$ is \emph{exposed} if there is $A\in \mathcal{X}$ such
that the real part of the weak$^{\ast }$-continuous functional $\hat{A}%
:\sigma \mapsto \sigma (A)$ from $\mathcal{X}^{\ast }$ to $\mathbb{K}=%
\mathbb{R},\mathbb{C}$ (cf. (\ref{sdfsdfkljsdlfkj})) takes its \emph{unique}
maximum on $C$ at $\sigma _{0}\in C$. Considering exposed points instead of
general extreme points is technically convenient. So, we denote by $\mathcal{%
E}_{0}(U)\subseteq \mathcal{E}(U)$ the set of all exposed points of any
element $U\in \mathbf{CU}_{\mathcal{U}}(\mathcal{X}^{\ast })$ and define%
\begin{equation}
\mathcal{D}_{0,\mathcal{U}}\doteq \left\{ U\in \mathbf{CU}_{\mathcal{U}%
}\left( \mathcal{X}^{\ast }\right) :U=\overline{\mathcal{E}_{0}\left(
U\right) }\right\} \subseteq \mathcal{D}_{\mathcal{U}}\subseteq \mathbf{CU}_{%
\mathcal{U}}\left( \mathcal{X}^{\ast }\right)  \label{Zbisexposed}
\end{equation}%
for any fixed $0$-neighborhood $\mathcal{U}\in \mathcal{O}$.

In order to have the weak$^{\ast }$-Hausdorff density of $\mathcal{D}_{0,%
\mathcal{U}}\subseteq \mathcal{D}_{\mathcal{U}}$ in $\mathbf{CU}_{\mathcal{U}%
}(\mathcal{X}^{\ast })$, the absolute polar $\mathcal{U}^{\circ }$ must be
infinite-dimensional, that is, its (linear) span $\mathcal{U}^{\circ }$ is
an infinite-dimensional subspace of $\mathcal{X}^{\ast }$. This refers to
Condition \ref{condition} on the $0$-neighborhood $\mathcal{U}\subseteq
\mathcal{X}$. Then, like in the proof of \cite[Theorem 4.3]%
{FonfLindenstrauss} and in contrast with \cite{Klee}, we design elements of $%
\mathcal{D}_{0,\mathcal{U}}$ that approximate $U\in \mathbf{CU}_{\mathcal{U}%
}(\mathcal{X}^{\ast })$ by using a procedure that is very similar to the
construction of the Poulsen simplex \cite{Poulsen}. Note however that
Poulsen used the existence of orthonormal bases in infinite-dimensional
Hilbert spaces\footnote{%
In \cite{Poulsen}, Poulsen uses the Hilbert space $\ell ^{2}(\mathbb{N})$ to
construct his example of a convex compact set (in fact a simplex) with dense
extreme boundary.}. Here, the Hahn-Banach separation theorem \cite[Theorem
3.4 (b)]{Rudin} replaces the orthogonality property coming from the Hilbert
space structure. In all previous results \cite{Klee,FonfLindenstrauss} on
the density of convex compact sets with dense extreme boundary, the norm
topology is used, while the primordial topology is here the weak$^{\ast }$
topology. In this context, the metrizability of weak$^{\ast }$ and weak$%
^{\ast }$-Hausdorff topologies on absolute polars is pivotal. See Theorem %
\ref{Solution selfbaby copy(4)+1}. We give now the precise assertion along
with its proof:

\begin{theorem}[Weak$^{\ast }$-Hausdorf density of $\mathcal{D}_{0,\mathcal{U%
}}$]
\label{Solution selfbaby copy(5)+0}\mbox{ }\newline
Let $\mathcal{X}$ be a separable topological $\mathbb{K}$-vector space with $%
\mathbb{K}=\mathbb{R},\mathbb{C}$ and $\mathcal{U}\in \mathcal{O}$ satisfy
Condition \ref{condition}. Then, $\mathcal{D}_{0,\mathcal{U}}\subseteq
\mathcal{D}_{\mathcal{U}}$ is a weak$^{\ast }$-Hausdorff dense subset of $%
\mathbf{CU}_{\mathcal{U}}(\mathcal{X}^{\ast })$.
\end{theorem}

\begin{proof}
Let $\mathcal{X}$ be a separable (infinite-dimensional) topological $\mathbb{%
K}$-vector space ($\mathbb{K}=\mathbb{R},\mathbb{C}$) and fix once and for
all $\mathcal{U}\in \mathcal{O}$ satisfying Condition \ref{condition}, as
well as a convex weak$^{\ast }$-compact subset $U\in \mathbf{CU}_{\mathcal{U}%
}(\mathcal{X}^{\ast })$. The construction of convex weak$^{\ast }$-compact
sets in $\mathcal{D}_{0,\mathcal{U}}$ approximating $U$ is done in several
steps:\medskip

\noindent \underline{Step 0:} Let $d$ be the metric defined by (\ref%
{metrics0}) and generating the weak$^{\ast }$ topology on the absolute polar
$\mathcal{U}^{\circ }$. Then, for any $\varepsilon \in \mathbb{R}^{+}$,
there is a finite set $\{\omega _{j}\}_{j=1}^{n_{\varepsilon }}\subseteq U$,
$n_{\varepsilon }\in \mathbb{N}$, such that
\begin{equation}
U\subseteq \bigcup\limits_{j=1}^{n_{\varepsilon }}B\left( \omega
_{j},\varepsilon \right) \ ,  \label{dense1}
\end{equation}%
where%
\begin{equation}
B\left( \omega ,r\right) \doteq \left\{ \sigma \in \mathcal{U}^{\circ
}:d\left( \omega ,\sigma \right) <r\right\} \subseteq \mathcal{U}^{\circ }
\label{ball weak}
\end{equation}%
denotes the weak$^{\ast }$-open ball of radius $r\in \mathbb{R}^{+}$,
centered at $\omega \in \mathcal{U}^{\circ }$. We then define the convex weak%
$^{\ast }$-compact set
\begin{equation}
U_{0}\doteq \mathrm{co}\left\{ \omega _{1},\ldots ,\omega _{n_{\varepsilon
}}\right\} \subseteq \mathrm{span}\{\omega _{1},\ldots ,\omega
_{n_{\varepsilon }}\}\ .  \label{K0}
\end{equation}%
By (\ref{dense1}), note that
\begin{equation}
\mathfrak{d}_{H}(U,U_{0})\leq \varepsilon \ ,  \label{K0bis}
\end{equation}%
where $\mathfrak{d}_{H}$ is the Hausdorff distance associated with the
metric $d$, as defined by (\ref{metric1}).\smallskip

\noindent \underline{Step 1:} Observe that the absolute polar $\mathcal{U}%
^{\circ }$ is weak$^{\ast }$-separable, by its weak$^{\ast }$ compactness
\cite[Theorem 3.15]{Rudin} and its metrizability \cite[Theorem 3.16]{Rudin}.
Take any weak$^{\ast }$-dense countable set $\{\varrho _{0,k}\}_{k\in
\mathbb{N}}$ of $U_{0}$. By Condition \ref{condition} together with a simple
rescaling argument, note also the existence of an infinite set of linearly
independent, continuous linear functionals within $\mathcal{U}^{\circ }$. As
a consequence, there is $\sigma _{1}\in \mathcal{U}^{\circ }\backslash
\mathrm{span}\{\omega _{1},\ldots ,\omega _{n_{\varepsilon }}\}$ satisfying
\begin{equation}
\sup \left\vert \sigma _{1}\left( \mathcal{U}\right) \right\vert =1\ .
\label{toto0}
\end{equation}%
As in the proof of Proposition \ref{convexity lemma copy(1)}, recall that $%
\mathcal{X}^{\ast }$, endowed with the weak$^{\ast }$ topology, is a locally
convex space with $\mathcal{X}$ being its dual space. Since $\{\sigma _{1}\}$
is a convex weak$^{\ast }$-compact set and $\mathrm{span}\{\omega
_{1},\ldots ,\omega _{n_{\varepsilon }}\}$ is convex and weak$^{\ast }$%
-closed \cite[Theorem 1.42]{Rudin}, we infer from the Hahn-Banach separation
theorem \cite[Theorem 3.4 (b)]{Rudin} the existence of $A_{1}\in \mathcal{X}$
such that
\begin{equation*}
\sup \left\{ \mathrm{Re}\left\{ \sigma \left( A_{1}\right) \right\} :\sigma
\in \mathrm{span}\{\omega _{1},\ldots ,\omega _{n_{\varepsilon }}\}\right\} <%
\mathrm{Re}\left\{ \sigma _{1}\left( A_{1}\right) \right\} \ .
\end{equation*}%
Since $\mathrm{span}\{\omega _{1},\ldots ,\omega _{n_{\varepsilon }}\}$ is a
linear space, observe that
\begin{equation}
\mathrm{Re}\left\{ \sigma \left( A_{1}\right) \right\} =0\ ,\qquad \sigma
\in \mathrm{span}\{\omega _{1},\ldots ,\omega _{n_{\varepsilon }}\}\ .
\label{cool1}
\end{equation}%
Additionally, by rescaling $A_{1}\in \mathcal{X}$, we can assume without
loss of generality that%
\begin{equation}
\mathrm{Re}\left\{ \sigma _{1}\left( A_{1}\right) \right\} =1\ .
\label{cool2}
\end{equation}%
Let
\begin{equation}
\omega _{n_{\varepsilon }+1}\doteq \left( 1-\lambda _{1}\right) \varpi
_{1}+\lambda _{1}\sigma _{1}\ ,\qquad \text{with}\qquad \lambda _{1}\doteq
\min \left\{ 1,2^{-2}\varepsilon \right\} ,\ \varpi _{1}\doteq \varrho
_{0,1}\in U_{0}\ .  \label{omega1}
\end{equation}%
In contrast with the proof of \cite[Theorem 4.3]{FonfLindenstrauss}, we use
a convex combination to automatically ensure that $\omega _{n_{\varepsilon
}+1}\in \mathcal{U}^{\circ }$, by convexity of the absolute polar $\mathcal{U%
}^{\circ }$. By (\ref{metrics0}), the inequality $\lambda _{1}\leq
2^{-2}\varepsilon $ yields%
\begin{equation}
d\left( \omega _{n_{\varepsilon }+1},\varpi _{1}\right) \leq
2^{-1}\varepsilon \ .  \label{dense2dense2}
\end{equation}%
Define the new convex weak$^{\ast }$-compact set%
\begin{equation*}
U_{1}\doteq \mathrm{co}\left\{ \omega _{1},\ldots ,\omega _{n_{\varepsilon
}+1}\right\} \subseteq \mathrm{span}\{\omega _{1},\ldots ,\omega
_{n_{\varepsilon }+1}\}\ .
\end{equation*}%
Observe that $\omega _{n_{\varepsilon }+1}$ is an exposed point of $U_{1}$,
by (\ref{cool1}) and (\ref{cool2}). By (\ref{metric1}), (\ref{K0}) and (\ref%
{dense2dense2}), note also that $\mathfrak{d}_{H}(U_{0},U_{1})\leq
2^{-1}\varepsilon $, which, by the triangle inequality and (\ref{K0bis}),
yields
\begin{equation}
\mathfrak{d}_{H}(U,U_{1})\leq \left( 1+2^{-1}\right) \varepsilon
\label{K1bis}
\end{equation}%
for an arbitrary (but previously fixed) $\varepsilon \in \mathbb{R}^{+}$%
.\medskip

\noindent \underline{Step 2:} Take any weak$^{\ast }$ dense countable set $%
\{\varrho _{1,k}\}_{k\in \mathbb{N}}$ of $U_{1}$. By Condition \ref%
{condition}, there is $\sigma _{2}\in \mathcal{U}^{\circ }\backslash \mathrm{%
span}\{\omega _{1},\ldots ,\omega _{n_{\varepsilon }+1}\}$ with
\begin{equation}
\sup \left\vert \sigma _{2}\left( \mathcal{U}\cup \{A_{1}\}\right)
\right\vert \leq \min \left\{ 1,2^{-1}\lambda _{1}\right\} \ .
\label{sdfsdf}
\end{equation}%
As before, we deduce from the Hahn-Banach separation theorem \cite[Theorem
3.4 (b)]{Rudin} the existence of $A_{2}\in \mathcal{X}$ such that
\begin{equation}
\mathrm{Re}\left\{ \sigma _{2}\left( A_{2}\right) \right\} =1\qquad \text{and%
}\qquad \mathrm{Re}\left\{ \sigma \left( A_{2}\right) \right\} =0\ ,\qquad
\sigma \in \mathrm{span}\{\omega _{1},\ldots ,\omega _{n_{\varepsilon
}+1}\}\ .  \label{cool3}
\end{equation}%
Let%
\begin{equation}
\omega _{n_{\varepsilon }+2}\doteq \left( 1-\lambda _{2}\right) \varpi
_{2}+\lambda _{2}\sigma _{2}\ ,\qquad \text{with}\qquad \lambda _{2}\doteq
\min \left\{ 1,2^{-3}\varepsilon \right\} ,\ \varpi _{2}\doteq \varrho
_{1,1}\in U_{1}\ .  \label{omega2}
\end{equation}%
In this case, similar to Inequality (\ref{dense2dense2}),
\begin{equation}
d\left( \omega _{n_{\varepsilon }+2},\varpi _{2}\right) \leq
2^{-2}\varepsilon \ .  \label{dense2dense3}
\end{equation}%
Define the new convex weak$^{\ast }$-compact set%
\begin{equation*}
U_{2}\doteq \mathrm{co}\left\{ \omega _{1},\ldots ,\omega _{n_{\varepsilon
}+2}\right\} \subseteq \mathrm{span}\{\omega _{1},\ldots ,\omega
_{n_{\varepsilon }+2}\}\ .
\end{equation*}%
By (\ref{cool3}), $\omega _{n_{\varepsilon }+2}$ is an exposed point of $%
U_{2}$, but it is not obvious that the exposed point $\omega
_{n_{\varepsilon }+1}$ of $U_{1}$ is still an exposed point of $U_{2}$, with
respect to $A_{1}\in \mathcal{X}$. This property is a consequence of%
\begin{equation*}
\mathrm{Re}\left\{ \omega _{n_{\varepsilon }+2}\left( A_{1}\right) \right\}
=\left( 1-\lambda _{2}\right) \mathrm{Re}\left\{ \varpi _{2}\left(
A_{1}\right) \right\} +\lambda _{2}\mathrm{Re}\left\{ \sigma _{2}\left(
A_{1}\right) \right\} <\mathrm{Re}\left\{ \omega _{n_{\varepsilon }+1}\left(
A_{1}\right) \right\} =\lambda _{1}\ ,
\end{equation*}%
(see (\ref{cool1}), (\ref{omega1}) and (\ref{omega2})), which holds true
because of Equation (\ref{sdfsdf}). By (\ref{metric1}), (\ref{K1bis}) and (%
\ref{dense2dense3}) together with the triangle inequality,%
\begin{equation*}
\mathfrak{d}_{H}(U,U_{2})\leq \left( 1+2^{-1}+2^{-2}\right) \varepsilon
\end{equation*}%
for an arbitrary (but previously fixed) $\varepsilon \in \mathbb{R}^{+}$%
.\medskip

\noindent \underline{Step $n\rightarrow \infty $:} We now iterate the above
procedure, ensuring, at each step $n\geq 3$, that the addition of the
element
\begin{equation}
\omega _{n_{\varepsilon }+n}\doteq \left( 1-\lambda _{n}\right) \varpi
_{n}+\lambda _{n}\sigma _{n}\ ,\qquad \text{with}\qquad \lambda _{n}\doteq
\min \left\{ 1,2^{-(n+1)}\varepsilon \right\} \ ,  \label{definition omegan}
\end{equation}%
in order to define the convex weak$^{\ast }$-compact set
\begin{equation}
U_{n}\doteq \mathrm{co}\left\{ \omega _{1},\ldots ,\omega _{n_{\varepsilon
}+n}\right\} \subseteq \mathrm{span}\{\omega _{1},\ldots ,\omega
_{n_{\varepsilon }+n}\}\ ,  \label{Kn}
\end{equation}%
does not destroy the property\ of the elements $\omega _{n_{\varepsilon
}+1},\ldots ,\omega _{n_{\varepsilon }+n-1}$ being exposed. To this end, for
any $n\geq 2$, we choose $\sigma _{n}\in \mathcal{U}^{\circ }\backslash
\mathrm{span}\{\omega _{1},\ldots ,\omega _{n_{\varepsilon }+n-1}\}$ such
that%
\begin{equation}
\sup \left\vert \sigma _{n}\left( \mathcal{U}\cup \{A_{1}\}\cup \cdots \cup
\{A_{n-1}\}\right) \right\vert \leq \min \left\{ 1,2^{-1}\lambda _{1},\ldots
,2^{-1}\lambda _{n-1}\right\} \ .  \label{toto}
\end{equation}%
Here, for any integer $n\geq 2$ and $j\in \{1,\ldots ,n-1\}$, $A_{j}\in
\mathcal{X}$ satisfies
\begin{equation}
\mathrm{Re}\left\{ \sigma _{j}\left( A_{j}\right) \right\} =1\qquad \text{and%
}\qquad \mathrm{Re}\left\{ \sigma \left( A_{j}\right) \right\} =0\ ,\qquad
\sigma \in \mathrm{span}\{\omega _{1},\ldots ,\omega _{n_{\varepsilon
}+j-1}\}\ .  \label{totototo}
\end{equation}%
We also have to conveniently choose $\varpi _{n}\in U_{n-1}$ in order to get
the asserted weak$^{\ast }$ density. Like in the proof of \cite[Theorem 4.3]%
{FonfLindenstrauss} the sequence $(\varpi _{n})_{n\in \mathbb{N}}$ is chosen
such that
\begin{equation*}
\{\varpi _{n}\}_{n\in \mathbb{N}}=\left\{ \varrho _{n,k}\right\} _{n\in
\mathbb{N}_{0},k\in \mathbb{N}}
\end{equation*}%
and all the functionals $\varrho _{n,k}$ appear infinitely many times in the
sequence $(\varpi _{n})_{n\in \mathbb{N}}$. In this case, we obtain a weak$%
^{\ast }$-dense set $\{\omega _{n}\}_{n\in \mathbb{N}}$ in the convex weak$%
^{\ast }$-compact set
\begin{equation}
U_{\infty }\doteq \overline{\mathrm{co}\left\{ \{\omega _{n}\}_{n\in \mathbb{%
N}}\right\} }\in \mathbf{CU}_{\mathcal{U}}\left( \mathcal{X}^{\ast }\right)
\ ,  \label{equaion}
\end{equation}%
which, by construction, satisfies
\begin{equation*}
\mathfrak{d}_{H}(U,U_{\infty })\leq \sum_{n=0}^{\infty }2^{-n}\varepsilon
=2\varepsilon
\end{equation*}%
for an arbitrary (but previously fixed) $\varepsilon \in \mathbb{R}^{+}$.
\medskip

\noindent \underline{Step $n=\infty $:} It remains to verify that $\omega
_{n_{\varepsilon }+j}$, $j\in \mathbb{N}$, are exposed points of $U_{\infty
} $, whence $U_{\infty }\in \mathcal{D}$. By (\ref{definition omegan}) with $%
\varpi _{n}\in U_{n-1}$ (see (\ref{Kn})), for each natural number $n\geq j+1$%
, there are $\alpha _{n,j-1}^{(j)},\ldots ,\alpha _{n,n}^{(j)}\in \lbrack
0,1]$ and $\rho _{n}^{(j)}\in \mathrm{co}\left\{ \omega _{1},\ldots ,\omega
_{n_{\varepsilon }+j-1}\right\} $ such that
\begin{equation}
\alpha _{n,j-1}^{(j)}+\alpha _{n,j}^{(j)}+\sum_{k=j+1}^{n}\alpha
_{n,k}^{(j)}\lambda _{k}=1\quad \text{and}\quad \omega _{n_{\varepsilon
}+n}=\alpha _{n,j-1}^{(j)}\rho _{n}^{(j)}+\alpha _{n,j}^{(j)}\omega
_{n_{\varepsilon }+j}+\sum_{k=j+1}^{n}\alpha _{n,k}^{(j)}\lambda _{k}\sigma
_{k}\ .  \label{inequality ddfdfinequality ddfdf}
\end{equation}%
Additionally, define $\alpha _{n,k}^{(j)}\doteq 1$ for all natural numbers $%
k\geq n$ while $\alpha _{n,k}^{(j)}\doteq 0$ for $k\in \mathbb{N}_{0}$ such
that $k\leq j-2$. Using (\ref{toto}), (\ref{totototo}) and (\ref{inequality
ddfdfinequality ddfdf}), at fixed $j\in \mathbb{N}$, we thus obtain that%
\begin{eqnarray}
\mathrm{Re}\left\{ \omega _{n_{\varepsilon }+n}\left( A_{j}\right) \right\}
&=&\alpha _{n,j}^{(j)}\mathrm{Re}\left\{ \omega _{n_{\varepsilon }+j}\left(
A_{j}\right) \right\} +\sum_{k=j+1}^{n}\alpha _{n,k}^{(j)}\lambda _{k}%
\mathrm{Re}\left\{ \sigma _{k}\left( A_{j}\right) \right\}  \notag \\
&\leq &\alpha _{n,j}^{(j)}\lambda _{j}+\sum_{k=j+1}^{n}\alpha
_{n,k}^{(j)}\lambda _{k}\left( 2^{-1}\lambda _{j}\right) \leq \lambda
_{j}\left( 1-2^{-1}\sum_{k=j+1}^{n}\alpha _{n,k}^{(j)}\lambda _{k}\right)
\label{inequality ddfdf}
\end{eqnarray}%
for any $n\geq j+1$, while, for any natural number $n\leq j-1$,
\begin{equation*}
\mathrm{Re}\left\{ \omega _{n_{\varepsilon }+n}\left( A_{j}\right) \right\}
=0\ ,
\end{equation*}%
using (\ref{totototo}). Fix $j\in \mathbb{N}$ and let $\omega _{\infty }\in
U_{\infty }$ be a solution to the variational problem%
\begin{equation}
\max_{\sigma \in U_{\infty }}\mathrm{Re}\left\{ \sigma \left( A_{j}\right)
\right\} =\mathrm{Re}\left\{ \omega _{\infty }\left( A_{j}\right) \right\}
\geq \mathrm{Re}\left\{ \omega _{n_{\varepsilon }+j}\left( A_{j}\right)
\right\} =\lambda _{j}\ .  \label{maximum}
\end{equation}%
($U_{\infty }$ is weak$^{\ast }$-compact, by weak$^{\ast }$-compactness of $%
\mathcal{U}^{\circ }\supseteq U_{\infty }$.) By weak$^{\ast }$-density of $%
\{\omega _{n}\}_{n\in \mathbb{N}}$ in $U_{\infty }$, there is a sequence $%
(\omega _{n_{\varepsilon }+n_{l}})_{l\in \mathbb{N}}$ converging to $\omega
_{\infty }$ in the weak$^{\ast }$ topology. Since $U_{j}$ is weak$^{\ast }$%
-compact and $\alpha _{n,k}^{(j)}\in \lbrack 0,1]$ for all $k\in \mathbb{N}%
_{0}$ and $n,j\in \mathbb{N}$, by a standard argument with a so-called
diagonal subsequence, we can choose the sequence $(n_{l})_{l\in \mathbb{N}}$
such that $(\rho _{n_{l}}^{(j)})$ weak$^{\ast }$-converges to $\rho _{\infty
}^{(j)}\in U_{j-1}$, and $(\alpha _{n_{l},k}^{(j)})_{l\in \mathbb{N}}$ has a
limit $\alpha _{\infty ,k}^{(j)}\in \lbrack 0,1]$ for any fixed $k\in
\mathbb{N}_{0}$ and $j\in \mathbb{N}$. Using (\ref{definition omegan}), (\ref%
{inequality ddfdf}) and the inequality
\begin{equation*}
\sum_{k=j+1}^{n_{l}}\alpha _{n_{l},k}^{(j)}\lambda _{k}\leq \varepsilon
\sum_{k=j+1}^{\infty }2^{-(k+1)}=2^{-(j+1)}\varepsilon
\end{equation*}%
together with Lebesgue's dominated convergence theorem, we thus arrive at
\begin{equation*}
\mathrm{Re}\left\{ \omega _{\infty }\left( A_{j}\right) \right\}
=\lim_{l\rightarrow \infty }\mathrm{Re}\left\{ \omega _{n_{\varepsilon
}+n_{l}}\left( A_{j}\right) \right\} \leq \lambda _{j}\left(
1-2^{-1}\sum_{k=j+1}^{\infty }\lambda _{k}\lim_{l\rightarrow \infty }\alpha
_{n_{l},k}^{(j)}\right) \ .
\end{equation*}%
Because of (\ref{maximum}), it follows that%
\begin{equation}
\alpha _{\infty ,k}^{(j)}\doteq \lim_{l\rightarrow \infty }\alpha
_{n_{l},k}^{(j)}=0\ ,\qquad k\in \{j+1,\ldots ,\infty \}\ .  \label{limit}
\end{equation}%
As absolute polars are weak$^{\ast }$ compact \cite[Theorem 3.15]{Rudin},
for any $A\in \mathcal{X}$, the continuous function $\hat{A}$ defined by (%
\ref{sdfsdfkljsdlfkj}) satisfies
\begin{equation*}
\sup |\hat{A}(\mathcal{U}^{\circ })|<\infty \ .
\end{equation*}%
Therefore, by (\ref{definition omegan}), (\ref{inequality ddfdfinequality
ddfdf})\ and Lebesgue's dominated convergence theorem, for all $A\in
\mathcal{X}$,
\begin{equation*}
\lim_{l\rightarrow \infty }\omega _{n_{\varepsilon }+n_{l}}\left( A\right)
=\lim_{l\rightarrow \infty }\alpha _{n_{l},j-1}^{(j)}\rho
_{n_{l}}^{(j)}\left( A\right) +\omega _{n_{\varepsilon }+j}\left( A\right)
\lim_{l\rightarrow \infty }\alpha _{n_{l},j}^{(j)}+\sum_{k=j+1}^{\infty
}\lambda _{k}\sigma _{k}\left( A\right) \lim_{l\rightarrow \infty }\alpha
_{n_{l},k}^{(j)}\ ,
\end{equation*}%
which combined with (\ref{limit}) implies that
\begin{equation*}
\lim_{l\rightarrow \infty }\omega _{n_{\varepsilon }+n_{l}}\left( A\right)
=\alpha _{\infty ,j-1}^{(j)}\rho _{\infty }^{(j)}\left( A\right) +\alpha
_{\infty ,j}^{(j)}\omega _{n_{\varepsilon }+j}\left( A\right) \ ,
\end{equation*}%
where $\rho _{\infty }^{(j)}\in U_{j-1}$, $\alpha _{\infty
,j-1}^{(j)},\alpha _{\infty ,j}^{(j)}\in \lbrack 0,1]$ and $\alpha _{\infty
,j-1}^{(j)}+\alpha _{\infty ,j}^{(j)}=1$. Hence, the sequence $(\omega
_{n_{\varepsilon }+n_{l}})_{l\in \mathbb{N}}$ weak$^{\ast }$ converges to an
element of $U_{j}$. (Recall that $U_{j}$ is defined by (\ref{Kn}) for $%
n=j\in \mathbb{N}$.) Since $\omega _{n_{\varepsilon }+j}$ is by construction
the unique maximizer of
\begin{equation*}
\max_{\sigma \in U_{j}}\mathrm{Re}\left\{ \sigma \left( A_{j}\right)
\right\} =\mathrm{Re}\left\{ \omega _{n_{\varepsilon }+j}\left( A_{j}\right)
\right\}
\end{equation*}%
and (\ref{maximum}) holds true with $\omega _{\infty }\in U_{j}$, we deduce
that $\omega _{\infty }=\omega _{n_{\varepsilon }+j}$, which is thus an
exposed point of $U_{\infty }$ for any $j\in \mathbb{N}$.
\end{proof}

Recall now that $\mathbf{CU}(\mathcal{X}^{\ast })$ is the set of all
nonempty, uniformly bounded in a $0$-neighborhood, weak$^{\ast }$-closed
subsets of $\mathcal{X}^{\ast }$ defined by (\ref{CU}), that is,
\begin{equation*}
\mathbf{CU}\left( \mathcal{X}^{\ast }\right) \doteq \bigcup\limits_{\mathcal{%
U}\in \mathcal{O}}\mathbf{CF}\left( \mathcal{X}^{\ast }\right) \cap \mathbf{U%
}_{\mathcal{U}}\left( \mathcal{X}^{\ast }\right) \doteq \bigcup\limits_{%
\mathcal{U}\in \mathcal{O}}\mathbf{CU}_{\mathcal{U}}\left( \mathcal{X}^{\ast
}\right) \ .
\end{equation*}%
Provided there is one $0$-neighborhood in $\mathcal{X}$ satisfying Condition %
\ref{condition}, Theorem \ref{Solution selfbaby copy(5)+0} directly implies
the weak$^{\ast }$-Hausdorff density in $\mathbf{CU}(\mathcal{X}^{\ast })$
of the sets $\mathcal{D}_{0}$ and $\mathcal{D}$, where%
\begin{equation}
\mathcal{D}_{0}\doteq \bigcup\limits_{\mathcal{U}\in \mathcal{O}}\mathcal{D}%
_{0,\mathcal{U}}\subseteq \mathcal{D}\doteq \bigcup\limits_{\mathcal{U}\in
\mathcal{O}}\mathcal{D}_{\mathcal{U}}  \label{Zbis-full}
\end{equation}%
are the sets of all $U\in \mathbf{CU}\left( \mathcal{X}^{\ast }\right) $
with weak$^{\ast }$-dense exposed, respectively extreme, boundary\ (see (\ref%
{Zbis}) and (\ref{Zbisexposed})):

\begin{corollary}[Weak$^{\ast }$-Hausdorf density of $\mathcal{D}_{0}$]
\label{Solution selfbaby copy(8)}\mbox{ }\newline
Let $\mathcal{X}$ be a separable topological $\mathbb{K}$-vector space with $%
\mathbb{K}=\mathbb{R},\mathbb{C}$ and assume the existence of one $0$%
-neighborhood in $\mathcal{X}$ satisfying Condition \ref{condition}. Then, $%
\mathcal{D}_{0}\subseteq \mathcal{D}$ is a weak$^{\ast }$-Hausdorff dense
subset of $\mathbf{CU}(\mathcal{X}^{\ast })$.
\end{corollary}

\begin{proof}
By (\ref{polar}), observe that, for any $\mathcal{U}_{1},\mathcal{U}_{2}\in
\mathcal{O}$,%
\begin{equation*}
\mathcal{U}_{1}\cap \mathcal{U}_{2}\in \mathcal{O\qquad }\text{and}\mathcal{%
\qquad U}_{1}^{\circ }\cup \mathcal{U}_{2}^{\circ }\subseteq \left( \mathcal{%
U}_{1}\cap \mathcal{U}_{2}\right) ^{\circ }\ .
\end{equation*}%
Using this together with Theorem \ref{Solution selfbaby copy(5)+0}, one
deduces the assertion.
\end{proof}

Our proof of Theorem \ref{Solution selfbaby copy(5)+0} differs in several
important aspects from the one of \cite[Theorem 4.3]{FonfLindenstrauss},
even if it has the same general structure, inspired by Poulsen's
construction \cite{Poulsen}, as already mentioned. To be more precise, as
compared to the proof of \cite[Theorem 4.3]{FonfLindenstrauss}, \emph{Step 0}
is new and is a direct consequence of the weak$^{\ast }$-compactness of $U$,
a property not assumed in \cite[Theorem 4.3]{FonfLindenstrauss}. \emph{Step 1%
} to \emph{Step }$n\rightarrow \infty $ are similar to what is done in \cite%
{FonfLindenstrauss}, but with the essential difference that convex
combinations are used to produce new (strongly)\ exposed points and the
required bounds on $\{\lambda _{n},\sigma _{n}\}_{n\in \mathbb{N}}$ are very
different. Compare Equations (\ref{definition omegan}) and (\ref{toto}) with
the bounds on $\upsilon _{1},\upsilon _{2},\upsilon _{3}$ given in \cite[p.
27-29]{FonfLindenstrauss}, at parameters $r_{1}(t),r_{2}(t),r_{3}(t)=1$. In
particular, there is no norm on $\mathcal{X}^{\ast }$ ($\mathcal{X}$ is not
necessarily a Banach space and, in any case, $\mathcal{X}^{\ast }$ is
endowed with the weak$^{\ast }$ topology) and we use estimates on convex
combinations that completely differ from what is done in \cite[Theorem 4.3]%
{FonfLindenstrauss}. This corresponds to \emph{Step }$n=\infty $.

Note that \cite[Theorem 4.3]{FonfLindenstrauss} shows the density of convex
norm-closed sets with dense set of \emph{strongly} exposed points. A
strongly exposed point $\sigma _{0}$ in some convex set $C\subseteq \mathcal{%
X}^{\ast }$ is an\emph{\ }exposed point for some $A\in \mathcal{X}$ with the
additional property that any minimizing net of the real part of $\hat{A}$
(cf. (\ref{sdfsdfkljsdlfkj})) has to converge to $\sigma _{0}$ in the weak$%
^{\ast }$ topology\footnote{%
One should not mistake the notion of strongly exposed points discussed here
for the notion of weak$^{\ast }$ strongly exposed points of \cite[Definition
5.8]{Phelps-Asplund}, where $\mathcal{X}$ is always a Banach space and a weak%
$^{\ast }$ strongly exposed point is a (weak$^{\ast }$) exposed point with
the additional property that any minimizing net of the real part of $\hat{A}$
has to converge to $\sigma _{0}$ \emph{in the norm topology} of $\mathcal{X}%
^{\ast }$.}. Note that such a minimizing net is weak$^{\ast }$-convergent
with the exposed point $\sigma _{0}$ being its limit, by weak$^{\ast }$
continuity of $\hat{A}$. If $C=K$ is weak$^{\ast }$-compact, this yields
that any minimizing net converges to $\sigma _{0}$ in the weak$^{\ast }$
topology. In other words, any exposed point is \emph{automatically} strongly
exposed in all convex weak$^{\ast }$-compact sets $K\in \mathbf{CK}(\mathcal{%
X}^{\ast })\supseteq \mathbf{CU}(\mathcal{X}^{\ast })$.

\subsection{Extension of the Straszewicz Theorem}

In this section, we study the relations between the set $\mathcal{D}_{%
\mathcal{U}}$ of convex weak$^{\ast }$-compact sets with weak$^{\ast }$%
-dense set of \emph{extreme} points and the set $\mathcal{D}_{0,\mathcal{U}}$
of convex weak$^{\ast }$-compact sets with weak$^{\ast }$-dense set of \emph{%
exposed} points, for any fixed $0$-neighborhood $\mathcal{U}\in \mathcal{O}$%
. See (\ref{Zbis}) and (\ref{Zbisexposed}). In fact, we give a very general
condition on the topological vector space $\mathcal{X}$ leading to the
equality $\mathcal{D}_{\mathcal{U}}=\mathcal{D}_{0,\mathcal{U}}$ for all $%
\mathcal{U}\in \mathcal{O}$. This result is used in Section \ref{Gdelta}.

Such a study is reminiscent of the Straszewicz theorem: In 1935, Straszewicz
proves \cite{Straszewicz} that the set of exposed points of a convex compact
space of a finite-dimensional space ($\mathbb{R}^{n}$) is \emph{dense} in
the set of extreme points. See, e.g., \cite[Theorem 7.89]{Aliprantis}. An
extension of this result\footnote{%
There is also a result \cite{Asplund} of Asplund in 1963 generalizing the
Straszewicz theorem to so-called $k$-exposed and $k$-extreme points in the
finite-dimensional space $\mathbb{R}^{n}$.}\ to convex (locally)
norm-compact (closed) subsets of an infinite-dimensional normed space was
performed by Klee in 1958, see \cite[Theorems (2.1), (2.3)]{Klee-exposed}.
In 1976, Bair \cite[Theorem 1]{Bair} proves the Straszewicz theorem in an
arbitrary real vector space for algebraically closed convex sets with
so-called finite \textquotedblleft copointure\textquotedblright , see \cite[%
Section II.5.1]{BairI}. This last condition cannot be satisfied by weak$%
^{\ast }$-compact sets. In fact, such studies on dual spaces $\mathcal{X}%
^{\ast }$ have been performed by Larman and Phelps in \cite{LarmanPhelps}
for special Banach spaces $\mathcal{X}$, named \emph{%
Gateaux-differentiability space }\cite[Definition 6.1]{Phelps-Asplund}. For
topological vector spaces, it means the following:

\begin{definition}[Gateaux-differentiability space]
\label{Gateaux Differentiability space}\mbox{ }\newline
A Gateaux-differentiability space $\mathcal{X}$ is a topological vector
space on which every continuous convex real-valued function with a nonempty
open convex subset as domain is Gateaux-differentiable on a dense set in
that domain.
\end{definition}

\noindent Compare with \cite[Definition 6.1]{Phelps-Asplund}. Recall that a
\emph{weak Asplund} space is a topological vector space on which every
continuous convex real-valued function with a nonempty open convex subset as
domain is Gateaux-differentiable on a \emph{generic} set in that domain.
See, e.g., \cite[Definition 1.22]{Phelps-Asplund} for the Banach case or
\cite[p. 203]{Mazur-extended} in the general case. Phelps explains in \cite[%
p. 95]{Phelps-Asplund} that the new space class of Definition \ref{Gateaux
Differentiability space} is obviously \textquotedblleft formally larger than
the class of weak Asplund spaces, but in some ways is a more natural object
of study.\textquotedblright\ All these spaces are reminiscent of the
celebrated Mazur theorem \cite[Satz 2]{Mazur} (see also \cite[Theorem 10.44]%
{BruPedra2}) proven for separable real Banach spaces.

If the Banach space $\mathcal{X}$ is a Gateaux-differentiability space then,
by \cite[Theorem 6.2]{Phelps-Asplund}, every convex weak$^{\ast }$-compact
set $K\in \mathbf{CK}(\mathcal{X}^{\ast })$ is the weak$^{\ast }$-closed
convex hull of its exposed points. By the Milman theorem \cite[Theorem 10.13]%
{BruPedra2} it follows that the set $\mathcal{E}_{0}(K)$ of exposed points
is weak$^{\ast }$-dense in the set $\mathcal{E}(K)$. \cite[Theorem 6.2]%
{Phelps-Asplund} refers to Banach spaces and we give here another extension
of the Straszewicz theorem to all (possibly non-Banach)
Gateaux-differentiability spaces, because the weak$^{\ast }$-density of the
set of exposed points is an important ingredient in the next subsection. Our
proof is quite direct and thus, relatively pedagogical while being very
general:

\begin{theorem}[Extension of the Straszewicz theorem - I]
\label{stra theorem}\mbox{ }\newline
Let $\mathcal{X}$ be a Gateaux-differentiability $\mathbb{K}$-vector space
with $\mathbb{K}=\mathbb{R},\mathbb{C}$. Then, for any convex weak$^{\ast }$%
-compact set $K\in \mathbf{CK}(\mathcal{X}^{\ast })$, the set of exposed
points of $K$ is weak$^{\ast }$-dense in $\mathcal{E}(K)$.
\end{theorem}

\begin{proof}
Fix all parameters of the theorem. Assume without loss of generality that $%
\mathcal{X}$ is a $\mathbb{C}$-vector space. (The case $\mathbb{K}=\mathbb{R}
$ is even slightly simpler.) Denote by $\mathcal{E}_{0}(K)\subseteq \mathcal{%
E}(K)$ the set of all exposed points of $K$. Assume that%
\begin{equation}
\overline{\mathrm{co}}\left( \mathcal{E}_{0}(K)\right) \subsetneq K\ ,
\label{contraduction}
\end{equation}%
i.e., there is an element $\sigma _{0}\in K\backslash \overline{\mathrm{co}}%
\left( \mathcal{E}_{0}(K)\right) $. We thus infer from the Hahn-Banach
separation theorem \cite[Theorem 3.4 (b)]{Rudin} the existence of $A_{0}\in
\mathcal{X}$ such that
\begin{equation}
\max \left\{ \mathrm{Re}\left\{ \sigma \left( A_{0}\right) \right\} :\sigma
\in \overline{\mathrm{co}}\left( \mathcal{E}_{0}(K)\right) \right\} <\mathrm{%
Re}\left\{ \sigma _{0}\left( A_{0}\right) \right\} \ .
\label{sdkljsdlkfjsdlkj}
\end{equation}%
Since $K$ is weak$^{\ast }$-compact,
\begin{equation}
h\left( A\right) \doteq \max_{\sigma \in K}\mathrm{Re}\{\hat{A}\left( \sigma
\right) \}\ ,\qquad A\in \mathcal{X}\ ,  \label{mazur magique}
\end{equation}%
defines a continuous mapping $h:\mathcal{X}\rightarrow \mathbb{R}$. Observe
that $h$ is a convex function because it is the maximum of a family of
linear functions. Recall that a tangent $\mathbb{R}$-linear functional $%
\mathrm{d}h(A)$ at $A\in \mathcal{Y}$ satisfies%
\begin{equation*}
h\left( A+B\right) -h\left( A\right) \geq \left[ \mathrm{d}h\left( A\right) %
\right] \left( B\right) \ ,\qquad B\in \mathcal{X}\ ,
\end{equation*}%
by definition. Any maximizer $\sigma _{A}\in K$ of the variational problem (%
\ref{mazur magique}) yields a (continuous) tangent ($\mathbb{R}$-linear)
functional:
\begin{equation*}
h(A+B)-h(A)\geq \mathrm{Re}\{\sigma _{A}\left( B\right) \}\ ,\qquad B\in
\mathcal{X}\ .
\end{equation*}%
Therefore, if there is a unique continuous tangent ($\mathbb{R}$-linear)
functional $\mathrm{d}h(A)$ at $A\in \mathcal{Y}$, then the solution $\sigma
_{A}\in K$ of the variational problem (\ref{mazur magique}) is unique:
\begin{equation*}
\sigma _{A}\left( B\right) =\left[ \mathrm{d}h\left( A\right) \right] \left(
B\right) +i\left[ \mathrm{d}h\left( A\right) \right] \left( -iB\right) \
,\qquad B\in \mathcal{X}\ ,
\end{equation*}%
by $\mathbb{C}$-linearity of $\sigma _{A}$. In particular, $\sigma _{A}$ is,
in this case, an exposed point of $K$. By Definition \ref{Gateaux
Differentiability space}, there is a net $(\sigma _{j})_{j\in J}$ of exposed
points of $K$ as well as a net $(A_{j})_{j\in J}\in \mathcal{X}$ converging
to $A_{0}$ such that%
\begin{equation}
\mathrm{Re}\{\sigma _{j}\left( A_{j}\right) \}=\max_{\sigma \in K}\mathrm{Re}%
\{\hat{A}_{j}\left( \sigma \right) \}\ ,\qquad j\in J\ .  \label{toto1}
\end{equation}%
By taking any maximizer of $\mathrm{Re}\{\hat{A}_{0}\left( \sigma \right) \}$
over $\sigma \in K$, note that%
\begin{equation}
\underset{J}{\lim \inf }\mathrm{Re}\{\sigma _{j}\left( A_{j}\right) \}=%
\underset{J}{\lim \inf }\max_{\sigma \in K}\mathrm{Re}\{\hat{A}_{j}\left(
\sigma \right) \}\geq \max_{\sigma \in K}\mathrm{Re}\{\hat{A}_{0}\left(
\sigma \right) \}  \label{toto2}
\end{equation}%
while, by compactness of $K$, we can assume without loss of generality that $%
(\sigma _{j})_{j\in J}$ converges to some $\sigma _{\infty }\in K$. Now, by
the Banach-Steinhaus theorem \cite[Theorem 2.5]{Rudin} together with (\ref%
{toto1})-(\ref{toto2}), it follows that
\begin{equation*}
\mathrm{Re}\{\sigma _{\infty }\left( A_{0}\right) \}=\lim_{J}\mathrm{Re}%
\{\sigma _{j}\left( A_{0}\right) \}=\lim_{J}\mathrm{Re}\{\sigma _{j}\left(
A_{j}\right) \}=\lim_{J}\max_{\sigma \in K}\mathrm{Re}\{\hat{A}_{j}\left(
\sigma \right) \}=\max_{\sigma \in K}\mathrm{Re}\{\hat{A}_{0}\left( \sigma
\right) \}\ .
\end{equation*}%
In fact, the second equality above, i.e.,
\begin{equation}
\lim_{J}\mathrm{Re}\{\sigma _{j}\left( A_{0}\right) \}=\lim_{J}\mathrm{Re}%
\{\sigma _{j}\left( A_{j}\right) \},  \label{eqsup}
\end{equation}%
results from the Banach-Steinhaus theorem: The family $(\hat{A}_{j})_{j\in J}
$, as defined by (\ref{sdfsdfkljsdlfkj}) for $A=A_{j}$, $j\in J$, is a
collection of continuous linear mappings from $\mathcal{X}^{\ast }$ to $%
\mathbb{K}$ and the set $\{\hat{A}_{j}\left( \sigma \right) :j\in J\}$ is
bounded for any $\sigma \in \mathcal{X}^{\ast }$, by convergence of $%
(A_{j})_{j\in J}$ to $A_{0}$. By \cite[Theorem 2.5]{Rudin}, $(\hat{A}%
_{j})_{j\in J}$ is equicontinuous and Equation (\ref{eqsup}) follows. As a
consequence, by (\ref{sdkljsdlkfjsdlkj}), there is an exposed point outside $%
\overline{\mathrm{co}}\left( \mathcal{E}_{0}(K)\right) $, which contradicts (%
\ref{contraduction}). Therefore, $\overline{\mathrm{co}}\left( \mathcal{E}%
_{0}(K)\right) =K$ and, by the Milman theorem \cite[Theorem 10.13]{BruPedra2}%
, $\mathcal{E}(K)\subseteq \overline{\mathcal{E}_{0}(K)}$.
\end{proof}

\begin{corollary}[Extension of the Straszewicz theorem - II]
\label{Solution selfbaby}\mbox{ }\newline
Let $\mathcal{X}$ be a Gateaux-differentiability $\mathbb{K}$-vector space ($%
\mathbb{K}=\mathbb{R},\mathbb{C}$). Then, for any $\mathcal{U}\in \mathcal{O}
$, $\mathcal{D}_{0,\mathcal{U}}=\mathcal{D}_{\mathcal{U}}$ with $\mathcal{D}%
_{\mathcal{U}}$ and $\mathcal{D}_{0,\mathcal{U}}$ being respectively defined
by (\ref{Zbis}) and (\ref{Zbisexposed}).
\end{corollary}

\begin{proof}
The assertion is an obvious consequence of Theorem \ref{stra theorem}.
\end{proof}

Obviously, Theorem \ref{stra theorem} and Corollary \ref{Solution selfbaby}
hold true for \emph{all} weak Asplund spaces. Well-known examples of such
spaces are separable Baire spaces: Recall that a topological space is a
Baire space if every non-empty open subset of this space is nonmeager. See,
e.g., \cite[Prerequisites B.9]{Schaefer}. Both completely metrizable spaces
and locally compact Hausdorff spaces are Baire spaces, by Baire's theorem
(also known as the category theorem) \cite[Theorem 2.2]{Rudin}. In
particular, Banach spaces are very specific Baire spaces. By \cite[Theorem
2.1]{Mazur-extended}, the Mazur theorem can be extended to all \emph{%
separable} topological vector spaces which are Baire spaces. We thus obtain
the following corollary:

\begin{corollary}[Extension of the Straszewicz theorem - Separable case]
\label{Solution selfbaby copy(1)}\mbox{ }\newline
If the topological vector space $\mathcal{X}$ is a separable Baire space
then Theorem \ref{stra theorem} and Corollary \ref{Solution selfbaby} hold
true.
\end{corollary}

\begin{proof}
Combine \cite[Theorem 2.1]{Mazur-extended} with Theorem \ref{stra theorem}
and Corollary \ref{Solution selfbaby}.
\end{proof}

Separability is not a necessary condition in Corollary \ref{Solution
selfbaby copy(1)}. In nonseparable Banach space theory, a pivotal role is
played by the so-called \emph{weakly compactly generated} Banach spaces $%
\mathcal{X}$, meaning that $\mathcal{X}$ is the closed linear span of a
weakly compact subset. See \cite[Definition 3.1]{Zizler} or \cite[Definition
2.41]{Phelps-Asplund}. Weakly compactly generated Banach spaces have been
proven to be weak Asplund spaces\footnote{%
In fact, Asplund proves that if $\mathcal{X}$ admits an equivalent norm
which has a strictly convex dual norm then $\mathcal{X}$ is a weak Asplund
space. See \cite[Corollary 2.39]{Phelps-Asplund}.} \cite[Theorem 2]%
{Asplund1968}. We thus obtain the following assertion:

\begin{corollary}[Extension of the Straszewicz theorem - Nonseparable case]
\label{Solution selfbaby copy(3)}\mbox{ }\newline
If $\mathcal{X}$ is a weakly compactly generated Banach space then Theorem %
\ref{stra theorem} and Corollary \ref{Solution selfbaby} hold true.
\end{corollary}

\begin{proof}
Combine \cite[Theorem 2]{Asplund1968} with Theorem \ref{stra theorem} and
Corollary \ref{Solution selfbaby}.
\end{proof}

\noindent For more details on weak Asplund spaces, see for instance \cite%
{Phelps-Asplund,Fabian}.

In 1979, Larman and Phelps in \cite{LarmanPhelps} raised the question
whether every Gateaux-differentiability space is a weak Asplund space. Known
examples \cite{Coban,Talagrand,Fabian} of individual convex continuous
functions that are Gateaux-differentiable on a dense, \emph{but non-residual}%
, subset of their domain suggest that a Gateaux-differentiability space is
not necessarily a weak Asplund space. A first answer to Larman and Phelps's
question\ has been given in 2006, in \cite{Moors} where a
Gateaux-differentiability space $\mathcal{X}$ that is \emph{not} a weak
Asplund space is constructed.

This shows that Theorem \ref{stra theorem} and Corollary \ref{Solution
selfbaby} are very general, probably \emph{optimum}, results on the weak$%
^{\ast }$-density of the set of exposed points in the extreme boundary of a
convex weak$^{\ast }$-compact set $K\in \mathbf{CK}(\mathcal{X}^{\ast })$.
If $\mathcal{X}$ is a Banach space, note that \cite[Theorem 6.2]%
{Phelps-Asplund} already tells us that a convex weak$^{\ast }$-compact set $%
K\in \mathbf{CK}(\mathcal{X}^{\ast })$ is the weak$^{\ast }$-closed convex
hull of its exposed points \emph{iff} $\mathcal{X}$ is a
Gateaux-differentiability space. This equivalence holds probably true for
more general (not necessarily Banach) topological vector spaces.

\subsection{$G_{\protect\delta }$ Subsets of Convex Weak$^{\ast }$-Compact
Sets\label{Gdelta}}

Recall that
\begin{equation*}
\mathbf{CU}_{\mathcal{U}}\left( \mathcal{X}^{\ast }\right) \doteq \mathbf{CF}%
\left( \mathcal{X}^{\ast }\right) \cap \mathbf{U}_{\mathcal{U}}\left(
\mathcal{X}^{\ast }\right) \doteq \left\{ U\in \mathbf{CF}\left( \mathcal{X}%
^{\ast }\right) :U\subseteq \mathcal{U}^{\circ }\right\} \ .
\end{equation*}%
See (\ref{polar})-(\ref{U0}) and (\ref{CU0}). Because of Corollary \ref%
{convexity corrolary copy(2)} (iii), it is a weak$^{\ast }$-Hausdorff-closed
space. If $\mathcal{X}$ is separable then it is even weak$^{\ast }$%
-Hausdorff-compact and completely metrizable, by Theorem \ref{Solution
selfbaby copy(4)+1}. From the Banach-Alaoglu theorem \cite[Theorem 3.15]%
{Rudin}, every $U\in \mathbf{CU}_{\mathcal{U}}(\mathcal{X}^{\ast })$ is a
convex weak$^{\ast }$-compact set and is thus the weak$^{\ast }$-closure of
the convex hull of the (nonempty) set $\mathcal{E}(U)$ of its extreme points
(cf. the Krein-Milman theorem \cite[Theorem 3.23]{Rudin}).

Having all this information in mind together with the Straszewicz theorem
(Corollary \ref{Solution selfbaby}), we are in a position to show that, for
general separable topological vector spaces $\mathcal{X}$ and each $\mathcal{%
U}\in \mathcal{O}$, the subset%
\begin{equation*}
\mathcal{D}_{\mathcal{U}}\doteq \left\{ U\in \mathbf{CU}_{\mathcal{U}}\left(
\mathcal{X}^{\ast }\right) :U=\overline{\mathcal{E}\left( U\right) }\right\}
\subseteq \mathbf{CU}_{\mathcal{U}}\left( \mathcal{X}^{\ast }\right) \ ,
\end{equation*}%
already defined in (\ref{Zbis}), is a $G_{\delta }$ subset of $\mathbf{CU}_{%
\mathcal{U}}(\mathcal{X}^{\ast })$ endowed with the weak$^{\ast }$-Hausdorff
hypertopology:

\begin{theorem}[$\mathcal{D}_{\mathcal{U}}$ as a $G_{\protect\delta }$ set]
\label{Solution selfbaby copy(5)+00}\mbox{ }\newline
Let $\mathcal{X}$ be a separable Gateaux-differentiability space (Definition %
\ref{Gateaux Differentiability space}). Then, for any $\mathcal{U}\in
\mathcal{O}$, $\mathcal{D}_{\mathcal{U}}$ is a $G_{\delta }$ subset of $%
\mathbf{CU}_{\mathcal{U}}(\mathcal{X}^{\ast })$.
\end{theorem}

\begin{proof}
Let $\mathcal{X}$ be a separable Gateaux-differentiability space. Assume
without loss of generality that $\mathcal{X}$ is a $\mathbb{C}$-vector
space. (The case $\mathbb{K}=\mathbb{R}$ is even slightly simpler.) For any $%
\mathcal{U}\in \mathcal{O}$, we use the metric $d$ defined by (\ref{metrics0}%
) and generating the weak$^{\ast }$ topology on the absolute polar $\mathcal{%
U}^{\circ }$ defined by (\ref{polar}). For any $\mathcal{U}\in \mathcal{O}$%
,\ recall that we denote by $B\left( \omega ,r\right) \subseteq \mathcal{U}%
^{\circ }$ the weak$^{\ast }$-open ball of radius $r\in \mathbb{R}^{+}$
centered at $\omega \in \mathcal{U}^{\circ }$, defined by (\ref{ball weak}).
Fix once and for all $\mathcal{U}\in \mathcal{O}$. Then, for any $m\in
\mathbb{N}$, let $\mathcal{F}_{\mathcal{U},m}$ be the set of all nonempty
convex weak$^{\ast }$-compact subsets $U\subseteq \mathcal{U}^{\circ }$ such
that $B\left( \omega ,1/m\right) \cap \mathcal{E}(U)=\emptyset $ for some $%
\omega \in U$, i.e.,
\begin{equation}
\mathcal{F}_{\mathcal{U},m}\doteq \left\{ U\in \mathbf{CU}_{\mathcal{U}%
}\left( \mathcal{X}^{\ast }\right) :\exists \omega \in U,\ B\left( \omega
,1/m\right) \cap \mathcal{E}\left( U\right) =\emptyset \right\} \subseteq
\mathbf{CU}_{\mathcal{U}}\left( \mathcal{X}^{\ast }\right) \ .  \label{Fm}
\end{equation}%
Recall that $\mathcal{E}(U)$ is the nonempty set of extreme points of $U$.
Now, observe that the complement of $\mathcal{D}_{\mathcal{U}}$ in $\mathbf{%
CU}_{\mathcal{U}}(\mathcal{X}^{\ast })$ equals%
\begin{equation}
\mathbf{CU}_{\mathcal{U}}\left( \mathcal{X}^{\ast }\right) \backslash
\mathcal{D}_{\mathcal{U}}=\bigcup_{m\in \mathbb{N}}\mathcal{F}_{\mathcal{U}%
,m}\ .  \label{union}
\end{equation}%
Therefore, $\mathcal{D}_{\mathcal{U}}$ is a $G_{\delta }$ subset of $\mathbf{%
CU}_{\mathcal{U}}(\mathcal{X}^{\ast })$ if $\mathcal{F}_{\mathcal{U},m}$ is
a weak$^{\ast }$-Hausdorff-closed set for any $m\in \mathbb{N}$.

By Theorem \ref{Solution selfbaby copy(4)+1}, the weak$^{\ast }$-Hausdorff
hypertopology of $\mathbf{CU}_{\mathcal{U}}(\mathcal{X}^{\ast })$ is
metrizable and $\mathbf{CU}_{\mathcal{U}}(\mathcal{X}^{\ast })$ is a weak$%
^{\ast }$-Hausdorff-closed set (see Corollary \ref{convexity corrolary
copy(2)} (iii)). So, fix $m\in \mathbb{N}$ and take any sequence $%
(U_{n})_{n\in \mathbb{N}}\subseteq \mathcal{F}_{\mathcal{U},m}$ converging
with respect to the weak$^{\ast }$-Hausdorff hypertopology to $U_{\infty
}\in $ $\mathbf{CU}_{\mathcal{U}}(\mathcal{X}^{\ast })$. For any $n\in
\mathbb{N}$, there is $\omega _{n}\in U_{n}$ such that $B\left( \omega
_{n},1/m\right) \cap \mathcal{E}(U_{n})=\emptyset $. By metrizability and
weak$^{\ast }$\ compactness of the absolute polar $\mathcal{U}^{\circ }$ and
Corollary \ref{Solution selfbaby copy(4)}, there is a subsequence $(\omega
_{n_{k}})_{k\in \mathbb{N}}$ converging to some $\omega _{\infty }\in
U_{\infty }$. Assume that, for some $\varepsilon \in (0,1/m)$, there is $%
\sigma _{\infty }\in \mathcal{E}(U_{\infty })$ such that%
\begin{equation*}
d\left( \omega _{\infty },\sigma _{\infty }\right) \leq \frac{1}{m}%
-\varepsilon \ .
\end{equation*}%
Recall meanwhile that $\mathbf{CU}_{\mathcal{U}}(\mathcal{X}^{\ast
})\subseteq \mathbf{CK}(\mathcal{X}^{\ast })$. Then, by weak$^{\ast }$%
-density of the set of exposed points in $\mathcal{E}(U_{\infty })$ (Theorem %
\ref{stra theorem}), we can assume without loss of generality that $\sigma
_{\infty }$ is an exposed point. In particular, there is $A\in \mathcal{X}$
such that
\begin{equation}
\max_{\sigma \in U_{\infty }}\mathrm{Re}\{\hat{A}(\sigma )\}=\hat{A}\left(
\sigma _{\infty }\right) \ ,  \label{itupeva1}
\end{equation}%
with $\sigma _{\infty }$ being the \emph{unique} maximizer in $U_{\infty }$
and where $\hat{A}$ is the mapping $\sigma \mapsto \sigma (A)$ from $%
\mathcal{X}^{\ast }$ to $\mathbb{C}$ defined by (\ref{sdfsdfkljsdlfkj}).
Consider now the sets%
\begin{equation*}
\mathcal{M}_{n}\doteq \left\{ \tilde{\sigma}\in U_{n}:\max_{\sigma \in U_{n}}%
\mathrm{Re}\{\hat{A}(\sigma )\}=\hat{A}\left( \tilde{\sigma}\right) \right\}
\ ,\qquad n\in \mathbb{N}\ .
\end{equation*}%
By linearity and weak$^{\ast }$-continuity of the function $\hat{A}$,
together with the weak$^{\ast }$-compactness of $U_{n}$, the set $\mathcal{M}%
_{n}$ is a convex weak$^{\ast }$-compact subset of $U_{n}$ for any $n\in
\mathbb{N}$. In fact, $\mathcal{M}_{n}$ is a (weak$^{\ast }$-closed)\ face%
\footnote{%
It means that if $\sigma \in \mathcal{M}_{n}$ is a finite convex combination
of elements $\sigma _{j}\in U_{n}$ then all $\sigma _{j}\in \mathcal{M}_{n}$.%
} of $U_{n}$ and thus, any extreme point of $\mathcal{M}_{n}$ belongs to $%
\mathcal{E}(U_{n})$. So, pick any extreme point $\sigma _{n}\in \mathcal{E}%
(U_{n})$ of $\mathcal{M}_{n}$ for each $n\in \mathbb{N}$. Since
\begin{eqnarray*}
\max_{\sigma \in U_{n}}\mathrm{Re}\{\hat{A}(\sigma )\}-\max_{\tilde{\sigma}%
\in U_{\infty }}\mathrm{Re}\{\hat{A}(\tilde{\sigma})\} &=&\max_{\sigma \in
U_{n}}\min_{\tilde{\sigma}\in U_{\infty }}\mathrm{Re}\{\hat{A}(\sigma -%
\tilde{\sigma})\}\leq \max_{\sigma \in U_{n}}\min_{\tilde{\sigma}\in
U_{\infty }}\left\vert \left( \sigma -\tilde{\sigma}\right) \left( A\right)
\right\vert \ , \\
\max_{\tilde{\sigma}\in U_{\infty }}\mathrm{Re}\{\hat{A}(\tilde{\sigma}%
)\}-\max_{\sigma \in U_{n}}\mathrm{Re}\{\hat{A}(\sigma )\} &=&\max_{\tilde{%
\sigma}\in U_{\infty }}\min_{\sigma \in U_{n}}\mathrm{Re}\{\hat{A}(\tilde{%
\sigma}-\sigma )\}\leq \max_{\tilde{\sigma}\in U_{\infty }}\min_{\sigma \in
U_{n}}\left\vert \left( \sigma -\tilde{\sigma}\right) \left( A\right)
\right\vert \ ,
\end{eqnarray*}%
we deduce from Definition \ref{hypertopology0} and the weak$^{\ast }$%
-Hausdorff convergence of $(U_{n})_{n\in \mathbb{N}}$ to $U_{\infty }$ that%
\begin{equation*}
\lim_{n\rightarrow \infty }\mathrm{Re}\{\hat{A}(\sigma
_{n})\}=\lim_{n\rightarrow \infty }\max_{\sigma \in U_{n}}\mathrm{Re}\{\hat{A%
}(\sigma )\}=\max_{\sigma \in U_{\infty }}\mathrm{Re}\{\hat{A}(\sigma )\}=%
\hat{A}\left( \sigma _{\infty }\right) \ .
\end{equation*}%
Therefore, keeping in mind the convergence of the subsequence $(\omega
_{n_{k}})_{k\in \mathbb{N}}$ towards $\omega _{\infty }\in U_{\infty }$,
there is a subsequence $(\sigma _{n_{k(l)}})_{l\in \mathbb{N}}$ of $(\sigma
_{n_{k}})_{k\in \mathbb{N}}$ (itself being a subsequence of $(\sigma
_{n})_{n\in \mathbb{N}}$) converging to $\sigma _{\infty }$, as it is the
\emph{unique} maximizer of (\ref{itupeva1}) and $\hat{A}$ is weak$^{\ast }$%
-continuous. Since, for any $l\in \mathbb{N}$,
\begin{multline*}
d(\sigma _{n_{k(l)}},\omega _{n_{k(l)}})\leq d(\sigma _{\infty },\omega
_{\infty })+d(\omega _{\infty },\omega _{n_{k(l)}})+d(\sigma
_{n_{k(l)}},\sigma _{\infty }) \\
\leq \frac{1}{m}-\varepsilon +d(\omega _{\infty },\omega
_{n_{k(l)}})+d(\sigma _{n_{k(l)}},\sigma _{\infty })
\end{multline*}%
with $\varepsilon \in (0,1/m)$ and $\sigma _{n}\in \mathcal{E}(U_{n})$ for $%
n\in \mathbb{N}$, we thus arrive at a contradiction. Therefore, $U_{\infty
}\in \mathcal{F}_{\mathcal{U},m}$. This means that $\mathcal{F}_{\mathcal{U}%
,m}$ is a weak$^{\ast }$-Hausdorff-closed set for any $m\in \mathbb{N}$ and
hence, the countable union (\ref{union}) is a $F_{\sigma }$ set with
complement being $\mathcal{D}_{\mathcal{U}}$. The assertion follows, as the
complement of an $F_{\sigma }$ set is a $G_{\delta }$ set.
\end{proof}

\begin{corollary}[$\mathcal{D}_{\mathcal{U}}$ as a $G_{\protect\delta }$ set]

\label{Solution selfbaby copy(6)}\mbox{ }\newline
If the topological vector space $\mathcal{X}$ is a separable Baire space
then, for any $\mathcal{U}\in \mathcal{O}$, $\mathcal{D}_{\mathcal{U}}$ is a
$G_{\delta }$ subset of $\mathbf{CU}_{\mathcal{U}}(\mathcal{X}^{\ast })$.
\end{corollary}

\begin{proof}
Combine \cite[Theorem 2.1]{Mazur-extended} with Theorem \ref{Solution
selfbaby copy(5)+00}.
\end{proof}

To conclude this section, recall now that $\mathbf{CU}(\mathcal{X}^{\ast })$
is the set of all nonempty, uniformly bounded in a $0$-neighborhood, weak$%
^{\ast }$-closed subsets of $\mathcal{X}^{\ast }$ defined by (\ref{CU}),
that is,
\begin{equation*}
\mathbf{CU}\left( \mathcal{X}^{\ast }\right) \doteq \bigcup\limits_{\mathcal{%
U}\in \mathcal{O}}\mathbf{CU}_{\mathcal{U}}\left( \mathcal{X}^{\ast }\right)
\ .
\end{equation*}%
By considering the special case of completely metrizable topological vector
spaces $\mathcal{X}$, which are separable Baire spaces, one can also prove
that the subset
\begin{equation*}
\mathcal{D}\doteq \bigcup\limits_{\mathcal{U}\in \mathcal{O}}\mathcal{D}_{%
\mathcal{U}}
\end{equation*}%
of all $U\in \mathbf{CU}\left( \mathcal{X}^{\ast }\right) $ with weak$^{\ast
}$-dense extreme boundary\ (see (\ref{Zbis})) is also a $G_{\delta }$ subset
of $\mathbf{CU}(\mathcal{X}^{\ast })$, endowed with the weak$^{\ast }$%
-Hausdorff hypertopology:

\begin{corollary}[$\mathcal{D}$ as a $G_{\protect\delta }$ set]
\label{Solution selfbaby copy(7)}\mbox{ }\newline
If $\mathcal{X}$ is a separable, completely metrizable topological vector
space then $\mathcal{D}$ is a $G_{\delta }$ subset of $\mathbf{CU}(\mathcal{X%
}^{\ast })$.
\end{corollary}

\begin{proof}
By assumption, there is a metric $d_{\mathcal{X}}$ generating the topology
on $\mathcal{X}$ and, since $\mathcal{U}_{1}\subseteq \mathcal{U}_{2}$
yields $\mathcal{U}_{2}^{\circ }\subseteq \mathcal{U}_{1}^{\circ }$ (see (%
\ref{polar})), we observe that%
\begin{equation}
\mathbf{CU}\left( \mathcal{X}^{\ast }\right) =\bigcup\limits_{D\in \mathbb{N}%
}\mathbf{CU}_{B\left( 0,D^{-1}\right) }\left( \mathcal{X}^{\ast }\right)
\qquad \text{and}\qquad \mathbf{CU}\left( \mathcal{X}^{\ast }\right)
\backslash \mathcal{D=}\bigcup_{D,m\in \mathbb{N}}\mathcal{F}_{B\left(
0,D^{-1}\right) ,m}\ ,  \label{sdfkljsdflkj}
\end{equation}%
where $\mathcal{F}_{\mathcal{U},m}$ is defined by (\ref{Fm}) for any $%
\mathcal{U}\in \mathcal{O}$ and
\begin{equation}
B\left( 0,R\right) \doteq \left\{ A\in \mathcal{X}:d_{\mathcal{X}}\left(
0,A\right) <R\right\} \in \mathcal{O}  \label{norm ball}
\end{equation}%
is the open ball of radius $R\in \mathbb{R}^{+}$ in $\mathcal{X}$. Any
completely metrizable topological vector space is a
Gateaux-differentiability space, by \cite[Theorem 2.1]{Mazur-extended}. As
shown in the proof of Theorem \ref{Solution selfbaby copy(5)+00}, $\mathcal{F%
}_{\mathcal{U},m}$ is thus a weak$^{\ast }$-Hausdorff-closed set for any $%
m\in \mathbb{N}$ and $\mathcal{U}\in \mathcal{O}$. By (\ref{sdfkljsdflkj}),
the assertion follows.
\end{proof}

\subsection{Generic Convex Weak$^{\ast }$-Compact Sets}

Outcomes of Sections \ref{Dense Subsets} and \ref{Gdelta} for separable,
infinite-dimensional topological vector spaces directly yield, as stated in
the next theorem, that the convex weak$^{\ast }$-compact sets with weak$%
^{\ast }$-dense extreme boundary are generic convex weak$^{\ast }$-compact
sets, answering the main question raised in the introduction of Section \ref%
{Generic section}. Recall that $\mathcal{D}_{\mathcal{U}}$ and $\mathcal{D}$
are the sets of convex weak$^{\ast }$-compact sets with weak$^{\ast }$-dense
extreme boundary, respectively defined by (\ref{Zbis}) and (\ref{Zbis-full}).

\begin{theorem}[Generic convex weak$^{\ast }$-compact sets for topological
spaces]
\label{genereic theorem}\mbox{ }\newline
\emph{(i)} Let $\mathcal{X}$ be a separable Gateaux-differentiability space
(Definition \ref{Gateaux Differentiability space}) and a $0$-neighborhood $%
\mathcal{U}\in \mathcal{O}$ satisfying Condition \ref{condition}. Then, $%
\mathcal{D}_{\mathcal{U}}$ is a weak$^{\ast }$-Hausdorff-dense $G_{\delta }$
subset of $\mathbf{CU}_{\mathcal{U}}(\mathcal{X}^{\ast })$.\newline
\emph{(ii)} Let $\mathcal{X}$ be a separable, completely metrizable space
and assume the existence of some $0$-neighborhood in $\mathcal{X}$
satisfying Condition \ref{condition}. Then $\mathcal{D}$ is a weak$^{\ast }$%
-Hausdorff-dense $G_{\delta }$ subset of $\mathbf{CU}(\mathcal{X}^{\ast })$.
\end{theorem}

\begin{proof}
In order to prove the first statement (i), combine Theorem \ref{Solution
selfbaby copy(5)+0} with Theorem \ref{Solution selfbaby copy(5)+00}.
Assertion (ii) is a direct consequence of Corollaries \ref{Solution selfbaby
copy(8)} and \ref{Solution selfbaby copy(7)}.
\end{proof}

As a matter of fact, the Hausdorff metric topology is very fine, as compared
to various standard hypertopologies (apart from the Vietoris\footnote{%
Vietoris and Hausdorff metric topologies are not comparable.}
hypertopology). Consequently, the weak$^{\ast }$-Hausdorff hypertopology can
be seen as a very fine, weak$^{\ast }$-type, topology on $\mathbf{CU}(%
\mathcal{X}^{\ast })$. It means that the density of the subset of all convex
weak$^{\ast }$ compact sets with weak$^{\ast }$-dense set of extreme points
stated in Theorem \ref{genereic theorem} is a very strong property. The
study of generic convex weak$^{\ast }$-compact sets can also be performed
within special weak$^{\ast }$-closed subsets of the dual space $\mathcal{X}%
^{\ast }$. An important example is given by so-called \emph{positive}
continuous functionals.

In order to present this example, we consider a topological $\mathbb{R}$%
-vector space $\mathcal{X}$. As is usual, $\mathcal{X}^{+}\subseteq \mathcal{%
X}$ is said to be a \emph{convex cone}\ in $\mathcal{X}$ if, for all $%
\lambda \in \mathbb{R}^{+}$ and $A_{1},A_{2}\in \mathcal{X}^{+}$, $\lambda
A_{1}\in \mathcal{X}^{+}$ (cf. cone) and $A_{1}+A_{2}\in \mathcal{X}^{+}$
(cf. convex). Such a $\mathcal{X}^{+}$ is of course a convex set in the
usual sense. Recall that any convex cone $\mathcal{X}^{+}\subseteq \mathcal{X%
}$ naturally defines a preorder relation\footnote{%
I.e., a binary relation that is reflexive and transitive, but not necessarly
antisymmetric.} $\preceq $ in $\mathcal{X}$ as follows: $A_{1}\preceq A_{2}$
iff $A_{2}-A_{1}\in \mathcal{X}^{+}$. This preoder is compatible\footnote{%
That is, if $A_{1}\preceq A_{2}$ then, for any $\alpha \in \mathbb{R}^{+}$
and $A\in \mathcal{X}$, $\alpha A_{1}\preceq \alpha A_{2}$ and $%
A_{1}+A\preceq A_{2}+A$.} with the \emph{real} vector space structure of $%
\mathcal{X}$. Given such a preorder, intervals in $\mathcal{X}$ are defined
to be the sets
\begin{equation}
\left[ A_{1},A_{2}\right] \doteq \left\{ A\in \mathcal{X}:A_{1}\preceq
A\preceq A_{2}\right\} \subseteq \mathcal{X}\ ,\qquad A_{1},A_{2}\in
\mathcal{X}\ .  \label{interval}
\end{equation}

A linear functional on $\mathcal{X}$ is defined to be \emph{positive} if its
values on $\mathcal{X}^{+}$ are non-negative. The set of all continuous
positive functionals is denoted by
\begin{equation*}
\mathcal{X}^{\ast ,+}\doteq \bigcap\limits_{A\in \mathcal{X}^{+}}\left\{
\sigma \in \mathcal{X}^{\ast }:\sigma \left( A\right) \geq 0\right\}
\subseteq \mathcal{X}^{\ast }\ .
\end{equation*}%
Being the intersection of weak$^{\ast }$-closed sets, $\mathcal{X}^{\ast ,+}$
is weak$^{\ast }$-closed. Therefore, we consider the hyperspaces defined,
for any $0$-neighborhood $\mathcal{U}\in \mathcal{O}$, by
\begin{equation}
\mathbf{CU}_{\mathcal{U}}\left( \mathcal{X}^{\ast ,+}\right) \doteq \left\{
U\in \mathbf{CU}_{\mathcal{U}}\left( \mathcal{X}^{\ast }\right) :U\subseteq
\mathcal{X}^{\ast ,+}\right\} \mathcal{\quad }\text{and}\mathcal{\quad }%
\mathbf{CU}\left( \mathcal{X}^{\ast ,+}\right) \doteq \bigcup\limits_{%
\mathcal{U}\in \mathcal{O}}\mathbf{CU}_{\mathcal{U}}\left( \mathcal{X}^{\ast
,+}\right) \ .  \label{definition positive}
\end{equation}%
Similar to $\mathbf{CU}_{\mathcal{U}}(\mathcal{X}^{\ast })$ (Corollary \ref%
{convexity corrolary copy(2)} (iii)), $\mathbf{CU}_{\mathcal{U}}(\mathcal{X}%
^{\ast ,+})$ is a weak$^{\ast }$-Hausdorff-closed space, at least under the
separability condition:

\begin{lemma}[Complete metrizability of hyperspaces for positive functionals]

\label{connected sub-hyperspace copy(1)}\mbox{ }\newline
Let $\mathcal{X}$ be a separable topological $\mathbb{R}$-vector space.
Then, for any $\mathcal{U}\in \mathcal{O}$, $\mathbf{CU}_{\mathcal{U}}(%
\mathcal{X}^{\ast ,+})$ is a weak$^{\ast }$-Hausdorff-closed subset of the
compact and completely metrizable hyperspace $\mathbf{CU}_{\mathcal{U}}(%
\mathcal{X}^{\ast })$.
\end{lemma}

\begin{proof}
Fix any $0$-neighborhood $\mathcal{U}\in \mathcal{O}$. By (\ref{definition
positive}) and Theorem \ref{Solution selfbaby copy(4)+1}, $\mathbf{CU}_{%
\mathcal{U}}(\mathcal{X}^{\ast ,+})$ belongs to the weak$^{\ast }$%
-Hausdorff-compact and completely metrizable hyperspace $\mathbf{CU}_{%
\mathcal{U}}(\mathcal{X}^{\ast })$. By Corollary \ref{Solution selfbaby
copy(4)} and because $\mathcal{X}^{\ast ,+}$ is a weak$^{\ast }$-closed set,
$\mathbf{CU}_{\mathcal{U}}(\mathcal{X}^{\ast ,+})$ is weak$^{\ast }$%
-Hausdorff-closed.
\end{proof}

Using (\ref{Zbis}), (\ref{Zbisexposed}) and (\ref{Zbis-full}), we naturally
define the subsets
\begin{equation*}
\mathcal{D}_{\mathcal{U}}^{+}\doteq \mathcal{D}_{\mathcal{U}}\cap \mathbf{CU}%
_{\mathcal{U}}\left( \mathcal{X}^{\ast ,+}\right) \mathcal{\qquad }\text{and}%
\mathcal{\qquad D}_{0,\mathcal{U}}^{+}\doteq \mathcal{D}_{0,\mathcal{U}}\cap
\mathbf{CU}_{\mathcal{U}}\left( \mathcal{X}^{\ast ,+}\right) \subseteq
\mathcal{D}_{\mathcal{U}}^{+}
\end{equation*}%
for any $\mathcal{U}\in \mathcal{O}$, as well as
\begin{equation}
\mathcal{D}^{+}\doteq \mathcal{D}\cap \mathbf{CU}\left( \mathcal{X}^{\ast
,+}\right) \mathcal{\qquad }\text{and}\mathcal{\qquad D}_{0}^{+}\doteq
\mathcal{D}_{0}\cap \mathbf{CU}\left( \mathcal{X}^{\ast ,+}\right) \subseteq
\mathcal{D}^{+}\ .  \label{D+}
\end{equation}%
Then, similar to Theorem \ref{Solution selfbaby copy(5)+0} and Corollary \ref%
{Solution selfbaby copy(8)}, we then get the following density property
within the set of positive continuous functionals:

\begin{theorem}[Weak$^{\ast }$-Hausdorf density of $\mathcal{D}_{0,\mathcal{U%
}}^{+}$]
\label{Solution selfbaby copy(9)}\mbox{ }\newline
Let $\mathcal{X}$ be a separable locally convex $\mathbb{R}$-vector space
with $\mathcal{X}^{+}\subseteq \mathcal{X}$ being a convex cone. Assume
that, for all $\mathcal{U}\in \mathcal{O}$, there is an interval $%
[A_{1},A_{2}]\subseteq \mathcal{X}$ and some $\mathcal{\tilde{U}}\in
\mathcal{O}$ such that $\mathcal{\tilde{U}}\subseteq \lbrack
A_{1},A_{2}]\subseteq \mathcal{U}$. Then,\emph{\ }for any $\mathcal{U}\in
\mathcal{O}$ satisfying Condition \ref{condition} and $U\in \mathbf{CU}_{%
\mathcal{U}}(\mathcal{X}^{\ast ,+})$, there is $\mathcal{\tilde{U}\in }$ $%
\mathcal{O}$ with $\mathcal{\tilde{U}}\subseteq \mathcal{U}$ and a sequence $%
(U_{n})_{n\in \mathbb{N}}\subseteq \mathcal{D}_{0,\mathcal{\tilde{U}}%
}^{+}\subseteq \mathcal{D}_{\mathcal{\tilde{U}}}^{+}$ converging to $U$ in
the weak$^{\ast }$-Hausdorff topology.
\end{theorem}

\begin{proof}
In order to show the assertion, it suffices to reproduce the proof of
Theorem \ref{Solution selfbaby copy(5)+0}, with the addition of one
essential ingredient: the decomposition of \emph{equicontinuous} linear
functionals into positive \emph{equicontinuous} components \cite%
{decomposition}, as originally proven by Grosberg and Krein \cite%
{decomposition1} for normed spaces and by Bonsall \cite{decomposition2} for
locally convex $\mathbb{R}$-vector spaces. For any $\mathcal{U}\in \mathcal{O%
}$, this means that there is $\mathcal{\tilde{U}\in }$ $\mathcal{O}$ with $%
\mathcal{\tilde{U}}\subseteq \mathcal{U}$ such that \emph{every} continuous
linear functional $\sigma \in \mathcal{U}^{\circ }$ can be decomposed as
\begin{equation}
\sigma =\rho _{1}-\rho _{2}\ ,\qquad \rho _{1},\rho _{2}\in \mathcal{X}%
^{\ast ,+}\cap \mathcal{\tilde{U}}^{\circ }\ .  \label{decoposition0}
\end{equation}%
In order to prove this assertion, all conditions of the theorem are
necessary, except Condition \ref{condition} and the separability of $%
\mathcal{X}$. Note in particular that the existence of arbitrarily small
neighborhoods of the origin which are intervals, as assumed here, is a
necessary and sufficient condition to obtain the decomposition of
equicontinuous linear functionals in any real locally convex vector space,
as shown in \cite{decomposition}.

Then, at \emph{Step 1} of the proof of Theorem \ref{Solution selfbaby
copy(5)+0}, because of (\ref{decoposition0}), there is%
\begin{equation*}
\sigma _{1}\in (\mathcal{\tilde{U}}^{\circ }\backslash \mathrm{span}\{\omega
_{1},\ldots ,\omega _{n_{\varepsilon }}\})\cap \mathcal{X}^{\ast ,+}\ .
\end{equation*}%
So, we proceed by taking $\sigma _{1}$ as a positive continuous functional
in $\mathcal{\tilde{U}}^{\circ }$ instead of a general functional of $%
\mathcal{U}^{\circ }$. One then iterates the arguments, as explained in the
proof of Theorem \ref{Solution selfbaby copy(5)+0}, always taking a positive
continuous functional $\sigma _{n}\in \mathcal{X}^{\ast ,+}\cap \mathcal{%
\tilde{U}}^{\circ }$ appearing in the positive decomposition (\ref%
{decoposition0}) ($\rho _{1}$ or $\rho _{2}$) of a continuous functional $%
\sigma \in \mathcal{U}^{\circ }$, as already explained. In doing so, we
ensure that the convex weak$^{\ast }$-compact set $U_{\infty }$ of Equation (%
\ref{equaion}) belongs to $\mathbf{CU}_{\mathcal{\tilde{U}}}(\mathcal{X}%
^{\ast ,+})$. Of course, the neighborhood $\mathcal{U}$ in Equations (\ref%
{toto0}), (\ref{sdfsdf}) and (\ref{toto}) in the proof of Theorem \ref%
{Solution selfbaby copy(5)+0} has to be replaced by $\mathcal{\tilde{U}}$.
\end{proof}

Observe that the last theorem is not a direct consequence of Theorem \ref%
{Solution selfbaby copy(5)+0} because the complement of $\mathcal{D}_{0,%
\mathcal{U}}^{+}$ is generally open and dense in $\mathcal{D}_{0,\mathcal{U}%
} $.

\begin{corollary}[Weak$^{\ast }$-Hausdorf density of $\mathcal{D}_{0}^{+}$]
\label{Solution selfbaby copy(10)}\mbox{ }\newline
Assume conditions of Theorem \ref{Solution selfbaby copy(9)} together with
the existence of one $0$-neighborhood satisfying Condition \ref{condition}.
Then, $\mathcal{D}_{0}^{+}\subseteq \mathcal{D}^{+}$ is a weak$^{\ast }$%
-Hausdorff dense subset of $\mathbf{CU}(\mathcal{X}^{\ast ,+})$.
\end{corollary}

\begin{proof}
The intersection of any $0$-neighborhood $\mathcal{U}\in \mathcal{O}$ with
the $0$-neighborhood satisfying Condition \ref{condition} leads to a new $0$%
-neighborhood satisfying Condition \ref{condition}. In particular, we can
assume without loss of generality that any element $U\in \mathbf{CU}(%
\mathcal{X}^{\ast ,+})$ belongs to some $\mathbf{CU}_{\mathcal{U}}(\mathcal{X%
}^{\ast ,+})$ with $\mathcal{U}\in \mathcal{O}$ satisfying Condition \ref%
{condition}. The assertion is then a direct consequence of Theorem \ref%
{Solution selfbaby copy(9)}.
\end{proof}

By Lemma \ref{connected sub-hyperspace copy(1)} combined with Theorem \ref%
{Solution selfbaby copy(4)+1} and Corollary \ref{Solution selfbaby copy(4)},
note that the proof of Theorem \ref{Solution selfbaby copy(5)+00} can be
straightforwardly adapted in order to prove that, for any $0$-neighborhood $%
\mathcal{U}\in \mathcal{O}$, $\mathcal{D}_{\mathcal{U}}^{+}$ is a $G_{\delta
}$ subset of $\mathbf{CU}_{\mathcal{U}}(\mathcal{X}^{\ast ,+})$, provided $%
\mathcal{X}$ is a separable Gateaux-differentiability space (Definition \ref%
{Gateaux Differentiability space}). If $\mathcal{X}$ is a separable,
completely metrizable topological vector space then $\mathcal{D}^{+}$ is a $%
G_{\delta }$ subset of $\mathbf{CU}(\mathcal{X}^{\ast ,+})$. Cf. proof of
Corollary \ref{Solution selfbaby copy(7)}. Using this together with
Corollary \ref{Solution selfbaby copy(10)}, one directly gets:

\begin{theorem}[Generic convex weak$^{\ast }$-compact sets - Positive
functionals]
\label{Solution selfbaby copy(12)}\mbox{ }\newline
Assume the conditions of Corollary \ref{Solution selfbaby copy(10)} with $%
\mathcal{X}$ being a completely metrizable topological vector space. Then, $%
\mathcal{D}^{+}$ is a weak$^{\ast }$-Hausdorff-dense $G_{\delta }$ subset of
$\mathbf{CU}(\mathcal{X}^{\ast ,+})$.
\end{theorem}

\begin{proof}
See proofs of Theorem \ref{Solution selfbaby copy(5)+00} and Corollary \ref%
{Solution selfbaby copy(7)} together with Theorem \ref{Solution selfbaby
copy(4)+1}, Corollary \ref{Solution selfbaby copy(4)}, Lemma \ref{connected
sub-hyperspace copy(1)} and Corollary \ref{Solution selfbaby copy(10)}. We
omit the details.
\end{proof}

Such results for sets of positive functionals is, e.g., useful in the
context of separable $C^{\ast }$-algebras, as explained in Section \ref%
{Application}, because closed balls centered at the origin in the (real)
Banach space of all self-adjoint elements of such algebras are intervals.

\subsection{Application to Separable Banach Spaces\label{Application}}

If $\mathcal{X}$ is a Banach space then recall that $\mathbf{B}(\mathcal{X}%
^{\ast })$ equals the set of all nonempty norm-bounded weak$^{\ast }$-closed
subsets of $\mathcal{X}^{\ast }$ and a set in $\mathcal{X}^{\ast }$ is
norm-bounded and weak$^{\ast }$-closed iff it is weak$^{\ast }$-compact.
See, e.g., \cite[Proposition 1.2.9]{Beer}. In particular, in this situation,
absolute polars can be replaced with norm-closed balls in $\mathcal{X}^{\ast
}$ and%
\begin{equation}
\mathbf{CU}\left( \mathcal{X}^{\ast }\right) =\mathbf{CK}\left( \mathcal{X}%
^{\ast }\right) =\mathbf{CB}\left( \mathcal{X}^{\ast }\right)
\label{eqaulity banach}
\end{equation}%
is a convex, path-connected, weak$^{\ast }$-Hausdorff-closed subset of $%
\mathbf{CF}(\mathcal{X}^{\ast })$, by Equation (\ref{Banach}) and Corollary %
\ref{convexity corrolary copy(2)} (ii). In this case, Theorem \ref{genereic
theorem} (ii) can be rephrased as follows:

\begin{theorem}[Generic convex weak$^{\ast }$-compact sets]
\label{theorem dense cool1}\mbox{ }\newline
Let $\mathcal{X}$ be an infinite-dimensional separable Banach space. Then,
the set $\mathcal{D}$ of all nonempty convex weak$^{\ast }$-compact sets $K$
with a weak$^{\ast }$-dense set $\mathcal{E}(K)$ of extreme points is a weak$%
^{\ast }$-Hausdorff-dense $G_{\delta }$ subset of the weak$^{\ast }$%
-Hausdorff-closed space $\mathbf{CK}(\mathcal{X}^{\ast })=\mathbf{CB}(%
\mathcal{X}^{\ast })$.
\end{theorem}

\begin{proof}
If $\mathcal{X}$ is an infinite-dimensional separable Banach space, then any
open ball of center $0$ in $\mathcal{X}$ satisfies Condition \ref{condition}%
. Therefore, the assertion follows from Theorem \ref{genereic theorem} (ii).
\end{proof}

\noindent As a consequence, $\mathcal{D}$ is generic in the hyperspace $%
\mathbf{CK}(\mathcal{X}^{\ast })$, that is, the complement of a meagre set,
i.e., a nowhere dense set. In other words, $\mathcal{D}$ is of second
category in $\mathbf{CK}(\mathcal{X}^{\ast })$.

In the Banach case, the weak$^{\ast }$-Hausdorff hypertopology on $\mathbf{CK%
}(\mathcal{X}^{\ast })$ is finer than the scalar topology \cite[Section 4.3]%
{Beer} restricted to weak$^{\ast }$-closed sets. The linear topology on the
set of nonempty closed convex subsets is the supremum of the scalar and
Wijsman topologies. Since the Wijsman topology \cite[Definition 2.1.1]{Beer}
requires a metric space, one has to use the norm on $\mathcal{X}^{\ast }$
and the linear topology is not comparable with the weak$^{\ast }$-Hausdorff
hypertopology. If one uses the metric (\ref{metrics0}) generated the weak$%
^{\ast }$ topology on balls of $\mathcal{X}^{\ast }$ for a separable Banach
space $\mathcal{X}$, then the Wijsman and linear topologies for norm-closed
balls of $\mathcal{X}^{\ast }$ are coarser than the weak$^{\ast }$-Hausdorff
hypertopology, by Theorem \ref{Solution selfbaby copy(4)+1}. In fact, as
already mentioned, the weak$^{\ast }$-Hausdorff hypertopology can be seen as
a very fine, weak$^{\ast }$-type, topology on $\mathbf{CK}(\mathcal{X}^{\ast
})$ and the density of the subset of all convex weak$^{\ast }$ compact sets
with weak$^{\ast }$-dense set of extreme points stated in Theorem \ref%
{theorem dense cool1} is a very strong property.

Meanwhile, in the Banach case, the weak$^{\ast }$-Hausdorff density property
within the set of positive continuous functionals, as stated in Theorems \ref%
{Solution selfbaby copy(9)} and \ref{Solution selfbaby copy(12)}, can be
strengthened. In order to present these outcomes, by fixing a convex cone $%
\mathcal{X}^{+}\subseteq \mathcal{X}$\ in a real Banach space $\mathcal{X}$
and using the usual norm $\Vert \cdot \Vert _{\mathcal{X}^{\ast }}$ on
continuous functionals of $\mathcal{X}^{\ast }$, we define the spaces%
\begin{eqnarray}
\mathbf{CK}\left( \mathcal{X}^{\ast ,+}\right) &\doteq &\left\{ K\in \mathbf{%
CK}(\mathcal{X}^{\ast }):K\subseteq \mathcal{X}^{\ast ,+}\right\}
\label{space 1} \\
\mathbf{CK}_{\leq 1}\left( \mathcal{X}^{\ast ,+}\right) &\doteq &\left\{
K\in \mathbf{CK}\left( \mathcal{X}^{\ast ,+}\right) :\forall \sigma \in K,\
\left\Vert \sigma \right\Vert _{\mathcal{X}^{\ast }}\leq 1\right\}
\label{space 2} \\
\mathbf{CK}_{1}\left( \mathcal{X}^{\ast ,+}\right) &\doteq &\left\{ K\in
\mathbf{CK}\left( \mathcal{X}^{\ast ,+}\right) :\forall \sigma \in K,\
\left\Vert \sigma \right\Vert _{\mathcal{X}^{\ast }}=1\right\}
\label{space 3}
\end{eqnarray}%
as well as the set%
\begin{equation}
\mathcal{D}^{+}\doteq \left\{ K\in \mathbf{CK}\left( \mathcal{X}^{\ast
,+}\right) :K=\overline{\mathcal{E}\left( K\right) }\right\}  \label{space 4}
\end{equation}%
of all nonempty convex weak$^{\ast }$-compact sets of $\mathcal{X}^{\ast ,+}$
with weak$^{\ast }$-dense extreme boundary. Compare with (\ref{D+}).
Additionally, define
\begin{equation*}
E_{\mathcal{X}}\doteq \left\{ \sigma \in \mathcal{X}^{\ast ,+}:\left\Vert
\sigma \right\Vert _{\mathcal{X}^{\ast }}=1\right\} ,
\end{equation*}%
the set of positive normalized functionals of $\mathcal{X}^{\ast }$. We
arrive at the following result:

\begin{theorem}[Generic convex weak$^{\ast }$-compact sets - Positive
functionals]
\label{Solution selfbaby copy(11)}\mbox{ }\newline
Let $\mathcal{X}$ be an infinite-dimensional separable real Banach space
with $\mathcal{X}^{+}\subseteq \mathcal{X}$ being a convex cone. Assume
that, for any open ball $B\left( 0,R\right) $ of radius $R\in \mathbb{R}^{+}$%
, centered at $0\in \mathcal{X}$, there is an interval $[A_{1},A_{2}]%
\subseteq \mathcal{X}$ (see (\ref{interval}))\ and some $r\in \mathbb{R}^{+}$
such that $B\left( 0,r\right) \subseteq \lbrack A_{1},A_{2}]\subseteq
B\left( 0,R\right) $. Then, one has:\newline
\emph{(i)} $\mathcal{D}^{+}$ is a weak$^{\ast }$-Hausdorff-dense $G_{\delta
} $ subset of the weak$^{\ast }$-Hausdorff-closed space $\mathbf{CK}(%
\mathcal{X}^{\ast ,+})$. \newline
\emph{(ii)} $\mathcal{D}^{+}\cap \mathbf{CK}_{\leq 1}(\mathcal{X}^{\ast ,+})$%
\emph{\ }is a weak$^{\ast }$-Hausdorff-dense $G_{\delta }$ subset of the weak%
$^{\ast }$-Hausdorff-compact and completely metrizable space $\mathbf{CK}%
_{\leq 1}(\mathcal{X}^{\ast ,+})$. \newline
\emph{(iii)} If $E_{\mathcal{X}}\in \mathbf{CK}_{1}(\mathcal{X}^{\ast ,+})$
then $\mathcal{D}^{+}\cap \mathbf{CK}_{1}(\mathcal{X}^{\ast ,+})$ is a weak$%
^{\ast }$-Hausdorff-dense $G_{\delta }$ subset of the weak$^{\ast }$%
-Hausdorff-compact and completely metrizable space $\mathbf{CK}_{1}(\mathcal{%
X}^{\ast ,+})$.
\end{theorem}

\begin{proof}
First of all, note that $\mathbf{CK}(\mathcal{X}^{\ast ,+})$ belongs to the
weak$^{\ast }$-Hausdorff-closed hyperspace $\mathbf{CK}(\mathcal{X}^{\ast })$%
, by (\ref{eqaulity banach}) and Corollary \ref{convexity corrolary copy(2)}
(ii). By Corollary \ref{Solution selfbaby copy(4)} and because $\mathcal{X}%
^{\ast ,+}$ is a weak$^{\ast }$-closed set, $\mathbf{CK}(\mathcal{X}^{\ast
,+})$ is also weak$^{\ast }$-Hausdorff-closed. By Lemma \ref{connected
sub-hyperspace copy(1)},
\begin{equation*}
\mathbf{CK}_{\leq 1}\left( \mathcal{X}^{\ast ,+}\right) =\mathbf{CU}%
_{B\left( 0,1\right) }\left( \mathcal{X}^{\ast }\right) \cap \mathbf{CK}(%
\mathcal{X}^{\ast ,+})
\end{equation*}%
is a weak$^{\ast }$-Hausdorff-compact and completely metrizable space. By
Corollary \ref{Solution selfbaby copy(4)} and the weak$^{\ast }$-closedness
of $E_{\mathcal{X}}\in \mathbf{CK}_{1}(\mathcal{X}^{\ast ,+})$, $\mathbf{CK}%
_{1}(\mathcal{X}^{\ast ,+})$ is also a weak$^{\ast }$-Hausdorff-closed
subspace of $\mathbf{CK}_{\leq 1}(\mathcal{X}^{\ast ,+})$, and is thus weak$%
^{\ast }$-Hausdorff-compact and completely metrizable. Now, we prove
Assertions (i)-(iii):\medskip

\noindent \underline{(i):} The first assertion is simply Theorem \ref%
{Solution selfbaby copy(12)} applied to the Banach case, keeping in mind
that all open balls of center $0$ in $\mathcal{X}$ satisfy Condition \ref%
{condition} when $\mathcal{X}$ is an infinite-dimensional separable Banach
space. \medskip

\noindent \underline{(ii):} In the proof of Theorem \ref{Solution selfbaby
copy(5)+00}, redefine the set $\mathcal{F}_{\mathcal{U},m}$ of Equation (\ref%
{Fm}) by
\begin{equation}
\mathcal{F}_{m}\doteq \left\{ K\in \mathbf{CK}_{\leq 1}\left( \mathcal{X}%
^{\ast ,+}\right) :\exists \omega \in K,\ B\left( \omega ,1/m\right) \cap
\mathcal{E}\left( K\right) =\emptyset \right\} \ .  \label{totot ala con}
\end{equation}%
Then, by following the arguments of the proof of Theorem \ref{Solution
selfbaby copy(5)+00}, it follows that $\mathcal{D}^{+}\cap \mathbf{CK}_{\leq
1}\left( \mathcal{X}^{\ast ,+}\right) $ is a $G_{\delta }$ subset of $%
\mathbf{CK}_{\leq 1}(\mathcal{X}^{\ast ,+})$, because $\mathbf{CK}_{\leq 1}(%
\mathcal{X}^{\ast ,+})$ is weak$^{\ast }$-Hausdorff-closed and any separable
real Banach space $\mathcal{X}$ is a Gateaux-differentiability space, by the
Mazur theorem \cite[Theorem 10.44]{BruPedra2}. It remains to prove the
asserted density. This follows from a straightforward adaptation of the
proof of Theorem \ref{Solution selfbaby copy(9)}: Use the same arguments
with the additional condition that $\Vert \sigma _{n}\Vert _{\mathcal{X}%
^{\ast }}\leq 1$. \medskip

\noindent \underline{(iii):} Recall that, by assumption, $E_{\mathcal{X}}$
is convex and weak$^{\ast }$-compact. We thus proceed in the same way as in
the proof of Assertion (ii). The proof of Theorem \ref{Solution selfbaby
copy(5)+00} can be adapted to this situation, by redefining the set $%
\mathcal{F}_{\mathcal{U},m}$ of Equation (\ref{Fm}): replace $\mathbf{CK}%
_{\leq 1}(\mathcal{X}^{\ast ,+})$ with $\mathbf{CK}_{1}(\mathcal{X}^{\ast
,+})$ in (\ref{totot ala con}). Hence, $\mathcal{D}^{+}\cap \mathbf{CK}_{1}(%
\mathcal{X}^{\ast ,+})$ is a $G_{\delta }$ subset of $\mathbf{CK}_{1}(%
\mathcal{X}^{\ast ,+})$. Like in\ Assertion (ii), in order to prove the
asserted density, we employ a straightforward adaptation of the proof of
Theorem \ref{Solution selfbaby copy(9)}: In addition to the condition $\Vert
\sigma _{n}\Vert _{\mathcal{X}^{\ast }}\leq 1$ used in the proof of
Assertion (ii), we also replace the definition of $\omega _{n_{\varepsilon
}+n}$ in (\ref{definition omegan}) by
\begin{equation*}
\omega _{n_{\varepsilon }+n}\doteq \left( 1-\lambda _{n}\Vert \sigma
_{n}\Vert _{\mathcal{X}^{\ast }}\right) \varpi _{n}+\lambda _{n}\sigma
_{n}\in E_{\mathcal{X}}
\end{equation*}%
at any step $n\in \mathbb{N}$.
\end{proof}

\noindent Note that Theorem \ref{Solution selfbaby copy(11)} (iii) does
\emph{not} directly follow from (i)-(ii) because the complement of $\mathbf{%
CK}_{1}(\mathcal{X}^{\ast ,+})$ in either $\mathbf{CK}(\mathcal{X}^{\ast
,+}) $ or $\mathbf{CK}_{\leq 1}(\mathcal{X}^{\ast ,+})$ is, in general, open
and dense in the weak$^{\ast }$-Hausdorff topology.

The assumptions of Theorem \ref{Solution selfbaby copy(11)} are particularly
relevant for $C^{\ast }$-algebra: Let $\mathcal{X}$ be a separable unital $%
C^{\ast }$-algebra, that is, a \emph{complex} Banach algebra having a unit $%
\mathfrak{1}\in \mathcal{X}$ for its product and endowed with an antilinear
involution $A\mapsto A^{\ast }$ such that%
\begin{equation*}
(AB)^{\ast }=B^{\ast }A^{\ast }\text{\qquad and\qquad }\left\Vert A^{\ast
}A\right\Vert _{\mathcal{X}}=\left\Vert A\right\Vert _{\mathcal{X}}^{2},%
\text{\qquad }A,B\in \mathcal{X}\ .
\end{equation*}%
The \emph{real} Banach (sub)space of all self-adjoint elements of $\mathcal{X%
}$ is denoted by
\begin{equation*}
\mathcal{X}^{\mathbb{R}}\doteq \left\{ A\in \mathcal{X}:A=A^{\ast }\right\}
\ .
\end{equation*}%
Recall that $\mathfrak{1}\in \mathcal{X}^{\mathbb{R}}$ and any element $A\in
\mathcal{X}$ can be decomposed as
\begin{equation*}
A=\mathrm{Re}\{A\}+i\mathrm{Im}\{A\}\qquad \text{with}\qquad \mathrm{Re}%
\{A\}=\frac{A+A^{\ast }}{2}\in \mathcal{X}^{\mathbb{R}},\qquad \mathrm{Im}%
\{A\}=\frac{A-A^{\ast }}{2i}\in \mathcal{X}^{\mathbb{R}}.
\end{equation*}%
Therefore, any continuous linear functional $\sigma \in \mathcal{X}^{\ast }$
is uniquely determined by its values on the real Banach space $\mathcal{X}^{%
\mathbb{R}}$ and any $\sigma \in (\mathcal{X}^{\mathbb{R}})^{\ast }$ can be
identified with an element of $\mathcal{X}^{\ast }$. As is well-known, such
functionals are the so-called hermitian elements\footnote{%
That is, functionals $\sigma \in \mathcal{X}^{\ast }$ satisfying $\sigma
\left( A^{\ast }\right) =\overline{\sigma \left( A\right) }$ for any $A\in
\mathcal{X}$.} of $\mathcal{X}^{\ast }$, which forms the \emph{real} vector
space%
\begin{equation}
(\mathcal{X}^{\ast })^{\mathbb{R}}\doteq \bigcap\limits_{A\in \mathcal{X}%
}\left\{ \sigma \in \mathcal{X}^{\ast }:\sigma \left( A^{\ast }\right) =%
\overline{\sigma \left( A\right) }\right\} \subseteq \mathcal{X}^{\ast }\ .
\label{hermitian functional}
\end{equation}%
The mapping from $(\mathcal{X}^{\mathbb{R}})^{\ast }$ to $(\mathcal{X}^{\ast
})^{\mathbb{R}}$ is denoted by $\daleth $ and, for any $\sigma \in (\mathcal{%
X}^{\mathbb{R}})^{\ast }$,%
\begin{equation}
\daleth \left( \sigma \right) \left( A\right) \doteq \sigma \left( \mathrm{Re%
}\{A\}\right) +i\sigma \left( \mathrm{Im}\{A\}\right) \ ,\qquad A\in
\mathcal{X}\ .  \label{tamel}
\end{equation}%
This mapping is a (real) linear homeomorphism, $(\mathcal{X}^{\mathbb{R}%
})^{\ast }$ and $(\mathcal{X}^{\ast })^{\mathbb{R}}$ being endowed with the
weak$^{\ast }$ topology\footnote{%
The same holds true with respect to the norm topology for continuous linear
functionals.}. It naturally induces a mapping, again denoted by $\daleth $,
from $\mathbf{CK}((\mathcal{X}^{\mathbb{R}})^{\ast })$ to%
\begin{equation*}
\mathbf{CK}\left( (\mathcal{X}^{\ast })^{\mathbb{R}}\right) \doteq \left\{
K\in \mathbf{CK}(\mathcal{X}^{\ast }):K\subseteq (\mathcal{X}^{\ast })^{%
\mathbb{R}}\right\} \ ,
\end{equation*}%
which is again a homeomorphism:

\begin{lemma}[Homeomorphism between complex and real structures in $C^{\ast
} $-algebras]
\label{lemma pas si trivial}\mbox{ }\newline
Let $\mathcal{X}$ be a separable unital $C^{\ast }$-algebra. $\daleth $\ is
a homeomorphism from $\mathbf{CK}((\mathcal{X}^{\mathbb{R}})^{\ast })$ to $%
\mathbf{CK}\left( (\mathcal{X}^{\ast })^{\mathbb{R}}\right) $, both spaces
being endowed with the weak$^{\ast }$-Hausdorff topology.
\end{lemma}

\begin{proof}
Any element of $(\mathcal{X}^{\ast })^{\mathbb{R}}$, i.e., any hermitian
continuous linear functional on $\mathcal{X}$, can be pushed forward,
through the restriction mapping, to an element of $(\mathcal{X}^{\mathbb{R}%
})^{\ast }$, i.e., a continuous linear functional on $\mathcal{X}^{\mathbb{R}%
}$. This obviously yields a weak$^{\ast }$-Hausdorff continuous bijective
mapping $\beth $ from $\mathbf{CK}\left( (\mathcal{X}^{\ast })^{\mathbb{R}%
}\right) $ to $\mathbf{CK}((\mathcal{X}^{\mathbb{R}})^{\ast })$. Observe
also that the union of any weak$^{\ast }$-Hausdorff convergent net $%
(K_{j})_{j\in J}\subseteq \mathbf{CK}((\mathcal{X}^{\ast })^{\mathbb{R}})$
is norm-bounded. To prove this, one uses an argument by contradiction and
the uniform boundedness principle (see, e.g., \cite[Theorems 2.4 and 2.5]%
{Rudin}). Therefore, we restrict without loss of generality our study on
\begin{equation*}
\mathbf{CK}_{B\left( 0,R\right) }\left( (\mathcal{X}^{\ast })^{\mathbb{R}%
}\right) \doteq \mathbf{CU}_{B\left( 0,R\right) }\left( \mathcal{X}^{\ast
}\right) \cap \mathbf{CK}((\mathcal{X}^{\ast })^{\mathbb{R}})
\end{equation*}%
at fixed $R\in \mathbb{R}^{+}$. Here, $B\left( 0,R\right) \subseteq \mathcal{%
X}$ is the open ball of radius $R\in \mathbb{R}^{+}$ centered at $0$, as
defined by (\ref{norm ball}). Since $\mathcal{X}$ is separable, we infer
from Theorem \ref{Solution selfbaby copy(4)+1} that $\mathbf{CU}_{B\left(
0,R\right) }\left( \mathcal{X}^{\ast }\right) $ is weak$^{\ast }$%
-Hausdorff-compact which, combined with Corollary \ref{Solution selfbaby
copy(4)}, in turn implies that $\mathbf{CK}_{B\left( 0,R\right) }\left( (%
\mathcal{X}^{\ast })^{\mathbb{R}}\right) $ is weak$^{\ast }$%
-Hausdorff-closed and, hence, weak$^{\ast }$-Hausdorff-compact. $\beth $ is
thus a weak$^{\ast }$-Hausdorff continuous bijective mapping from the weak$%
^{\ast }$-Hausdorff-compact space $\mathbf{CK}_{B\left( 0,R\right) }\left( (%
\mathcal{X}^{\ast })^{\mathbb{R}}\right) $ to $\mathbf{CK}_{B\left(
0,R\right) }((\mathcal{X}^{\mathbb{R}})^{\ast })$. Consequently, the inverse
of $\beth $ is also weak$^{\ast }$-continuous. This inverse is nothing else
than the mapping $\daleth $, which is thus a homeomorphism.
\end{proof}

\noindent The fact that $\daleth $\ defined on $(\mathcal{X}^{\mathbb{R}%
})^{\ast }$ by (\ref{tamel}) is a homeomorphism is obvious, but as a mapping
on $\mathbf{CK}((\mathcal{X}^{\mathbb{R}})^{\ast })$, this property,
asserted in Lemma \ref{lemma pas si trivial}, is a priori not clear at all.
It would be obvious if one used the hypertopology (Definition \ref{weak
hypertopology}) induced by the family of pseudometrics $\tilde{d}_{H}^{(A)}$
defined, for all $A\in \mathcal{X}$ and $F,\tilde{F}\in \mathbf{F}\left(
\mathcal{X}^{\ast }\right) $, by%
\begin{equation*}
\tilde{d}_{H}^{(A)}(F,\tilde{F})\doteq \max \left\{ \sup_{\sigma \in F}\inf_{%
\tilde{\sigma}\in \tilde{F}}\left\vert \mathrm{Re}\left\{ \left( \sigma -%
\tilde{\sigma}\right) \left( A\right) \right\} \right\vert ,\sup_{\tilde{%
\sigma}\in \tilde{F}}\inf_{\sigma \in F}\left\vert \mathrm{Re}\left\{ \left(
\sigma -\tilde{\sigma}\right) \left( A\right) \right\} \right\vert \right\}
\in \mathbb{R}^{+}\cup \left\{ \infty \right\} \ .
\end{equation*}%
Compare with Definition \ref{hypertopology0}: This hypertopology is \emph{%
coarser} than the weak$^{\ast }$-Hausdorff hypertopology, because $%
\left\vert \mathrm{Re}\left\{ z\right\} \right\vert \leq \left\vert
z\right\vert $ for any $z\in \mathbb{C}$. Note that the inequality $%
\left\vert z\right\vert \leq \left\vert \mathrm{Re}\left\{ z\right\}
\right\vert +\left\vert \mathrm{Im}\left\{ z\right\} \right\vert $ for $z\in
\mathbb{C}$ is not really useful to prove Lemma \ref{lemma pas si trivial}
because of the combination of suprema and infima in Definition \ref%
{hypertopology0}. To show Lemma \ref{lemma pas si trivial} we use the
separability of the $C^{\ast }$-algebra and abstract compactness arguments,
like in the proof of Theorem \ref{Solution selfbaby copy(4)+1}.

Meanwhile, it is also well-known \cite[Proposition 1.6.1]{Dixmier} that
\begin{equation*}
\mathcal{X}^{+}\equiv \left( \mathcal{X}^{\mathbb{R}}\right) ^{+}\doteq
\left\{ AA^{\ast }:A\in \mathcal{X}\right\} \subseteq \mathcal{X}^{\mathbb{R}%
}\subseteq \mathcal{X}
\end{equation*}%
is a closed convex cone, which is additionally pointed\footnote{%
That is, the only element $A\in (\mathcal{X}^{\mathbb{R}})^{+}$ such that $%
-A,A\in (\mathcal{X}^{\mathbb{R}})^{+}$ is $A=0$.}. In other words, the cone
$\mathcal{X}^{+}$ is the set of positive elements of the $C^{\ast }$-algebra
$\mathcal{X}$ \cite[Definition 1.6.5]{Dixmier}. It yields a partial order%
\footnote{%
I.e., a binary relation that is reflexive, transitive and antisymmetric.} $%
\preceq $ in $\mathcal{X}$ (and not only a preorder), and in $\mathcal{X}^{%
\mathbb{R}}$, which is compatible with the real vector space structure of $%
\mathcal{X}$. The weak$^{\ast }$-closed set of all positive (continuous)
linear functionals is denoted by
\begin{equation*}
\mathcal{X}^{\ast ,+}\doteq \bigcap\limits_{A\in \mathcal{X}^{+}}\left\{
\sigma \in \mathcal{X}^{\ast }:\sigma \left( A\right) \geq 0\right\}
\subseteq (\mathcal{X}^{\ast })^{\mathbb{R}}\ .
\end{equation*}%
It belongs to the set of hermitian elements of $\mathcal{X}^{\ast }$, by
\cite[Proposition 2.1.5]{Dixmier}. Using this set, we define by (\ref{space
1})-(\ref{space 4}), respectively, the spaces $\mathbf{CK}(\mathcal{X}^{\ast
,+})$, $\mathbf{CK}_{\leq 1}(\mathcal{X}^{\ast ,+})$, $\mathbf{CK}_{1}(%
\mathcal{X}^{\ast ,+})$ as well as $\mathcal{D}^{+}$. Then, one obtains an
analogous version of Theorem \ref{Solution selfbaby copy(11)} for separable
unital $C^{\ast }$-algebras:

\begin{corollary}[Generic convex weak$^{\ast }$-compact sets - $C^{\ast }$%
-algebras]
\label{Solution selfbaby copy(13)}\mbox{ }\newline
Let $\mathcal{X}$ be a infinite-dimensional separable unital $C^{\ast }$%
-algebra. Then, Assertions (i)-(iii) of Theorem \ref{Solution selfbaby
copy(11)} hold true.
\end{corollary}

\begin{proof}
First, since $\mathcal{X}$ is a (unital) $C^{\ast }$-algebra,
\begin{equation*}
\left\Vert A\right\Vert _{\mathcal{X}}=\max_{\sigma \in \mathcal{X}^{\ast
,+}:\ \left\Vert \sigma \right\Vert _{\mathcal{X}^{\ast }}=1}\left\vert
\sigma \left( A\right) \right\vert \ ,\qquad A\in \mathcal{X}^{\mathbb{R}}\ ,
\end{equation*}%
and, by the functional calculus in $C^{\ast }$-algebra \cite[Section 1.5]%
{Dixmier} combined with \cite[Proposition 1.6.1]{Dixmier}, for all $A\in
\mathcal{X}^{\mathbb{R}}$ and $R\in \mathbb{R}^{+}$,%
\begin{equation*}
0\leq \Vert A\Vert _{\mathcal{X}}\leq R\qquad \text{iff}\qquad -R\mathfrak{1}%
\preceq A\preceq R\mathfrak{1}\ .
\end{equation*}%
In particular, the closure of any open ball $B(0,R)$ of radius $R\in \mathbb{%
R}^{+}$ centered at $0\in \mathcal{X}$ in $\mathcal{X}^{\mathbb{R}}$ is the
interval $[-R\mathfrak{1},R\mathfrak{1}]\subseteq \mathcal{X}^{\mathbb{R}}$
(see (\ref{interval})). By \cite[Proposition 2.1.9]{Dixmier}, note
additionally that
\begin{equation}
E_{\mathcal{X}}\doteq \left\{ \sigma \in \mathcal{X}^{\ast ,+}:\left\Vert
\sigma \right\Vert _{\mathcal{X}^{\ast }}=1\right\} =\left\{ \sigma \in
\mathcal{X}^{\ast ,+}:\sigma \left( \mathfrak{1}\right) =1\right\} \in
\mathbf{CK}_{1}\left( \mathcal{X}^{\ast ,+}\right) \ .  \label{E}
\end{equation}%
We can thus apply Theorem \ref{Solution selfbaby copy(11)} for the real
Banach space $\mathcal{X}^{\mathbb{R}}$ to get Assertions (i)-(iii) with
\begin{equation*}
\left( \mathcal{X}^{\mathbb{R}}\right) ^{\ast ,+}=\daleth ^{-1}\left(
\mathcal{X}^{\ast ,+}\right)
\end{equation*}%
being the set of all positive (continuous) linear functionals of $(\mathcal{X%
}^{\mathbb{R}})^{\ast }$. By Lemma \ref{lemma pas si trivial}, the mapping $%
\daleth $ from $(\mathcal{X}^{\mathbb{R}})^{\ast }$ to $(\mathcal{X}^{\ast
})^{\mathbb{R}}\subseteq \mathcal{X}^{\ast }$ defined by (\ref{tamel})
yields a weak$^{\ast }$-Hausdorff homeomorphism between $\mathbf{CK}((%
\mathcal{X}^{\mathbb{R}})^{\ast })$ and $\mathbf{CK}((\mathcal{X}^{\ast })^{%
\mathbb{R}})$. Therefore, Assertions (i)-(iii) also hold true when we
replace $(\mathcal{X}^{\mathbb{R}})^{\ast ,+}$ with $\mathcal{X}^{\ast
,+}\subseteq (\mathcal{X}^{\ast })^{\mathbb{R}}$.
\end{proof}

If $\mathcal{X}$ is a separable unital $C^{\ast }$-algebra then the set $E_{%
\mathcal{X}}$ of all positive and normalized linear functionals of $\mathcal{%
X}^{\ast }$ (see (\ref{E})) is the so-called \emph{state space} in the
algebraic formulation of quantum mechanics. It is weak$^{\ast }$-closed and
thus a weak$^{\ast }$-compact subset of the unit ball of $\mathcal{X}^{\ast
} $, by the Banach-Alaoglu theorem \cite[Theorem 3.15]{Rudin}. Since $E_{%
\mathcal{X}}$ is also convex, $E_{\mathcal{X}}\in \mathbf{CK}_{1}(\mathcal{X}%
^{\ast ,+})$. For antiliminal\footnote{%
A $C^{\ast }$-algebra $\mathcal{X}$ is antiliminal if the zero ideal is its
only liminal closed two-sided ideal. A $C^{\ast }$-algebra $\mathcal{X}$ is
called liminal if, for every irreducible representation $\pi $ of $\mathcal{X%
}$\ and each $A\in \mathcal{X}$, $\pi (A)$ is compact.}\emph{\ }and simple%
\footnote{%
A $C^{\ast }$-algebra $\mathcal{X}$ is simple if the only closed two-sided
ideals of $\mathcal{X}$ are the trivial sets $\{0\}$ and $\mathcal{X}$.} $%
C^{\ast }$-algebras, $E_{\mathcal{X}}\in \mathcal{D}^{+}\cap \mathbf{CK}_{1}(%
\mathcal{X}^{\ast ,+})$, by \cite[Lemma 11.2.4]{Dixmier}. Important examples
of such $C^{\ast }$-algebras are the (even subalgebra of the) CAR $C^{\ast }$%
-algebras for (non-relativistic) fermions on the lattice. Quantum-spin
systems, i.e., infinite tensor products of copies of some elementary finite
dimensional matrix algebra, referring to a spin variable, are also important
examples.\ They are, for instance, widely used in quantum information theory
as well as in condensed matter physics. In all these physical situations,
the corresponding (non-commutative) $C^{\ast }$-algebra $\mathcal{X}$ is
separable and $E_{\mathcal{X}}$ is thus a metrizable weak$^{\ast }$-compact
convex set. It is \emph{not} a simplex \cite[Example 4.2.6]%
{BrattelliRobinsonI}, but
\begin{equation}
E_{\mathcal{X}}=\overline{\bigcup\limits_{n\in \mathbb{N}}\mathfrak{P}_{n}}
\label{property}
\end{equation}%
is the weak$^{\ast }$-closure of the union of a strictly increasing sequence
$(\mathfrak{P}_{n})_{n\in \mathbb{N}}\subseteq \mathcal{D}^{+}\cap \mathbf{CK%
}_{1}(\mathcal{X}^{\ast ,+})$ of Poulsen \cite{Poulsen} simplices. Equation (%
\ref{property}) is a consequence of well-known results (see, e.g., \cite%
{Israel,BruPedra2}) and we give its complete proof in \cite{BruPedra4}. In
other words, by Proposition \ref{Solution selfbaby copy(5)+1}, $E_{\mathcal{X%
}}$ is the \emph{weak}$^{\ast }$\emph{-Hausdorff limit} of the increasing
sequence $(\mathfrak{P}_{n})_{n\in \mathbb{N}}$ within $\mathcal{D}^{+}\cap
\mathbf{CK}_{1}(\mathcal{X}^{\ast ,+})$. Together with Corollary \ref%
{Solution selfbaby copy(13)}, this demonstrates the amazing structural
richness of the state space $E_{\mathcal{X}}$ for infinite-dimensional
quantum systems \bigskip

\noindent \textit{Acknowledgments:} This work is supported by CNPq
(308337/2017-4), FAPESP (2017/22340-9), as well as by the Basque Government
through the grant IT641-13 and the BERC 2018-2021 program, and by the
Spanish Ministry of Science, Innovation and Universities: BCAM Severo Ochoa
accreditation SEV-2017-0718, MTM2017-82160-C2-2-P. We finally thank the
reviewer for his/her thorough work and interest in the improvement of the
paper.

\noindent \textbf{Jean-Bernard Bru} \newline
Departamento de Matem\'{a}ticas\newline
Facultad de Ciencia y Tecnolog\'{\i}a\newline
Universidad del Pa\'{\i}s Vasco\newline
Apartado 644, 48080 Bilbao \medskip \newline
BCAM - Basque Center for Applied Mathematics\newline
Mazarredo, 14. \newline
48009 Bilbao\medskip \newline
IKERBASQUE, Basque Foundation for Science\newline
48011, Bilbao\bigskip

\noindent \textbf{Walter de Siqueira Pedra} \newline
Departamento de F\'{\i}sica Matem\'{a}tica\newline
Instituto de F\'{\i}sica,\newline
Universidade de S\~{a}o Paulo\newline
Rua do Mat\~{a}o 1371\newline
CEP 05508-090 S\~{a}o Paulo, SP Brasil


\begin{thebibliography}{99}
\bibitem{Beer} G. Beer, \textit{Topologies on closed and closed convex sets}%
, Mathematics and Its Applications 268, Springer Science+Business Media
Dordrecht, 1993.

\bibitem{Lucchetti} R. Lucchetti, \textit{Convexity and Well-Posed Problems}%
, CMS Books in Mathematics, Springer-Verlag New York, 2006.

\bibitem{Bankston} P. Bankston, Clopen Sets in Hyperspaces, Proceedings of
the AMS \textbf{54} (1976) 298-302.

\bibitem{Universal Algebra} S.N. Burris and H.P. Sankappanavar, \textit{A
Course in Universal Algebra}, Graduate Texts in Mathematics, Springer; 1
edition, 1981. Free online edition, revised in 2012.

\bibitem{topology} J.L. Kelley, \textit{General Topology}. Graduate Texts in
Mathematics, vol. 27, Springer-Verlag New York, 1975.

\bibitem{topology-painleve} K. Kuratowski, \textit{Topology, Volume I, new
edition, revised and augmented}, Translated from French by J. Jaworowski,
Academic Press New York and London, 1966.\

\bibitem{Dixmier} J. Dixmier, $C^{\ast }$\textit{-Algebras}, North-Holland
Publishing Company, Amsterdam - New York - Oxford, 1977.

\bibitem{Poulsen} E.T. Poulsen, A simplex with dense extreme points. Ann.
Inst. Fourier (Grenoble) \textbf{11} (1961) 83-87.

\bibitem{BruPedra4} J.-B. Bru and W. de Siqueira Pedra, Classical Dynamics
Generated by Long-Range Interactions for Lattice Fermions and Quantum Spins,
J. Math. Anal. Appl. \textbf{493}(1) (2021) 124434.

\bibitem{Simon} B. Simon, Operators with singular continuous spectrum: I.
General operators, Ann. Math. \textbf{141}(2) (1995) 131-145.

\bibitem{cesar1} S.L. Carvalho and C.R. de Oliveira, Correlation dimension
Wonderland theorems, J. Math. Phys. \textbf{57} (2016) 063501-1--19.

\bibitem{cesar2} S.L. Carvalho and C.R. de Oliveira, Generic Quasilocalized
and QuasiBallistic Discrete Schr\"{o}dinger Operators, Proceedings of the
AMS \textbf{144}(1) (2016) 129-141.

\bibitem{cesar3} S.L. Carvalho and C.R. de Oliveira, A characterization of
singular packing subspaces with an application to limit-periodic operators,
Forum Mathematicum \textbf{29} (2017) 31-40.

\bibitem{cesar4} S.L. Carvalho and C.R. de Oliveira, Generic zero-Hausdorff and  one-packing spectral measures, J. Math. Phys. \textbf{62} (2021) 013502.

\bibitem{Klee} V. Klee, Some new results on smoothness and rotundity in
normed linear spaces, Mathematische Annalen \textbf{139} (1959) 51-63.

\bibitem{FonfLindenstrauss} V.P. Fonf and J. Lindenstrauss, Some results on
infinite-dimensional convexity, Israel J. Math. \textbf{108} (1998) 13-32.

\bibitem{infinite dim convexity} V.P. Fonf, J. Lindenstrauss, R.R. Phelps,
Infinite dimensional convexity, Chapter 15, \textit{Handbook of the Geometry
of Banach Spaces}, Vol. 1, Edited by W.B. Johnson, J. Lindenstrauss, Issue
C, 2001.

\bibitem{LarmanPhelps} D. G. Larman and R. R. Phelps, Gateaux
differentiability of convex functions on Banach spaces, J. London Math. Soc.
\textbf{s2-20}(1) (1979) 115-127.

\bibitem{Phelps-Asplund} R.R. Phelps, \textit{Convex Functions, Monotone
Operators and Differentiability}, 2nd Edition, Lecture Notes in Mathematics
1364, Springer-Verlag Berlin Heidelberg 1989, 1993.

\bibitem{Straszewicz} S. Straszewicz: \"{U}ber exponierte Punkte
abgeschlossener Punktmengen. Fund. Math. \textbf{24} (1935) 139-143.

\bibitem{Klee-exposed} V. L. Klee, Extremal structure of convex sets. II,
Math. Zeitschr. \textbf{69} (1958) 90-104.

\bibitem{Bair} J. Bair, Extension du th\'{e}or\`{e}me de Straszewicz, Bull.
Soc. Roy. Se. Li\`{e}ge, \textbf{45} (1976) 166-168.

\bibitem{BairI} J. Bair and R. Fourneau, \textit{Etude G\'{e}om\'{e}trique
des Espaces Vectoriels II, Poly\`{e}dres et Polytopes Convexes, }%
Springer-Verlag Berlin Heidelberg New York 1980.

\bibitem{decomposition} J. D. Weston, The Decomposition of a continuous
linear functional into non-negative components, Math. Scand. \textbf{5}
(1957), 54-56.

\bibitem{decomposition1} J. Grosberg and M. Krein, Sur la d\'{e}composition
des fonctionnelles en composantes positives, C. R. (Doklady) Acad. Sci. URSS
\textbf{25} (1939) 723-726.

\bibitem{decomposition2} F. F. Bonsall, The decomposition of continuous
linear functionals into non-negative components, Proc. Durham Phil. Soc.
\textbf{13(A)} (1957) 6-11.

\bibitem{Rudin} W. Rudin, \textit{Functional Analysis}. McGraw-Hill Science,
1991.

\bibitem{Khaleelulla} S.M. Khaleelulla, \textit{Counterexamples in
Topological Vector Spaces}, Lecture Notes in Mathematics 936 (Eds: A. Dold
and B. Eckmann), Springer-Verlag Berlin Heidelberg New York, 1982.

\bibitem{BruPedra2} J.-B. Bru and W. de Siqueira Pedra, Non-cooperative
Equilibria of Fermi Systems With Long Range Interactions. Memoirs of the AMS
\textbf{224} (2013), no. 1052 (167 pages).

\bibitem{Aliprantis} C.D. Aliprantis and K.C. Border, \textit{Infinite
Dimensional Analysis, A Hitchhiker's Guide}, Third Edition, Springer-Verlag
Berlin Heidelberg 1999, 2006.

\bibitem{locally finite dim} B.R. Wenner, Finite--Dimensional Properties of
Infinite-dimensional spaces, Pacific Journal of Mathematics, \textbf{42}(1)
(1972) 267-276.

\bibitem{Asplund} E. Asplund, A $k$-extreme point is the limit of $k$%
-exposed points, Israel Journal of Mathematics, \textbf{1}(3) (1963) 161-162.

\bibitem{Mazur-extended} B. Sharp, The Differentiablility of Convex
Functions on Topological Linear Spaces, Bull. Austral. Math. Soc. \textbf{42}
(1990) 201-213

\bibitem{Mazur} S. Mazur, \"{U}ber konvexe Mengen in linearen normierten R%
\"{a}umen, Studia. Math. \textbf{4} (1933) 70-84

\bibitem{Schaefer} H.H. Schaefer with M.P. Wolff, \textit{Topological Vector
Spaces}, Second Edition, Springer Science+Business Media New York, 1999.

\bibitem{Zizler} V. Zizler, Nonseparable Banach spaces, Chapter 41, \textit{%
Handbook of the geometry of Banach spaces}, Vol. 2, Edited by W.B. Johnson,
J. Lindenstrauss, North-Holland, Amsterdam, 2003, pp. 1743--1816.

\bibitem{Asplund1968} E. Asplund, Fr\'{e}chet differentiability of convex
functions, Acta Math. \textbf{121} (1968) 31-47.

\bibitem{Fabian} M. Fabian, \textit{G\^{a}teaux Differentiability of Convex
Functions and Topology: Weak Asplund Spaces}, Wiley (1997).

\bibitem{Talagrand} M. Talagrand, Deux exemples de fonctions convexes, C. R.
Acad. Sc. Paris \textbf{288} (1979) 461-464.

\bibitem{Coban} M. Coban and P. S. Kenderov, Dense Gateaux differentiability
of the supnorm in C(T) and the topological properties of T, C. R. Acad.
Bulgare Sci. \textbf{38} (1985) 1603-1604.

\bibitem{Moors} W.B. Moors and S. Somasundaram, A G\^{a}teaux
Differentiability Space That Is Not Weak Asplund, Proceedings of the AMS
\textbf{134}(9) (2006) 2745-2754.

\bibitem{BrattelliRobinsonI} O. Bratteli and D.W. Robinson, \textit{Operator
Algebras and Quantum Statistical Mechanics, Vol. I, 2nd ed.} New York:
Springer-Verlag, 1987.

\bibitem{Israel} R.B. Israel, \textit{Convexity in the theory of lattice
gases}. Princeton: Princeton Series in Physics, Princeton Univ. Press, 1979.
\end{thebibliography}
\end{document}